\documentclass[reqno]{amsart}
\usepackage{amsmath,amssymb,stmaryrd}
\usepackage{amsrefs}
\usepackage{enumerate}
\usepackage{amsfonts}

\newtheorem{step}{Step}

\usepackage[colorlinks=true, pdfborder={0 0 0}]{hyperref}
\usepackage{cleveref}
\usepackage{upref}

\newtheorem{theorem}{Theorem}[section]
\newtheorem{rem}{Remark}

\newtheorem{proposition}{Proposition}
\newtheorem{coro}{Corollary}
\newtheorem{lemma}{Lemma}

\def \rr {\mathbb{R}}
\def \nn {\mathbb{N}}
\def \rn {\mathbb{R}^n}
\def \R {\mathbb{R}}
\def \eps {\epsilon}

\def \crit {2^\star}
\def \crits {2^\star(s)}
\def \ap {\alpha_+(\gamma)}
\def \am {\alpha_-(\gamma)}
\def \Ono { \Omega \setminus \{ 0 \} } 
\def \huno {H_{1,0}^2(\Omega)}
\def \Bn {\mathbb{B}^n}
\def \Omegabar {\overline{\Omega}}
\def \I {\mathcal{I}}
\def \rnp {\rn\setminus\{0\}}

\def \ctheta {C^1(\overline{\Omega}, |x|^{-\theta})}

\def \ue {u_\epsilon}
\def \ve {v_\epsilon}
\def \we {w_\eps}
\def \ze {z_\eps}
\def \re {r_\eps}
\def \he {h_\eps}
\def \pe {p_\eps}
\def \ye {y_\eps}
\def \elle {\ell_\eps}
\def \tv {\tilde{v}}
\def \tu {\tilde{u}}
\def \tve {\tilde{v}_\eps}
\def \tye {\tilde{y}_\eps}
\def \tue {\tilde{u}_\eps}
\def \tuei {\tilde{u}_{i,\eps}}
\def \tueun {\tilde{u}_{1,\eps}}
\def \tui {\tilde{u}_i}
\def \bue {\overline{u}_\eps}
\def \mei {\mu_{i,\eps}}
\def \kei {k_{i,\eps}}
\def \meun {\mu_{1,\eps}}
\def \keun {k_{1,\eps}}
\def \xeun {x_{1,\eps}}
\def \meN {\mu_{N,\eps}}
\def \keN {k_{N,\eps}}

\def \ds {\displaystyle}
\def \beq {\begin{eqnarray*}}
\def \eeq {\end{eqnarray*}}
\def \beqn {\begin{eqnarray}}
\def \eeqn {\end{eqnarray}}
\def \bequa {\begin{equation}}
\def \eequa {\end{equation}}

\title[Hardy-Schr\"odinger operator: Compactness and Multiplicity]{The Hardy--Schr\"odinger Operator on the Poincar\'e Ball: Compactness and Multiplicity}

\author{Nassif Ghoussoub}
\address{Nassif Ghoussoub: Department of Mathematics, University of British Columbia, Vancouver, V6T 1Z2 Canada}
\email{nassif@math.ubc.ca}
\author{Saikat Mazumdar}
\address{Saikat Mazumdar:
Department of Mathematics,
Indian Institute of Technology Bombay,
Mumbai 400076, India.}
\email{saikat@math.iitb.ac.in, saikat.mazumdar@iitb.ac.in}
\author{Fr\'ed\'eric Robert}
\address{Fr\'ed\'eric Robert, Institut \'Elie Cartan, Universit\'e de Lorraine, BP 70239, F-54506 Vand{\oe}uvre-l\`es-Nancy, France}
\email{frederic.robert@univ-lorraine.fr}

\date{March 31st 2021}
\thanks{This work was initiated when the second-named author held a postdoctoral position at the University of British Columbia under the supervision of the first-named author, that was partially supported by the Natural Sciences and Engineering Research Council of Canada (NSERC)}


\addtocontents{toc}{\protect\setcounter{tocdepth}{1}}

\begin{document}
\maketitle

\begin{abstract} 
Let $\Omega$
be  a compact smooth domain containing zero in the Poincar\'e ball model of the Hyperbolic space $\Bn$ ($n \geq 3$) and let $-\Delta_{\Bn}$ be the Laplace-Beltrami operator on $\Bn$,  associated with the metric 
$g_{\Bn}=  \frac{4}{(1-|x|^{2})^2}g_{_{\hbox{Eucl}}}$. We consider issues of non-existence, existence, and multiplicity of variational solutions 
for the borderline Dirichlet problem,
 \begin{eqnarray}\label{uno}
\left\{ \begin{array}{lll}
-\Delta_{\Bn}u-\gamma{V_2}u -\lambda u&=V_{\crits}|u|^{\crits-2}u &\hbox{ in }\Omega\\
\hfill u &=0   & \hbox{ on }  \partial \Omega,
\end{array} \right.
\end{eqnarray}
where $0\leq  \gamma \leq  \frac{(n-2)^2}{4}$, $0< s <2$, ${\crits}:=\frac{2(n-s)}{n-2}$ is the corresponding critical Sobolev exponent, $V_{2}$ (resp., $V_{\crits}$) is a Hardy-type potential (resp., Hardy-Sobolev weight) that is invariant under hyperbolic scaling and which behaves like $\frac{1}{r^{2}}$ (resp., $\frac{1}{r^{s}}$) at the origin. The bulk of this paper is a sharp blow-up analysis that we perform on approximate solutions of (\ref{uno}) with bounded but arbitrary high energies. When these approximate solutions are positive, our analysis leads to improvements of results in \cite{hypHS:min} regarding positive ground state solutions for \eqref{main.eq}, as we show that they exist whenever $n \geq 4$, $0 \leq \gamma \leq  \frac{(n-2)^2}{4}-1$ and $ \lambda > 0$.  The latter result also holds true for $n\geq 3$ and $\gamma >  \frac{(n-2)^2}{4}-1$ provided the domain has a positive ``hyperbolic mass". On the other hand, the same analysis yields that if $\gamma >  \frac{(n-2)^2}{4}-1$ and the mass is non vanishing, then there is a surprising stability of regimes where no variational positive solution exists. As for higher energy solutions to (\ref{uno}), we show that there are infinitely many of them provided $n\geq 5$, $0\leq \gamma<\frac{(n-2)^2}{4}-4$ and $ \lambda  > \frac{n-2}{n-4} \left(\frac{n(n-4)}{4}-\gamma \right)$.


\end{abstract}

\tableofcontents 

\section{Introduction}
Consider the  Poincar\'e ball model of the Hyperbolic space $\Bn$, $n \geq 3$, which is the Euclidean unit ball  $B_{1}(0):= \{ x \in \R^{n}: |x|<1\}$  endowed with the metric 
$g_{\Bn}=  \left( \frac{2}{1-|x|^{2}} \right)^{2}g_{_{\hbox{Eucl}}}$. Let
\begin{align}\label{hypgreen}
f(r):=\frac{(1-r^2)^{n-2}}{r^{n-1}} \ \hbox{and}\ {G}(r):=\int \limits_{r}^{1}f(t)dt, 
\end{align}
where $r=\sqrt{\sum \limits_{i=1}^{n} x_{i}^{2}}$  denotes  the Euclidean distance from a point $x$ to the origin.\\
The function $\frac{1}{n \omega_{n-1}}G(r)$ is then a fundamental solution of the Hyperbolic Laplacian $\Delta_{\Bn}u=   div_{\Bn}(\nabla_{\Bn}u)$. Moreover, if we consider the hyperbolic scaling of a given function $u: \Bn\to \R$, defined for $\lambda >0$ by
$${u}_\lambda(r)=\lambda^{-\frac12}u\left(G^{-1}(\lambda{G}(r))\right),$$ 
then for any radially symmetric $u\in{H}_r^1(\Bn)$ and $p\geq1$, one has the following invariance property:
\begin{align}
\int \limits_{\Bn} \left| \nabla_{\Bn} u_\lambda \right|^pdv_{g_{\Bn}}= \int \limits_{\Bn} |\nabla_{\Bn} u|^pdv_{g_{\Bn}} \ \hbox{ and } \ \int \limits_{\Bn}V_p |{u_\lambda}|^p dv_{g_{\Bn}}=\int \limits_{\Bn}V_p |{u}|^{p} dv_{g_{\Bn}},
\end{align}
where
\begin{align}\label{Definition-V_p}
V_p(r):= \frac{f(r)^2(1-r^2)^2}{ 4(n-2)^2 G(r)^{\frac{p+2}{2}}}.
\end{align}
The weights $V_{p}$ have the following  asymptotic behaviors: for $n\geq3$ and $p>1$, 
\begin{align}
&V_p(r)=\frac{c_0(n,p)}{r^{n(1-p/2^*)}}(1+o(1))\quad &\text{ as }{r\to0} \notag\\
&V_p(r)=\frac{c_1(n,p)}{(1-r)^{(n-1)(p-2)/2}}(1+o(1))\quad &\text{ as }{r\to1}.
\end{align}
In particular for $n \geq 3$, the weight $V_{2}(r)= \frac{1}{4(n-2)^{2}}\left( \dfrac{f(r)(1-r^{2})}{G(r)}\right)^{2} \sim_{ r \to0}\dfrac{1}{4r^{2}}$, while at $r=1$ it has a finite positive value. In other words, $V_{2}$ is qualitatively similar to the Euclidean Hardy potential, which led Sandeep--Tintarev  to establish  the following Hyperbolic Hardy inequality on $\Bn$ (Theorem 3.4  of  \cite{Sandeep-Tintarev}): 
\begin{align}
\frac{(n-2)^2}{4}\int \limits_{\Bn}V_2|{u}|^2~ dv_{g_{\Bn}} \leq \int\limits_{\Bn}|\nabla_{\Bn}{u}|^2 ~dv_{g_{\Bn}} \quad \hbox{for any } u\in{H}^1(\Bn).
\end{align}
They also show the following Hyperbolic Sobolev inequality: for some constant $C>0$.
\begin{align}
C\left(~\int \limits_{\Bn}V_{\crit} |{u}|^{\crit}~ dv_{g_{\Bn}} \right)^{2/\crit}\leq \int \limits_{\Bn}|\nabla_{\Bn}{u}|^2~ dv_{g_{\Bn}} \quad \hbox{for any $u\in{H}^1(\Bn)$.}
\end{align}
By interpolating between these two inequalities, one then easily obtains for $0\leq s\leq 2$,   the following Hyperbolic Hardy-Sobolev inequality \cite{hypHS:min}: \\
If $\gamma <\dfrac{(n-2)^2}{4}$, then there exists a constant $C>0$ such that for any $u\in{H}^1(\Bn)$,
\begin{align}\label{H-HS}
C\left(~\int \limits_{\Bn} V_{2^*(s)} |{u}|^{2^*(s)}~dv_{g_{\Bn}} \right)^{2/2^*(s)}\leq \int \limits_{\Bn} |\nabla_{\Bn}{u}|^2~dv_{g_{\Bn}}-\gamma\int \limits_{\Bn}V_2|{u}|^2~ dv_{g_{\Bn}},
\end{align}
\noindent where $\crits:= \frac{2(n-s)}{(n-2)}$. Note that  $V_{\crits}$ behaves like $\frac{1}{r^{s}}$ at the origin, making (\ref{H-HS}) the exact analogue of the Euclidean Hardy-Sobolev inequality:
$$C\left(~\int \limits_{\Omega}  \frac{|{u}|^{2^*(s)}}{|x|^{s}}~dx \right)^{2/2^*(s)}\leq \int \limits_{\Omega} |\nabla {u}|^2~dx -\gamma\int \limits_{\Omega} \frac{u^2}{|x|^{2}}~ dx \quad\hbox{for any $u\in C^{\infty}_{c}(\Omega)$}.$$
\medskip

Let $\Omega_{\Bn}$ be a compact smooth subdomain of  $\Bn$, $n\geq 3$,  such that $ 0 \in \Omega_{\Bn}$, but $\overline{\Omega}_{\Bn}$ does not touch the boundary of $\Bn$. We write $\Omega_{\Bn} \Subset \Bn$. In this paper, we are interested in the question of existence and multiplicity of solutions to the following Dirichlet  boundary value problem:
 \begin{eqnarray}\label{main.eq}
\left\{ \begin{array}{lll}
-\Delta_{\Bn}u-\gamma{V_2}u-\lambda u &=V_{\crits}|u|^{\crits-2}u &\hbox{ in }\Omega_{\Bn} \\
\hfill u &=0   & \hbox{ on }  \partial \Omega_{\Bn}.
\end{array} \right.
\end{eqnarray}
It is clear that \eqref{main.eq} is the Euler-Lagrange equation for the following energy functional on $H^1_0(\Omega_{\Bn})$:
\begin{align*}
\mathcal{J}_{\gamma, s,\lambda }(u):=\dfrac{ \int \limits_{\Omega_{\Bn}} \left( \left|\nabla_{\Bn} u \right|^2 -\gamma V_{2} u^2 -\lambda u^{2}\right)dv_{g_{\Bn}}  }{\left( ~ \int \limits_{\Omega_{\Bn}} |u|^{2^*(s)} V_{\crits}dv_{g_{\Bn}} \right)^{2/\crits }},
\end{align*}
where  $H_{0}^1(\Omega_{\Bn})$ is the completion of $C_{c}^{\infty}(\Omega_{\Bn})$ with respect to the norm  given by $\| u\|= \sqrt{\int \limits_{\Omega_{\Bn}} |\nabla_{\Bn} u|^{2} dv_{g_{\Bn}}}$.
\medskip
 The hyperbolic mass  $m_{\gamma,\lambda}(\Omega_{\Bn})$ of a domain $\Omega_{\Bn}$, associated to the operator $-\Delta_{\Bn}-\gamma{V_2}-\lambda  $ was defined in \cite{hypHS:min} and is recalled below.

\noindent The  existence of a positive ground state solution for (\ref{main.eq}) has already been addressed in \cite{hypHS:min}. 
While our main objective in this paper is to study higher energy solutions, our methods improve in a couple of ways their results about ground state solutions, in particular in the critical dimensions $3$ and $4$. We get the following.

\begin{theorem}\label{ground:sc:improv}
Let $\Omega_{\Bn} \Subset \Bn$ be a smooth compact domain containing $0$ and let $0< s < 2$.  Assume that $0\leq\gamma< \frac{(n-2)^2}{4}$Let $\lambda\in\rr$ be such
that the operator $-\Delta_{\Bn}-\gamma{V_2}-\lambda  $ is coercive. Then  
\begin{enumerate}
\item For $n \geq 4$ and $ \gamma \leq  \frac{(n-2)^2}{4}-1$, equation  \eqref{main.eq} has a positive ground state solution whenever $ \lambda > 0$. 
\item For $n \geq 3$ and $\gamma >  \frac{(n-2)^2}{4}-1$, equation \eqref{main.eq} has a positive ground state solution whenever the mass $m_{\gamma,\lambda}(\Omega_{\Bn})>0$. 
\ \end{enumerate}
\end{theorem}
\medskip\noindent  We note that statement (1) of the above theorem was established in \cite{hypHS:min} in the case where $n \geq 5$, $0 \leq \gamma \leq  \frac{(n-2)^2}{4}-1$ and $ \lambda >  \dfrac{n-2}{n-4} \left(\dfrac{n(n-4)}{4}-\gamma \right) $. As for statement (2), it was proved in \cite{hypHS:min} under additional assumptions on $\lambda$, at least in dimensions $3$ and $4$. 

\medskip\noindent The above theorem is actually a consequence of the following more general compactness result for positive approximate solutions.
\begin{theorem}\label{th:cpct:hyp}
Let $\Omega_{\Bn} \Subset \Bn$ be a smooth compact domain containing $0$ and let $0< s < 2$ and $\lambda\in\rr$ such that the operator $-\Delta_{\Bn} -\gamma V_2 -\lambda$ is coercive. Let $(\pe)_{\eps>0}$ be such that $0\leq p_\eps<\crits-2$ and $\lim \limits_{\eps\to 0}\pe=0$ and consider a sequence of functions $(\ue)_{\eps>0}$ that  is  uniformly bounded in $H^1_0(\Omega_{\Bn})$ such that for each $\epsilon >0$, $u_\epsilon$ is a solution to the equation:
\begin{align*}
\left\{ \begin{array}{clll}
-\Delta_{\Bn} \ue-\gamma V_2 \ue-\lambda\ue&=&V_{\crits} \ue^{\crits-1-p_\eps}\ \ &\text{in } \Omega_{\Bn},\\
\hfill \ue &>&0 &\text{in } \Omega_{\Bn},\\
\hfill \ue&=&0 &\text{on }\partial \Omega_{\Bn}.
\end{array} \right.
\end{align*}
Assume that $0\leq \gamma\leq \frac{(n-2)^2}{4}$. Assuming one of the following conditions
\begin{itemize}
\item $n\geq 4$, $  0\leq \gamma \leq \frac{(n-2)^{2}}{4}-1 $ and $\lambda>0$; 
\item  $\gamma > \frac{(n-2)^{2}}{4}-1 $ and $m_{\gamma,\lambda}(\Omega)_{\Bn}>0$. 
\end{itemize}
Then  the sequence  $(\ue)_{\eps>0}$ is  pre-compact in the space $H^1_0(\Omega_{\Bn})$. 
\end{theorem}
We now recall the notion of hyperbolic mass of a domain $\Omega$ associated with a coercive operator $-\Delta_{\Bn} -\gamma V_2 -\lambda$. This notion was introduced in the Euclidean case in \cite{gr1} and was extended to the hyperbolic setting in \cite{hypHS:min}. 

\begin{theorem}[The hyperbolic interior mass]\label{def:mass:hyp}
Let $\Omega_{\Bn} \Subset \Bn$ ($   n \geq 3$) be a smooth compact domain containing $0$. Let  $\gamma<\dfrac{(n-2)^2}{4}$ and let $\lambda\in\rr$  be such that the operator $-\Delta_{\Bn}-\gamma{V_2}-\lambda  $ is coercive on $\Omega_{\Bn}$. Then, there exists $H\in C^\infty(\overline{\Omega_{\Bn}}\setminus\{0\})$ such that 
\begin{align}\label{interior.mass.hyp} 
\qquad\qquad \left\{\begin{array}{cl}
-\Delta_{\Bn}H-\gamma{V_2}H-\lambda H   =0 &\hbox{ in } \Omega_{\Bn}\setminus\{0\}\\
 H>0&\hbox{ in } \Omega_{\Bn}\setminus\{0\}\\
 H=0&\hbox{ on }\partial\Omega_{\Bn}.
\end{array}\right.
\end{align}
These solutions are unique up to a positive multiplicative constant, and there exists $ c_{+}>0$ such that 
$$\displaystyle H(x)\simeq_{x\to 0}\dfrac{c_{+}}{|x|^{\alpha_+(\gamma)}}.$$ 
Moreover, when $\gamma>\max \left\{ \dfrac{n(n-4)}{4},0\right\} $, then there exists $c_{-}\in\rr$ such that
\begin{equation*}
H(x)=\dfrac{c_{+}}{|x|^{\ap}}+\dfrac{c_{-}}{|x|^{\am}}+o\left(\dfrac{1}{|x|^{\am}}\right)\hbox{ as }x\to 0
\end{equation*}
where $\ap,\am$ stand for
\begin{equation}\label{def:alpha}
\alpha_{\pm}(\gamma):=\dfrac{n-2}{2}\pm\sqrt{\dfrac{(n-2)^2}{4}-\gamma}.
\end{equation}
We define the interior mass as $m_{\gamma,\lambda}(\Omega_{\Bn}):=\dfrac{c_{+}}{c_{-}}$, which is independent of the choice of $H$.
\end{theorem}
This existence theorem will follow directly from Theorem \ref{def:mass} below.

\medskip\noindent In this paper we are  concerned with showing the multiplicity of higher energy solutions. Here is our main result.

\begin{theorem}\label{main.hyp:sc} Let $\Omega_{\Bn} \Subset \Bn$, $   n \geq 5$,  be a smooth compact domain containing $0$ and let $0< s < 2$.  Let $\lambda\in\rr$ be such
that the operator $-\Delta_{\Bn}-\gamma{V_2}-\lambda  $ is coercive. 

If $0\leq \gamma <  \dfrac{(n-2)^{2}}{4}-4$, then equation \eqref{main.eq} has an infinite number of solutions corresponding to higher energy critical levels for $\mathcal{J}_{\gamma, s}$ whenever $ \lambda  \geq  \dfrac{n-2}{n-4} \left(\dfrac{n(n-4)}{4}-\gamma \right) $.  
\end{theorem}

\noindent Note that -unlike the case of ground state solutions- we do not have a multiplicity result neither in {\it low dimensions}, nor when $\frac{(n-2)^2}{4}-4 \leq \gamma  \leq \frac{(n-2)^2}{4}-1$, even under the assumption of positivity of the {\it mass} of the compact subdomain. Note also the condition $ \lambda  \geq  \dfrac{n-2}{n-4} \left(\dfrac{n(n-4)}{4}-\gamma \right) $, which is more restrictive than $\lambda >0$.  
\medskip

In order to prove Theorem \ref{main.hyp:sc}, we shall proceed as in  \cite{hypHS:min} and use a conformal transformation  
$$\displaystyle g_{\Bn}= \varphi^{\frac{4}{n-2}}\hbox{Eucl}\quad \hbox{where \quad $\displaystyle \varphi= \left( \frac{2}{1-r^{2}} \right)^{\frac{n-2}{2}},$} $$
 to reduce equation  (\ref{main.eq}) to 
 a  Dirichlet  boundary value problem on Euclidean space. We shall denote by $\Omega$ the subdomain $\Omega_{\Bn}$ considered as a subset of $\R^n$. 
 The following sharpens a few estimates obtained in \cite{hypHS:min}. These will be needed when dealing with the critical (low dimensional) cases.

\begin{lemma}\label{Lemma:relation between Hyp and Euc solution on bounded domain}
 $u \in H_{0}^1(\Omega_{\Bn})$ satisfies  \eqref{main.eq} if and only if  ${ v:= \varphi u} \in H^1_0(\Omega_{\Omega})$ satisfies
\begin{eqnarray} \label{Problem on bounded domain:Euclidean}
\left\{ \begin{array}{lll}
-\Delta v  - \left( \dfrac{\gamma}{|x|^{2}} + h_{\gamma,\lambda}(x)\right)  v &=b(x) \dfrac{v^{\crits-1}}{|x|^{s}} &\hbox{ in }\Omega \\
\hfill v &=0   & \hbox{on }  \partial \Omega,
\end{array} \right.
\end{eqnarray}
where  
$x\mapsto b(x)$ is a positive function in $C^{0,1}(\overline{\Omega})$  with $b(0)=\dfrac{(n-2)^{\frac{  2-s}{n-2}}}{2^{2-s}}$ satisfying the following estimates:
\begin{itemize}
\item For $n=3$, we have that
\begin{equation}\label{est:nablab:3}
(x,\nabla b(x))= (\crits+2)b(0)|x|+O(|x|^2)\hbox{ as }x\to 0
\end{equation}
\item For $n=4$, we have $b\in C^1(\overline{\Omega})$, $\nabla b(0)=0$ and
\begin{equation}\label{est:nablab:4}
(x,\nabla b(x))=4(\crits+2) b(0)|x|^2\ln\frac{1}{|x|}+O(|x|^2)\hbox{ as }x\to 0
\end{equation}
\item For $n\geq 5$, we have $b\in C^2(\overline{\Omega})$, $\nabla b(0)=0$ and 
\begin{equation}\label{est:deltab}
\frac{\Delta b(0)}{b(0)}=\frac{4n(2n-2-s)}{n-4}.
\end{equation}
\end{itemize}
{We also have that $h_{\gamma,\lambda} \in C^{1}\left( \overline{\Omega} \setminus \{0\}\right)$ in such a way that there exist $c_3,c_4\in\mathbb{R}$ such that}
\begin{align}\label{def:h}
 h_{\gamma,\lambda}(x) =h_{\gamma,\lambda}(r)=
\left \{ \begin{array} {lc}
        \dfrac{4 \gamma}{r} +{c_3+O(r)}  &\hbox{ for }~ n=3, \\
           $~$\\
              {   8\gamma \ln \dfrac{1}{r}+c_4+O(r)}&\hbox{ for }~ n=4,\\
           $~$\\
4\lambda+\dfrac{4(n-2)}{n-4}\left(\gamma-\dfrac{n(n-4)}{4}\right)+O(r)&\hbox{ for }~ n\geq5.
            \end{array} \right. 
\end{align}
In addition, when $n=4$ and $\gamma=0$, we have that $c_4=4(\lambda-2)$.

Moreover, the  hyperbolic operator  $  \displaystyle -\Delta_{\Bn}-\gamma{V_2}-\lambda$ is coercive if and only if the corresponding Euclidean operator $-\Delta   - \left( \frac{\gamma}{|x|^{2}} + h_{\gamma,\lambda}(x)  \right)   $ is coercive.
\end{lemma}

\noindent We note that equation (\ref{Problem on bounded domain:Euclidean}) have been studied extensively in the case where  $h_{\gamma,\lambda}\equiv\lambda$ and $b\equiv 1$, that is 
\begin{align}\label{two}
\left\{ \begin{array}{llll}
-\Delta u-\gamma \dfrac{u}{|x|^2}-\lambda u&=& \dfrac{|u|^{\crits-2}u}{|x|^s}  \ \ &\text{in } \Ono,\\
\hfill u&=&0 &\text{on }\partial \Omega, 
\end{array} \right.
\end{align}
and more so in the non-singular case (i.e., when $\gamma=s=0$), with major contributions to the existence of ground state solutions by Brezis-Nirenberg \cite{bn} when $n\geq 4$ and $\lambda >0$, and by Druet \cite{d2, d3} and Druet-Laurain \cite{dl} when $n=3$. The multiplicity of solutions was established in that case by Devillanova-Solimini \cite{ds} in dimension $n\geq 7$ and $0<\lambda <\lambda_1(-\Delta)$. 

The existence of ground state solutions when $s=0$, $0<\gamma \leq \frac{(n-2)^2}{4}-1$,  $0<\lambda <\lambda_1(-\Delta -\frac{\gamma}{|x|^2})$, and $0\in \Omega$ was established by Janelli \cite{Jan}, while the multiplicity of solutions for when $0\leq \gamma < \frac{(n-2)^2}{4}-4$ and $n\geq 7$ was proved by Cao-Yan \cite{cy}. We also note that if $\Omega$ is the unit ball, then Catrina-Lavine \cite{cl}  showed that (\ref{two}) admits no radial positive solutions \cite{cl} for $\gamma>\frac{(N-2)^2}{4}-1$, while Esposito et al. show in \cite{egpv} that for $\lambda$ small, problem \eqref{two} has no radial sign-changing solutions in $B$ whenever $\gamma\geq \frac{(N-2)^2}{4}-4$. There are also many results in the case where $\lambda > \lambda_1(-\Delta)$. They are summarized in \cite{egpv}.

The cases where both $\gamma$ and $s$ are not zero, including the limiting cases were also studied by Ghoussoub-Robert \cite{gr3} and many other authors \cite{cp1,pw, ww, yy}.

 The case when $0\in \partial \Omega$ is more involved and was studied extensively in Ghoussoub-Robert \cite{gr4} and by the above authors \cite{gmr}. This paper is mostly concerned with the influence of the potential $h_{\gamma, \lambda}$ on existence and multiplicity issues related to the equation (\ref{Problem on bounded domain:Euclidean}).

\subsection{Existence and Multiplicity of solutions}
Our multiplicity result will therefore follow from the following more general Euclidean statement. 
We set
\begin{equation}
\ctheta:=\left\{\begin{array}{l}
h\in C^1(\overline{\Omega}\setminus\{0\});\ \hbox{ $\exists c\in\rr$ such that}\\
\\
\lim_{x \to 0}  |x|^{\theta}h(x) =  c\hbox{ \& }\lim_{x \to 0}  |x|^{\theta}( x, \nabla h(x) ) =  -c\theta.
\end{array}\right.
\end{equation}

\noindent
The main multiplicity result in this paper is the following theorem

\begin{theorem}\label{th:sc:theta}
Consider a  bounded, smooth domain $\Omega \subset \rn$ ($n\geq 3$)  such that $0\in  \Omega$ and  assume that   $0< s < 2$ and $0\leq \theta < 2$. Let $h\in\ctheta$ be such that $-\Delta - \frac{\gamma}{|x|^2}-h(x)$ is coercive in $\Omega$. Let $b\in C^{2}(\overline{\Omega})$  be such that $b\geq 0$,  $b(0) > 0$ and $\nabla b(0)=0$, while if $\theta =0$, we shall  assume in addition that $\nabla^{2}b(0)\geq 0$ in the sense of bilinear forms. 

If  $ 0\leq \gamma < \dfrac{(n-2)^{2}}{4}-(2-\theta)^{2}$ and $\lim_{x \to 0}  |x|^{\theta}h(x) =  K_{h}>0$, then the boundary value problem  
\begin{align}\label{one}
\left\{ \begin{array}{llll}
-\Delta u-\gamma \dfrac{u}{|x|^2}-h(x) u&=& b(x) \dfrac{|u|^{\crits-2}u}{|x|^s}  \ \ &\text{in } \Ono,\\
\hfill u&=&0 &\text{on }\partial \Omega, 
\end{array} \right.
\end{align}
has an infinite number of possibly sign-changing solutions  in $ \huno$. Moreover, these solutions belong to $ C^{2}(\overline{\Omega}\setminus\{0\})$ while around $0$ they behave like
\begin{align}
\lim_{x\to 0}|x|^{\frac{n-2}{2}-\sqrt{\frac{(n-2)^{2}}{4}-\gamma}}u(x) =K \in \R. 
\end{align}
\end{theorem} 
Note that when $h_{\gamma,\lambda}\equiv\lambda$ and $b\equiv 1$, the above yields an infinite number of solutions for equation (\ref{two}),  
provided $0< s < 2$, $n\geq 5$, $ 0\leq \gamma < \dfrac{(n-2)^{2}}{4}-4$ and $0<\lambda <\lambda_1(-\Delta -\frac{\gamma}{|x|^2})$. Note that the fact that $s>0$ allows for an improvement on the dimension $n\geq 7$ established for $s=0$ by Devillanova-Solimini \cite{ds} when $\gamma=0$ and Cao-Yan \cite{cy} when $\gamma > 0$.\\

\noindent
The multiplicity result  will follow from  standard min-max arguments once we prove the required compactness, which relies on blow-up analysis techniques. The proof consists of  analyzing  the asymptotic behaviour of a family of solutions to   the related  subcritical equations --potentially developing a singularity at zero-- as we approach the critical exponent.

\subsection{Compactness of approximate solutions}
\medskip

\noindent
For  $0\leq \theta < 2$, let $(h_\eps)_{\eps>0}$ be functions in $ C^1(\overline{\Omega} \setminus \{0\})$, and $h_0\in\ctheta$  such that 
\begin{equation}\label{hyp:he:0}
\lim_{\eps\to 0}\sup_{x\in\Omega}\left(|x|^\theta|\he(x)-h_0(x)|+|x|^{\theta+1}|\nabla(\he-h_0)(x)|\right)=0.
\end{equation}

\begin{theorem}\label{th:cpct:sc}
Consider a  bounded, smooth domain $\Omega \subset \rn$, $n\geq 3$,  such that $0\in  \Omega$ and  assume that   $0< s < 2$ and $0\leq \theta < 2$. Let $h_0\in\ctheta$ be such that $-\Delta - \frac{\gamma}{|x|^2}-h_0(x)$ is coercive in $\Omega$ and consider $(h_\eps)_{\eps>0}$ such that \eqref{hyp:he:0} hold. Let $b\in C^{2}(\overline{\Omega})$ be such that $b\geq 0$,  $b(0) > 0$ and $\nabla b(0)=0$, while if $\theta =0$, we shall  assume in addition that $\nabla^{2}b(0)\geq 0$ in the sense of bilinear forms.  

Let $(\pe)_{\eps>0}$ such that $0\leq p_\eps<\crits-2$ and $\lim \limits_{\eps\to 0}\pe=0$ and consider a sequence of functions $(\ue)_{\eps>0}$ that  is  uniformly bounded in $\huno$ and such that for each $\epsilon >0$, $u_\epsilon$ is a solution to the equation:
\begin{equation}
\left\{ \begin{array}{llll}
-\Delta \ue-\gamma \dfrac{\ue}{|x|^2}-\he(x) \ue&=&b(x) \dfrac{|\ue|^{\crits-2-\pe}\ue}{|x|^s}  \ \ &\text{in } \Ono,\\
\hfill \ue&=&0 &\text{on }\partial \Omega.
\end{array} \right.
\end{equation}
If $\gamma < \frac{(n-2)^{2}}{4}-(2-\theta)^{2}$ and $\lim \limits_{|x|\to 0} |x|^{\theta} h_0(x)=K_{h_0}>0$,   then  the sequence  $(\ue)_{\eps>0}$ is  pre-compact in the space $\huno$. 
\end{theorem}
\smallskip
\subsection{The interior mass and compactness for approximate positive solutions at any energy level}
\noindent
Note that the function $x\mapsto |x|^{-\alpha}$ is a solution of  
\begin{equation}\label{homo.Rn.1}
\left(-\Delta -\frac{\gamma}{|x|^2}\right)u=0 \qquad \hbox{ on }  \rn\setminus\{0\},
\end{equation}
if and only if $\alpha\in \{\alpha_-(\gamma),\alpha_+(\gamma)\}$, where $\alpha_-(\gamma),\alpha_+(\gamma)$ are as in  (\ref{def:alpha}). 
Actually, one can show that any non-negative solution $u\in C^2(\rn\setminus \{0\})$ of (\ref{homo.Rn.1}) is of the form 
\begin{equation}
 u(x)=C_- |x|^{-\am}+ C_+ |x|^{-\ap}\hbox{ for all }x\in \rnp,
\end{equation}
where $C_-,C_+\geq 0$, see \cite{gr3}.
\medskip

\noindent
For showing compactness in low-dimensions, the arguments are not any more local, but global. We are naturally led to introducing a notion of mass:
\begin{theorem}[The Euclidean interior mass, see \cite{gr3}]\label{def:mass}
Let $\Omega$ be a smooth bounded domain of $\rn$ such that $0\in\Omega$ is an interior point. Let  $\gamma<\dfrac{(n-2)^2}{4}$ and let $h \in \ctheta$ be such that the operator $-\Delta - \dfrac{\gamma}{|x|^2}-h(x)$ is coercive in $\Omega$. Then there exists $H\in C^\infty(\overline{\Omega}\setminus\{0\})$ such that 
\begin{align}\label{interior.mass} 
\qquad\qquad \left\{\begin{array}{cl}
-\Delta H-\dfrac{\gamma}{|x|^2}H-h(x)H=0 &\hbox{ in } \Omega\setminus\{0\}\\
 H>0&\hbox{ in } \Omega\setminus\{0\}\\
 H=0&\hbox{ on }\partial\Omega.
\end{array}\right.
\end{align}
These solutions are unique up to a positive multiplicative constant, and there exists $ c_{+}>0$ such that 
$$\displaystyle H(x)\simeq_{x\to 0}\dfrac{c_{+}}{|x|^{\alpha_+(\gamma)}}.$$ 
Moreover, when $\gamma>\dfrac{(n-2)^2}{4}-\left(1-\dfrac{\theta}{2}\right)^2$, there exists $c_{-}\in\rr$ such that
\begin{equation}
H(x)=\dfrac{c_{+}}{|x|^{\ap}}+\dfrac{c_{-}}{|x|^{\am}}+o\left(\dfrac{1}{|x|^{\am}}\right)\hbox{ as }x\to 0.
\end{equation}
We define the interior mass as $m_{\gamma,h}(\Omega):=\dfrac{c_{+}}{c_{-}}$, which is independent of the choice of $H$.
\end{theorem}
\noindent
We then establish the following compactness result for positive solutions.

\begin{theorem}\label{th:cpct:sc:pos}
Consider a  bounded, smooth domain $\Omega \subset \rn$, $n\geq 3$,  such that $0\in  \Omega$ and  assume that   $0< s < 2$ and $0\leq \theta < 2$. Let $h_0\in\ctheta$ be such that $-\Delta - \frac{\gamma}{|x|^2}-h_0(x)$ is coercive in $\Omega$ and consider $(h_\eps)_{\eps>0}$ such that \eqref{hyp:he:0} hold. Let $b\in C^{2}(\overline{\Omega})$ be such that $b\geq 0$,  $b(0) > 0$ and $\nabla b(0)=0$, while if $\theta =0$, we shall  assume in addition that $\Delta b(0)\geq0$. 

Let $(\pe)_{\eps>0}$ be such that $0\leq p_\eps<\crits-2$ and $\lim \limits_{\eps\to 0}\pe=0$ and consider a sequence of functions $(\ue)_{\eps>0}$ that  is  uniformly bounded in $\huno$ such that for each $\epsilon >0$, $u_\epsilon$ is a solution to the equation:
\begin{align}
\left\{ \begin{array}{clll}
-\Delta \ue-\gamma \dfrac{\ue}{|x|^2}-\he(x) \ue&=&b(x) \dfrac{\ue^{\crits-1-\pe}}{|x|^s}  \ \ &\text{in } \Ono,\\
\hfill \ue &>&0 &\text{in } \Omega,\\
\hfill \ue&=&0 &\text{on }\partial \Omega.
\end{array} \right.
\end{align}
Assuming one of the following conditions
\begin{itemize}
\item $0\leq \gamma \leq \frac{(n-2)^{2}}{4}-(1-\frac{\theta}{2})^{2}$ and $\lim \limits_{|x|\to 0} |x|^{\theta} h_0(x)=K_{h_0}>0$;

\item  $\gamma >\frac{(n-2)^{2}}{4}-(1-\frac{\theta}{2})^{2}$ and $m_{\gamma,h_0}(\Omega)>0$. 
\end{itemize}
Then  the sequence  $(\ue)_{\eps>0}$ is  pre-compact in the space $\huno$. 
\end{theorem}

\subsection{Non-existence: Stability of the Pohozaev obstruction.} To address issues of non-existence of solutions, we shall prove the following surprising stability of regimes where variational positive solutions do not exist. 

\begin{theorem}\label{thm:non:ter} Let $\Omega$ be a smooth bounded domain in $\rn$ $(n\geq 3)$ such that $0\in\Omega$ is an interior point. Assume that $0< s <2$, $0\leq \theta<2$ and $\gamma<(n-2)^2/4$. Let $h_0\in \ctheta$ be such that $-\Delta-\gamma|x|^{-2}-h_0$ is coercive and let  $b\in C^{2}(\overline{\Omega})$ be such that $b\geq 0$,  $b(0) > 0$ and $\nabla b(0)=0$.

Assume that $\gamma> \dfrac{(n-2)^2}{4}-\left(1-\dfrac{\theta}{2}\right)^2$, the mass  $m_{\gamma,h_0}(\Omega)$ is non-zero, and that there is no positive variational solution to the boundary value problem:
\begin{align}\label{eq:main:pos}
\left\{ \begin{array}{clll}
-\Delta u-\gamma \dfrac{u}{|x|^2}-h_{0}(x) u&=&b(x) \dfrac{u^{\crits-1-\pe}}{|x|^s}  \ \ &\text{in } \Ono,\\
\hfill u &>&0 &\text{in } \Omega,\\
 \hfill u&=&0 &\text{on }\partial \Omega.
\end{array} \right.
\end{align}
 Then, for all $\Lambda>0$, there exists $\epsilon:=\epsilon(\Lambda,h_0)>0$ such that for any $h\in\ctheta$ satisfying
\begin{equation}\label{hperturb}
\sup_{x\in\Omega}(|x|^\theta|h(x)-h_0(x)|+|x|^{\theta+1}|\nabla(h-h_0)(x)|)<\epsilon,
\end{equation}
there is no positive solution to \eqref{one} such that $\Vert \nabla u\Vert_2\leq\Lambda$. 
\end{theorem}

The above result is surprising for the following reason: Assuming $\Omega$ is starshaped with respect to $0$, then the classical Pohozaev obstruction (see Section \ref{sec:app:poho}) yields that  \eqref{one} has no positive variational solution whenever $h_0\in\ctheta$ satisfies
\begin{equation}\label{ineq:h:poho}
\hbox{ $h_0(x)+\dfrac{1}{2}( x,\nabla h_0(x)) \leq 0$ for all $x\in \Omega$.}
 \end{equation}
 We then get the following corollaries.

\begin{coro}\label{thm:non} Let $\Omega$ be a smooth bounded domain in $\rn$ $(n\geq 3)$ that is starshaped around $0$. Assume $0< s <2$, and 
$$ \frac{(n-2)^2}{4}-\left(1-\frac{\theta}{2}\right)^2<\gamma <\frac{(n-2)^2}{4}.$$
If $h_0\in\ctheta$ is a potential satisfying \eqref{ineq:h:poho}, then for all $\Lambda>0$, there exists $\epsilon(\Lambda,h_0)>0$ such that for all $h\in \ctheta$ satisfying (\ref{hperturb}), 
there is no positive solution to \eqref{eq:main:pos} with $b \equiv 1$, such that $\Vert \nabla u\Vert_2\leq\Lambda$. 
\end{coro}

\begin{coro}\label{thm:non:bis} Let $\Omega$ be a smooth bounded domain in $\rn$ $(n\geq 3)$ that is starshaped around $0$. Assume $0< s <2$, and 
$$ \frac{(n-2)^2}{4}-1<\gamma <\frac{(n-2)^2}{4}.$$
Then, for all $\Lambda>0$, there exists $\eps(\Lambda)>0$ such that for all $\lambda\in [0,\epsilon(\Lambda))$, there is no positive solution to equation (\ref{two}) 
with $\Vert \nabla u\Vert_2\leq\Lambda$. 
\end{coro}

\begin{rem} It is worth comparing these results to what happens in the nonsingular case, i.e., when  $\gamma=s=\theta=0$. Note that when $\gamma=\theta=0$, the condition $\gamma> \dfrac{(n-2)^2}{4}-\left(1-\dfrac{\theta}{2}\right)^2$ reads $n=3$. 
 In contrast to the singular case, a celebrated result of Brezis-Nirenberg \cite{bn} shows that, for $\gamma=s=\theta=0$, a variational solution to \eqref{two} always exists whenever $n\geq 4$ and $0<\lambda <\lambda_1(\Omega)$, with the geometry of the domain playing no role whatsoever. On the other hand, Druet-Laurain \cite{dl} showed that the geometry plays a role in dimension $n=3$, still for $\gamma=s=\theta=0$, by proving that when $\Omega$ is star-shaped, then there is no solution to \eqref{two} for all small values of $\lambda>0$ (with no apriori bound on $\Vert \nabla u\Vert_2$). Another point of view is that for $n=3$, the nonexistence of solutions persists under small perturbations, but it does not for $n\geq 4$, i.e., only for $n=3$ that the Pohozaev obstruction is stable in the nonsingular case.
\end{rem}

\noindent When $0\in\partial\Omega$, the situation is different. Indeed, the authors showed in \cite{gmr} that there is no low-dimensional phenomenon and that  in the singular case, the Pohozaev obstruction is stable in all dimensions.

\smallskip\noindent Here are some other extensions related to this phenomenon.

\begin{trivlist}
\item $\bullet$\, Our stability result still holds under an additional smooth perturbation of the domain $\Omega$, just as was done by Druet-Hebey-Laurain \cite{dhl} when $n=3$, $\gamma=s=0$. 
\item $\bullet$\, On the other hand, the stability result of Druet-Laurain \cite{dl} is not conditional on an  apriori bound $\Vert \nabla u\Vert_2$, while ours is. We expect however to be able to get rid of the apriori bound in the singular case, as long as $0\in\Omega$ and $s>0$. This is a project in progress.
\end{trivlist}

\section{Setting  the Blow-up }\label{sec:setup}
  
Throughout this paper, $\Omega$ will denote a bounded, smooth domain in $\rn$, $n\geq 3$, such that $0\in\Omega$.   We will always assume that  $\gamma <\dfrac{(n-2)^2}{4}$, $s\in (0,2)$. For $\eps>0$, we let $\pe\in [0,\crits-2)$ be such that
\begin{align}\label{lim:pe}
\lim_{\eps\to 0}\pe=0. 
\end{align}
We also consider $h_0 \in\ctheta$, $0\leq\theta<2$ such that there exists $K_{h_0}\in\rr$ with
\begin{equation}\label{hyp:h0}
\lim_{x \to 0}  |x|^{\theta}h_0(x) =  K_{h_0} \hbox{ and } \lim_{x \to 0}  |x|^{\theta}(x, \nabla h_0(x) )=  -\theta K_{h_0},
\end{equation}
and $-\Delta - \frac{\gamma}{|x|^2}-h_0(x)$ is coercive in $\Omega$. We also let  $(h_\eps)_{\eps>0}$ be  such that $h_\eps\in C^1(\overline{\Omega} \setminus \{0\})$ for all $\eps>0$ and 
\begin{equation}\label{hyp:he}
\lim_{\eps\to 0}\sup_{x\in\Omega}\left(|x|^\theta|\he(x)-h_0(x)|+|x|^{\theta+1}|\nabla(\he-h_0)(x)|\right)=0.
\end{equation}
Concerning $b$, we assume that
\begin{align}\label{hyp:b}
b(x) \geq 0  \hbox{ with } b\in C^2(\Omega) \hbox{ and }  b(0) > 0,  \nabla b(0)=0. 
\end{align}
The exponents $\alpha_{\pm}(\gamma)$ will denote  
\begin{align}
\alpha_{\pm}(\gamma):=\frac{n-2}{2}\pm\sqrt{\frac{(n-2)^2}{4}-\gamma}.
\end{align}
\medskip

\noindent
We consider a  sequence  of functions $(\ue)_{\eps>0}$ in $\huno$   such that for all $\eps >0$ the function $\ue$ is a solution to the Dirichlet boundary value problem:
$$
\left\{ \begin{array}{llll}
-\Delta \ue-\gamma \dfrac{\ue}{|x|^2}-h_{\epsilon}(x) \ue&=& b(x) \dfrac{|\ue|^{\crits-2-\pe}\ue}{|x|^s}  \ \ &\text{in } \huno ,\\
\hfill \ue&=&0 &\text{on }\partial \Omega.
\end{array}\right.\eqno{(E_\eps)}
$$
where $(\pe)$, $\he(x)$   and $b(x)$ is such that   \eqref{lim:pe},  \eqref{hyp:he} and \eqref{hyp:b}   holds. 
\medskip

\noindent
By the regularity Theorem \ref{th:hopf}, $\ue \in C^{2}(\overline{\Omega}\setminus\{0\})$ and there exists $K_{\eps} \in \R$  such  that  $ \ds \lim_{x\to 0}~ |x|^{\am} \ue (x)=K_{\epsilon}$. In  addition, we assume that the sequence $(\ue)_{\eps>0}$ is bounded in $\huno$ and we let  $\Lambda>0$ be such that
\begin{align}\label{bnd:ue}
\int \limits_{\Omega}  \frac{|\ue|^{\crits -\pe}}{|x|^{s}} dx \leq \Lambda  \qquad \hbox{ for all } \eps>0.
\end{align}
It then follows from the weak compactness of the unit ball of $\huno$ that there exists $u_0\in\huno$ such that
\bequa\label{weak:lim:ue}
\ue\rightharpoonup u_0 \qquad \hbox{ weakly in $\huno$ as $\eps\to 0$.}
\eequa
Then  $u_0$ is a weak solution to the Dirichlet boundary value problem
$$
\left\{ \begin{array}{llll}
-\Delta u_{0}-\gamma \dfrac{u_{0}}{|x|^2}-h_{0}(x) u_{0}&=&b(x) \dfrac{|u_{0}|^{\crits-2}u_{0}}{|x|^s}  \ \ & \hbox{in } \Ono ,\\
\hfill u_{0}&=&0 &\text{on }\partial \Omega.
\end{array}\right. \eqno{(E_0)}
$$
Again from the regularity Theorem \ref{th:hopf} it follows that $u_{0} \in C^{2,\theta}(\overline{\Omega}\setminus\{0\})$ and $ \ds \lim_{x\to 0} |x|^{\am}u_{0} (x)=K_{0} \in \R$.  Fix $\tau\in\rr$ such that
\begin{align}\label{tau}
\am<\tau<\frac{n-2}{2}.
\end{align}
The following proposition  shows  that the sequence $(\ue)_{\eps}$ is  pre-compact in $\huno$  if $x\mapsto|x|^{\tau}\ue$ is uniformly bounded in $L^\infty(\Omega)$. 

\begin{proposition}\label{prop:bounded}
Let $\Omega$ be a smooth bounded domain of $\rn$, $n\geq 3$,  such that $0\in  \Omega$ and  assume that  $0< s < 2$,   $\gamma<\frac{(n-2)^2}{4}$.  We let $(\ue)$, $(\pe)$, $h(x)$ and $b(x)$ be such that $(E_\eps)$,  \eqref{lim:pe}, \eqref{hyp:he},  \eqref{hyp:b}  and \eqref{bnd:ue} holds. Suppose there exists $\tau$ as in \eqref{tau} and $C>0$ such that $|x|^{\tau}|\ue(x)|\leq C$ for all $x\in\Omega\setminus\{0\}$ and for all $\eps >0$. Then up to a subsequence,
$\lim \limits_{\eps\to 0}\ue=u_0$ in $\huno$ where $u_0$ is as  in \eqref{weak:lim:ue}. 
\end{proposition}

\noindent{\it Proof:}  We have assumed that  $|x|^{\tau}|\ue(x)|\leq C$ for all $x\in\Omega$ and for all $\eps >0$. So  the sequence $(\ue)$ is uniformly bounded in $L^{\infty}(\Omega')$ for any $\Omega' \subset \subset \Omegabar\setminus\{0\}$. Then by standard elliptic estimates and from  \eqref{weak:lim:ue} it follows that $\ue \to u_{0}$ in $C^2_{loc}(\Omegabar\setminus\{0\})$. \\
Now since  $|x|^{\tau}|\ue(x)|\leq C$ for all $x\in\Omega$ and for all $\eps >0$ and  since $\tau<\frac{n-2}{2}$, we have 
\begin{align}\label{lim:int:0}
\lim \limits_{\delta\to 0}\lim \limits_{\eps \to 0} \int \limits_{B_{\delta}(0)} b(x) \frac{|\ue|^{\crits-\pe}}{|x|^{s}} dx=0 ~\hbox{ and }  ~\lim \limits_{\delta\to 0}\lim \limits_{\eps \to 0}  \int \limits_{B_{\delta}(0)}  \frac{|\ue|^{2}}{|x|^{2}} dx=0.
\end{align}
Therefore    
\begin{align*}
\lim \limits_{\eps \to 0} \int \limits_{\Omega} b(x)  \frac{|\ue|^{\crits-\pe}}{|x|^{s}} dx =  \int \limits_{\Omega } b(x)  \frac{|u_{0}|^{\crits}}{|x|^{s}} dx ~\hbox{ and }  ~
\lim \limits_{\eps \to 0} \int \limits_{\Omega}  \frac{|\ue|^{2}}{|x|^{2}} dx =  \int \limits_{\Omega }  \frac{|u_{0}|^{2}}{|x|^{2}} dx.
\end{align*}
From $(E_{\eps})$  and \eqref{weak:lim:ue} we  then obtain 
\begin{align*}
\lim \limits_{\eps \to 0} \int \limits_{\Omega} \left(  |\nabla \ue|^2-\gamma  \frac{\ue^2}{|x|^2}- h(x) \ue^2 \right)dx &= \int \limits_{\Omega} \left(  |\nabla u_{0}|^2-\gamma  \frac{u_{0}^2}{|x|^2}- h(x) u_{0}^2 \right)dx  \\
\hbox{  and so then } \lim \limits_{\eps \to 0} \int \limits_{\Omega}   |\nabla \ue|^2 &= \lim \limits_{\eps \to 0} \int \limits_{\Omega}   |\nabla u_{0}|^2.
\end{align*}
And hence  $\lim \limits_{\eps\to 0}\ue=u_0$ in    $\huno$.   \hfill$\Box$
\\

\noindent From now on we shall assume that 
\bequa\label{hyp:blowup}
\lim_{\eps\to 0}\Vert |x|^{\tau }\ue\Vert_{L^\infty(\Omega)}=+\infty ~\hbox{ where }  ~ \am<\tau<\frac{n-2}{2} , 
\eequa
and work towards obtaining a contradiction. We shall say that blow-up occurs whenever  (\ref{hyp:blowup}) holds.

\section{Some Scaling Lemmas}\label{sec:blowuplemmas}

\noindent
We start with two scaling lemmas which we shall use many times in our analysis.

\noindent
\begin{lemma}{\label{scaling lemma 1}}
Let $\Omega$ be a smooth bounded domain of $\rn$, $n\geq 3$,  such that $0\in  \Omega$ and  assume that  $0< s < 2$ and    $\gamma<\frac{(n-2)^2}{4}$. Let $(\ue)$, $(\pe)$, $h(x)$ and $b(x)$ be such that $(E_\eps)$,  \eqref{lim:pe}, \eqref{hyp:he},  \eqref{hyp:b}  and \eqref{bnd:ue} holds.
Let $(y_\eps)_\eps \in \Omega \setminus \{0 \}$ and  let
\begin{align*}
\nu_\eps^{-\frac{n-2}{2}}:= |u_{\eps}(y_{\eps})|, \quad \elle:=\nu_\eps^{1-\frac{\pe}{\crits-2}} ~ \hbox{ and } \quad \kappa_{\eps}:= \left| y_{\eps} \right|^{s/2} \elle^{\frac{2-s}{2}} \qquad \hbox{ for } \eps >0
\end{align*}
Suppose $\lim \limits_{\eps \to 0} y_{\eps}=0$ and $\lim \limits_{\eps \to 0} \nu_{\eps}=0$. Assume that for any $R > 0$ there exists   $C(R)>0$ such that  for all  $\eps >0$ 
\begin{align}\label{scale lem:hyp 1 on u}
|u_{\eps}(x)| \leq& C(R) ~\frac{|y_{\eps}|^{\tau }}{| x|^{\tau }} | u_{\eps}(y_{\eps})| \qquad \hbox{ for all } x\in B_{R \kappa_{\eps}}( y_{\eps}) \setminus \{ 0\}.
\end{align}
Then
\begin{align}
|\ye|=O(\elle) \qquad \hbox{ as } \eps \to 0.  
\end{align}
\end{lemma}

\noindent{\it Proof of Lemma \ref{scaling lemma 1}:} We proceed by contradiction and assume that
\begin{align}
\label{scaling lem 1: contradiction}
\lim_{\eps\to 0}\frac{|\ye|}{\elle}=+\infty.
\end{align}
Then it follows from the definition of $\kappa_{\eps}$ that
\begin{align}
\label{ppty:ke}
\lim \limits_{\eps\to 0}\kappa_{\eps}=0,\; \lim \limits_{\eps\to 0}\frac{\kappa_{\eps}}{\elle}=+\infty \hbox{ and }\lim \limits_{\eps\to 0}\frac{\kappa_{\eps}}{|\ye|}=0.
\end{align}
\medskip

\noindent
Fix a $\rho>0$. We define for all $\eps>0$
$$\ve(x):=\nu_\epsilon^{\frac{n-2}{2}} \ue(\ye+\kappa_{\eps} x) \qquad \hbox{ for } x \in B_{2\rho}(0)$$
Note that this is well defined since $\lim \limits_{\epsilon \to 0} |y_{\eps}|=0 \in \Omega$ and $\lim \limits_{\eps\to 0}\frac{\kappa_{\eps}}{|\ye|}=0$. It follows from \eqref{scale lem:hyp 1 on u} that there exists 
$C(\rho)>0$ such that all $\eps>0$
\begin{align}\label{bnd:lem:1:bis}
|\ve(x)|\leq C(\rho) \frac{1}{\left| \frac{\ye}{|y_{\epsilon}|}+\frac{\kappa_{\eps}}{|y_{\epsilon}|} x \right|^{\tau }} \qquad  \forall x \in B_{2\rho}(0)
\end{align}
using \eqref{ppty:ke} we  then get as $\eps \to 0$
\begin{align*}
|\ve(x)|\leq C(\rho)\left( 1+o(1)\right)   \qquad  \forall x \in B_{2\rho}(0).
\end{align*}
From  equation $(E_\eps)$ we obtain that  $v_{\eps}$ satisfies the equation 
$$-\Delta\ve-\frac{\kappa_{\eps}^2}{|y_{\eps}|^{2}} \frac{\gamma}{\left|\frac{\ye}{|\ye|}+\frac{\kappa_{\eps}}{|\ye|}x\right|^{2}}~\ve-\kappa_{\eps}^2~ \he(\ye+\kappa_{\eps} x)~\ve=b(\ye+\kappa_{\eps} x)\frac{|\ve|^{\crits-2-\pe}\ve}{\left|\frac{\ye}{|\ye|}+\frac{\kappa_{\eps}}{|\ye|}x\right|^s}$$
weakly in $B_{2\rho}(0)$ for all $\eps>0$.  With the help of \eqref{ppty:ke}, \eqref{hyp:he} and  standard elliptic theory it then  follows that there exists $v\in C^1(B_{2\rho}(0))$ such that
$$\lim \limits_{\eps\to 0} \ve = v \qquad \hbox{ in } C^1(B_{\rho}(0)).$$
In particular,
\begin{align}
\label{lim:v:case1:nonzero}
v(0)=\lim \limits_{\eps\to 0} \ve(0)=1
\end{align}
and therefore $v\not\equiv 0$. 
\medskip

\noindent
On the other hand, change of variables and the definition of $\kappa_{\eps}$ yields
\begin{align*}
\int \limits_{ B_{\rho \kappa_{\eps}}(\ye)}\frac{|\ue|^{\crits-\pe}}{|x|^s}~dx = &\frac{|\ue(\ye)|^{\crits-\pe}\kappa_{\eps}^n}{|\ye|^s}\int \limits_{B_{\rho}(0)}\frac{|\ve|^{\crits-\pe}}{\left|\frac{\ye}{|\ye|}+\frac{\kappa_{\eps}}{|\ye|} x\right|^s}~ dx\\
&= \elle^{-\left(1+\frac{2(2-s)}{\crits-2-\pe} \right)}\left(\frac{|\ye|}{\elle}\right)^{s\left(\frac{n-2}{2}\right)}\int \limits_{B_{\rho}(0)}\frac{|\ve|^{\crits-\pe}}{\left|\frac{\ye}{|\ye|}+\frac{\kappa_{\eps}}{|\ye|}x\right|^s}~dx \\
&\geq  \left(\frac{|\ye|}{\elle}\right)^{s \left( \frac{n-2}{2} \right)}\int \limits_{B_{\rho}(0)}\frac{|\ve|^{\crits-\pe}}{\left|\frac{\ye}{|\ye|}+\frac{\kappa_{\eps}}{|\ye|}x\right|^s}~dx.
\end{align*}
Using the equation $(E_\eps)$, \eqref{bnd:ue}, \eqref{scaling lem 1: contradiction}, \eqref{ppty:ke} and  passing to the limit $\eps\to 0$ we get that
$\ds \int \limits_{B_{\rho}(0)}|v|^{\crits}\, dx=0$, 
and so then  $v\equiv 0$ in $B_\rho(0)$, a contradiction with \eqref{lim:v:case1:nonzero}. Thus \eqref{scaling lem 1: contradiction} cannot hold. This proves that $\ye=O(\elle)$ when $\eps\to 0$, which proves the lemma.\hfill$\Box$
\medskip

\noindent
\begin{lemma}{\label{scaling lemma 2}}
Let $\Omega$ be a smooth bounded domain of $\rn$, $n\geq 3$,  such that $0\in  \Omega$ and  assume that  $0< s < 2$ and    $\gamma<\frac{(n-2)^2}{4}$.  Let $(\ue)$, $(\pe)$, $h(x)$ and $b(x)$ be such that $(E_\eps)$,  \eqref{lim:pe}, \eqref{hyp:he},  \eqref{hyp:b}  and \eqref{bnd:ue} holds.
Let $(y_\eps)_\eps \in \Omega \setminus \{0 \}$ and  let
\begin{align*}
\nu_\eps^{-\frac{n-2}{2}}:= |u_{\eps}(y_{\eps})| \quad  \hbox{ and }  \quad \elle:=\nu_\eps^{1-\frac{\pe}{\crits-2}} \qquad \hbox{ for } \eps >0
\end{align*}
Suppose $\nu_{\eps} \to 0$ and $|\ye|=O(\elle) $ as  $\eps \to 0$. 
\medskip

\noindent
For   $\eps>0$ we rescale and define 
\begin{align*}
\we(x):= \nu_{\eps}^{\frac{n-2}{2}} u_{\eps}( \elle x) \qquad \text{ for } x \in \elle^{-1}\Omega\setminus\{0\}.
\end{align*}
Assume that for any $R >\delta> 0$ there exists   $C(R,\delta)>0$ such that  for all  $\eps>0$ 
\begin{align}\label{scale lem:hyp 2 on u}
|\we(x)| \leq& C(R, \delta)  \qquad \text{ for all } x\in B_{R }(0) \setminus \overline{ B_{\delta}(0)}.
\end{align}
\noindent
Then there exists $w \in   {D}^{1,2}(\R^{n}) \cap C^1(\rnp)$ such that  
\begin{align*}
 \we \rightharpoonup  & ~w \qquad \hbox{ weakly in } {D}^{1,2}(\R^n) \quad \text{ as } \eps \rightarrow 0 \notag\\
  \we \rightarrow  & ~w  \qquad \hbox{in } C^{1}_{loc}(\rnp)  \qquad \text{ as } \eps \rightarrow 0 
\end{align*}
And  $ w$  satisfies weakly the equation 
$$-\Delta w-\frac{\gamma}{|x|^2}w= b(0)\frac{|w|^{\crits-2} w}{|x|^s}\hbox{ in }\rnp.$$

\noindent
Moreover if $w \not\equiv 0$, then 
$$\int \limits_{\R^{n}} \frac{|w|^{\crits}}{|x|^{s}} \geq \left(\frac{\mu_{\gamma, s,0}(\R^{n})}{b(0)}\right)^{\frac{\crits}{\crits-2}}$$ 
and  there exists $t \in (0,1]$ such that $\lim \limits_{\eps\to0} \nu_{\eps}^{\pe}=t$.

\end{lemma}

\noindent{\it Proof of Lemma \ref{scaling lemma 2}:}  
The proof proceeds in four steps.

 \noindent {\bf Step \ref{scaling lemma 2}.1:}
Let $\eta \in C^{\infty}_{c}(\R^{n})$. One has  that  $\eta \we \in H^{1}_{0}(\R^{n})$ for $\eps >0 $ sufficiently small. We claim that there exists $w_{\eta} \in D^{1,2}(\R^{n})$ such that upto a subsequence 
\begin{eqnarray*}
 \left \{ \begin{array} {lc}
          \eta \we \rightharpoonup  w_{\eta} \qquad  \quad \text{ weakly in }  D^{1,2} (\R^{n}) ~  \text{ as } \eps \to 0, \\
          \eta \we \rightarrow  w_{\eta}(x) \qquad  a.e ~  ~ \text{ in } \R^{n} ~  \text{ as } \eps \to 0. 
            \end{array} \right. 
\end{eqnarray*} 
\medskip

\noindent
We prove the claim. Let $x \in \R^{n}$, then  
\begin{align*}
\nabla \left( \eta \we \right)(x)=   \we (x) \nabla \eta(x)  + \nu_{\eps}^{\frac{n-2}{2}} \elle ~\eta(x) \nabla \ue(\elle x).
\end{align*}
\medskip\noindent Now for any $\theta >0$, there exists  $C({\theta}) >0$ such that for any $a,b >0$ 
\begin{align*}
(a+b)^{2} \leq C({\theta}) a^{2} + (1+ \theta) b^{2}
\end{align*}
With this inequality we then obtain 
\begin{align*}
\int \limits_{\R^{n}} \left|  \nabla \left( \eta  \we \right)\right|^{2}~ dx \leq C(\theta)  \int \limits_{\R^{n}} | \nabla \eta |^{2} \we^{2} ~ dx + (1 + \theta) \nu_{\epsilon}^{n-2}  \elle^{2} \int \limits_{\R^{n}} \eta^{2} 
\left|\nabla \ue (\elle x) \right|^{2} ~ dx
\end{align*}
With H\"{o}lder inequality and a change of variables  this becomes 
\begin{align}{\label{b'dd on the integral 1}}
\int \limits_{\R^{n}} \left|  \nabla \left( \eta  \we \right)\right|^{2} dx \leq &~C(\theta) \left\| \nabla \eta \right\|^{2}_{L^{n}}  \left( \frac{\nu_{\eps}}{\elle} \right)^{n-2} \left( ~\int \limits_{\Omega} |\ue|^{\crit} ~ dx \right)^{\frac{n-2}{n}}  \notag\\
&+ (1+ \theta)  \left( \frac{\nu_{\eps}}{\elle} \right)^{n-2}  \int \limits_{\Omega} \left( \eta\left(\frac{x}{\elle} \right) \right)^{2}  \left| \nabla \ue \right|^{2}~dx.
\end{align}
Since $\left\| u_{\epsilon} \right\|_{\huno}= O(1)$,  so for $\eps >0$ small enough
\begin{align*} 
\left\| \eta \we \right\|_{D^{1,2} (\R^{n})} \leq C_{\eta}
\end{align*}
Where $C_{\eta}$ is a constant depending on the function $\eta$. The claim then follows from the reflexivity of $D^{1,2} (\R^{n})$.
\medskip 

\noindent {\bf Step \ref{scaling lemma 2}.2:}
Let  ${\eta}_{1} \in C_{c}^{\infty}(\R^{n})$, $0 \leq  {\eta}_{1} \leq 1$ be a smooth cut-off function, such that
 \begin{eqnarray}
{\eta} _{1} = \left \{ \begin{array} {lc}
                 1 \quad \text{for } \ \   x \in B_0(1)\\
                 0 \quad \text{for } \ \   x \in {\R}^n \backslash B_0(2 )
          \end{array} \right.
\end{eqnarray}
For any $R>0$ we let $\eta_{R} = \eta_{1}(x/R)$. Then with a diagonal argument we can assume that upto a subsequence for any $R >0$ there exists  $w_{R} \in D^{1,2} (\R^{n})$ such  that 
\begin{eqnarray*}
 \left \{ \begin{array} {lc}
          \eta_{R} \we \rightharpoonup  w_{R} \qquad  \qquad \text{ weakly in }   D^{1,2} (\R^{n}) ~  \hbox{ as } \eps \to 0 \\
          \eta_{R} \we(x) \rightarrow  w_{R}(x) \qquad a.e ~ x~\text{ in } \R^{n} ~  \text{ as } \eps \to 0 
            \end{array} \right. 
\end{eqnarray*}
Since $\displaystyle{\left\| \nabla  \eta_{R} \right\|^{2}_{n} =\left\| \nabla  \eta_{1} \right\|^{2}_{n} }$ for all $R >0$,   letting $\eps \to 0$ in  $(\ref{b'dd on the integral 1})$ we obtain   that 
\begin{align*}
\int \limits_{\R^{n}_{-}} \left|\nabla w_{R} \right|^{2} dx \leq C \qquad \text{for all } R>0
\end{align*}
where $C$ is a constant independent of $R$. So there exists  $ w \in  D^{1,2} (\R^{n})  $ such  that 
\begin{eqnarray*}
 \left \{ \begin{array} {lc}
           w_{R} \rightharpoonup  w \qquad \qquad \text{ weakly in }  D^{1,2} (\R^{n}) ~  \text{ as } R \to +\infty \\
           w_{R}(x) \rightarrow  w(x)\qquad  a.e ~ x~\text{ in } \R^{n} ~  \text{ as } R \to +\infty 
            \end{array} \right. 
\end{eqnarray*}
\medskip 

\noindent {\bf Step \ref{scaling lemma 2}.3:}
We claim that $ w \in  C^{1}( \R^{n} \setminus \{ 0\})$ and  it satisfies weakly the equation 
$$-\Delta w-\frac{\gamma}{|x|^2}w= b(0)\frac{|w|^{\crits-2} w}{|x|^s}\hbox{ in }\rnp.$$

\noindent
We prove the claim. From   $(E_{\eps})$ it follows that   for any $\eps >0$ and $R>0$, $ \eta_{R} \we $ satisfies  weakly the equation 
\begin{align}{\label{blowup eqn 1}}
-\Delta \left(\eta_{R} \we \right)- \frac{\gamma}{|x|^{2}}\left(\eta_{R} \we \right) -\elle^2~\he(\elle x)\left(\eta_{R} \we \right)=b(\elle x)\frac{|\left(\eta_{R} \we \right)|^{\crits-2-\pe} \left(\eta_{R} \we \right)}{|x|^s}.
\end{align}
From \eqref{scale lem:hyp 2 on u} and \eqref{hyp:he}, using the  standard elliptic estimates  it follows that $w_{R} \in C^{1} \left(  B_{R}(0) \setminus \{ 0\} \right) $ and  that  up to a subsequence 
\begin{align*}
\lim \limits_{\epsilon \to 0} \eta_{R} \we = w_{R}  \qquad \hbox{ in } C^{1}_{loc} \left(  B_{R}(0) \setminus \{ 0\} \right).
\end{align*}
Letting $\eps \to 0$  in eqn $\eqref{blowup eqn 1}$ gives that $ w_{R} $  satisfies  weakly the equation 
$$-\Delta w_{R} - \frac{\gamma}{|x|^{2}}w_{R} = b(0)\frac{|w_{R}|^{\crits-2-\pe} w_{R} }{|x|^s}.$$
Again  we have that $|w_{R}(x)| \leq C(R, \delta)$ for all $x \in \overline{B_{R/2}(0)} \setminus  \overline{B_{2 \delta}(0)}$ and  then  again from standard elliptic estimates it follows that   $ w \in C^{1}( \R^{n} \setminus \{ 0\})$
and  $\lim \limits_{R\to + \infty} \tilde{w}_{R}= \tilde{w}$ in $C^{1}_{loc}(\R^{n} \setminus \{ 0\} )$, up to a subsequence.   Letting $R \to + \infty$  we obtain that $ w $  satisfies  weakly the equation 
$$-\Delta w-\frac{\gamma}{|x|^2}w= b(0)\frac{|w|^{\crits-2} w}{|x|^s}. $$
This proves our claim.

\medskip \noindent {\bf Step \ref{scaling lemma 2}.4:}
Coming back to  equation $\eqref{b'dd on the integral 1}$ we have  for $R>0$
\begin{align}{\label{b'dd on the integral 1b}}
\int \limits_{\R^{n}} \left|  \nabla (\eta_{R}  \we) \right|^{2}~ dx \leq& ~   C(\theta)   \left( ~\int \limits_{B_{0}(2R) \backslash B_{0}(R)} (\eta_{2R}  \we)^{2^{*}} ~ dx \right)^{\frac{n-2}{n}}   \notag\\
&+ (1+ \theta)  \left( \frac{\nu_{\eps}}{\elle} \right)^{n-2}  \int \limits_{\Omega}  \left| \nabla \ue \right|^{2}~dx.\end{align}
Since  the sequence $(\ue)_{\eps}$ is bounded in $\huno$,  letting $\eps \to 0$    and then  $R \to + \infty$  we  obtain  for some constant $C$
\begin{align}
\int \limits_{\R^{n}} \left|  \nabla  w \right|^{2}~ dx  \leq C\left(  \lim \limits_{\epsilon \to 0} \left( \frac{\nu_{\eps}}{\elle} \right)   \right)^{n-2}  .
\end{align}

\noindent
Now if  $w \not\equiv 0$ weakly satisfies  the equation
$$-\Delta w-\frac{\gamma}{|x|^2}w= b(0)\frac{|w|^{\crits-2} w}{|x|^s}$$
using the definition of  $\mu_{\gamma, s,0}(\R^{n})$ it then  follows that 
$$\int \limits_{\R^{n}} \frac{|w|^{\crits}}{|x|^{s}} \geq \left(\frac{\mu_{\gamma, s,0}(\R^{n})}{b(0)}\right)^{\frac{\crits}{\crits-2}}.$$
Hence $\displaystyle \lim \limits_{\epsilon \to 0} \left( \frac{\nu_{\eps}}{\elle} \right)  >0$ which implies that 
\begin{align}
t:= \lim\limits_{\epsilon \to 0} \nu_{\eps} ^{\pe} >0.
\end{align}
Since $\lim \limits_{\eps \to 0} \nu_{\eps}=0 $, therefore  we have that  $ 0 <t \leq 1.$ This completes the lemma.\hfill$\Box$

\section{Construction and Exhaustion of the Blow-up scales}\label{sec:exh}

\noindent
In this section we prove the following proposition:

\begin{proposition}\label{prop:exhaust} 
Let $\Omega$ be a smooth bounded domain of $\rn$, $n\geq 3$,  such that $0\in  \Omega$ and  assume that  $0< s < 2$ and    $\gamma<\dfrac{(n-2)^2}{4}$.  Let $(\ue)$, $(\pe)$, $(\he)$ and $b(x)$ be such that $(E_\eps)$,  \eqref{lim:pe}, \eqref{hyp:he},  \eqref{hyp:b}  and \eqref{bnd:ue} holds.  Assume that blow-up occurs, that is
$$\lim_{\eps\to 0}\Vert |x|^{\tau} \ue\Vert_{L^\infty(\Omega)}=+\infty ~\hbox{ where }  ~ \am<\tau<\frac{n-2}{2} .$$ 

\noindent
Then, there exists $N\in\nn^\star$  families of scales $(\mei)_{\eps>0}$ such that we have: \par

\begin{enumerate}
\item[{\bf(A1)}]
$\lim \limits_{\eps\to 0}\ue=u_0$ in $C^2_{loc}(\overline{\Omega}\setminus\{0\})$ where $u_0$ is as in\eqref{weak:lim:ue}. \\
\item[{\bf(A2)}]
$0< \meun< ...<\meN$, for all $\eps>0$. \\
\item[{\bf(A3)}]
$$\lim_{\eps\to 0}\meN=0\hbox{ and } \lim_{\eps\to 0}\frac{\mu_{i+1,\eps}}{\mei}=+\infty\hbox{ for all } 1 \leq  i \leq N-1.$$ \\
\item[{\bf(A4)}]
For any $1 \leq  i \leq N$ and for   $\eps >0$ we  rescale and define 
$$ \tuei(x):= \mei^{\frac{n-2}{2}}\ue(\kei x) \qquad \hbox{ for } x \in k_{i,\eps}^{-1}\Omega\setminus\{0\} $$
where $\kei=\mei^{1-\frac{\pe}{\crits-2}}$. 

\noindent
Then there exists $\tui \in {D}^{1,2}(\R^{n}) \cap C^1(\rnp)$, $\tui \not\equiv 0$  such that $\tui $ weakly solves the equation 
\begin{equation}\label{eq:tui}
-\Delta \tui-\frac{\gamma}{|x|^2}\tui= b(0)\frac{|\tui|^{\crits-2} \tui}{|x|^s}
\end{equation}
and
\begin{align*}
\tuei \longrightarrow &~ \tui \qquad \hbox{in } C^{1}_{loc}(\rnp)  \qquad \hbox{ as } \eps \to0, \notag\\
\tuei \rightharpoonup  &~ \tui \qquad \hbox{ weakly in } {D}^{1,2}(\R^n)  \quad \hbox{ as } \eps \to0.
\end{align*}\\
\item[{\bf(A5)}]
There exists
$C>0$ such that
$$|x|^{\frac{n-2}{2}}|\ue(x)|^{1-\frac{\pe}{\crits-2}}\leq C$$
for all $\eps>0$ and all $x\in\Omega \setminus \{ 0\}$.\\
\item[{\bf(A6)}]
$\lim \limits_{R\to +\infty} \lim \limits_{\eps\to 0}~\sup \limits_{ \Omega \setminus  B_{R k_{N,\eps}}(0)   } |x|^{\frac{n-2}{2}}|\ue(x)-u_0(x)|^{1-\frac{\pe}{\crits-2}}=0.$ \\
\item[{\bf(A7)}]
$\lim \limits_{R\to +\infty} \lim \limits_{\eps\to 0}~\sup_{  B_{\delta k_{1,\eps}}(0)\setminus \{ 0\}} |x|^{\frac{n-2}{2}}\left|\ue(x)-\mu_{1,\eps}^{-\frac{n-2}{2}} \tu_{1}\left( \frac{x}{k_{1,\eps} }  \right)\right|^{1-\frac{\pe}{\crits-2}}=0.$\\
\item[{\bf(A8)}]
For any $\delta>0$ and any $1 \leq  i \leq N-1$, we have 
$$\lim_{R\to +\infty}\lim_{\eps\to 0}\sup_{\delta k_{i+1, \eps}\geq |x|\geq R \kei}|x|^{\frac{n-2}{2}}\left|\ue(x)-\mu_{i+1,\eps}^{-\frac{n-2}{2}} \tu_{i+1}\left( \frac{x}{k_{i+1,\eps} }  \right)\right|^{1-\frac{\pe}{\crits-2}}=0.$$
\item[{\bf(A9)}]
For any $i\in\{1,...,N\}$ there exists $t_i\in (0,1]$ such that
$\lim_{\eps\to 0}\mei^{\pe}=t_i.$ 
\end{enumerate}
\end{proposition}
\noindent The proof of this proposition proceeds in five steps.\\
\noindent
Since  $s>0$, the subcriticality $\crits<\crit:=\crit(0)$  of equations $(E_{\eps})$ in $\Omega\setminus\{0\}$ along with \eqref{weak:lim:ue} yields that  $\ue \to u_{0}$ in $C^2_{loc}(\Omegabar\setminus\{0\})$. So the only blow-up point is the origin.

\medskip\noindent{\bf Step \ref{sec:exh}.1:}
 The construction of the $\mei$'s proceeds by induction. This step is the initiation.
\medskip

\noindent
By the regularity Theorem \ref{th:hopf} and the definition of $\tau$ it follows that for any $\eps >0$ there exists $\xeun \in \Ono$ such that 
\begin{align}\label{def: xe}
\sup \limits_{ x \in \Ono} |x|^{\tau} |\ue (x)|= |\xeun|^{\tau} |\ue (\xeun)|. 
\end{align}
We define  $\meun$ and $\keun>0$ as follows 
\begin{align}\label{def:me:ke}
\meun^{-\frac{n-2}{2}}:=|\ue(\xeun)| ~ \hbox{ and } ~\keun:=\meun^{1-\frac{\pe}{\crit-2}}.
\end{align}
Since blow-up occurs, that is  (\ref{hyp:blowup}) holds, we have
$$\lim \limits_{\eps \to 0} \meun =0$$
It follows that $\ue$ satisfies the hypothesis \eqref{scale lem:hyp 1 on u} of  Lemma \ref{scaling lemma 1} with $\ye=\xeun$, $\nu_{\eps}=\meun$. Therefore 
$$|\xeun| = O\left(\keun\right)~ \hbox{ as } \eps \to 0. $$
Infact, we claim that  there exists  $c_{1}>0$ such that
\begin{align}\label{def:c1}
\lim \limits_{\eps \to 0} \frac{|\xeun|}{\keun}= c_{1}.
\end{align}
 We argue by contradiction and we assume that $|\xeun|=o(\keun)$ as $\eps \to 0$. We define  for $\eps >0$
 $$\tilde{v}_{\eps}(x):=\meun^{\frac{n-2}{2}}\ue(|\xeun|x)\qquad \hbox { for } x\in |\xeun|^{-1} \Ono $$
Using  $(E_{\eps})$ we obtain  that $\tilde{v}_{\eps}$  weakly satisfies the equation in $\ds  |\xeun|^{-1} \Ono $
$$-\Delta\tilde{v}_{\eps}- \frac{\gamma}{|x|^2} \tilde{v}_{\eps}-|\xeun|^2 \he(|\xeun|x) ~\tilde{v}_{\eps}=b(|\xeun|x) \left(\frac{|\xeun|}{\keun}\right)^{2-s} \frac{|\tilde{v}_{\eps}|^{\crits-2}\tilde{v}_{\eps}}{|x|^s}.$$
The definition \eqref{def: xe} yields $|x|^\tau |\tilde{v}_{\eps}(x)|\leq 1$ for all $x\in |\xeun|^{-1} \Ono $. Standard elliptic theory then yield the existence of $\tilde{v}\in C^2(\rnp)$ such that $\tv_{\eps} \to \tv$ in $C^2_{loc}(\rnp)$ where
$$-\Delta\tilde{v}-\dfrac{\gamma}{|x|^2}\tilde{v}=0 \qquad \hbox{ in } \R^{n}\setminus \{ 0\}.$$ 
In addition, we have that  $ \left|\tilde{v}_{\eps} \left(|\xeun|^{-1}\xeun \right) \right|=1$ and  so $\tv \not \equiv 0$. Also  since $|x|^\tau |\tilde{v}(x)|\leq 1$  in $\rnp$,  we have the bound
$$|\tv(x)| < 2|x|^{-\ap}+ 2|x|^{-\am} \qquad \hbox{ in } \rnp.$$
The classification of positive solutions of $-\Delta v-\frac{\gamma}{|x|^2} v=0$ in $\rnp$ $($see  \eqref{prop:liouville}$)$   yields the existence of  $A,B \in \R$ such that $\tilde{v}(x)=A|x|^{-\ap}+B|x|^{-\am}$ in $\rnp$. Then the pointwise control $|x|^\tau |\tilde{v}(x)|\leq 1\hbox{ in }\rnp$ yields $A=B=0$, contradicting $\tv \not \equiv 0$. This proves the claim \eqref{def:c1}.
\medskip

\noindent
We rescale and define
$$\tu_{1,\eps}(x):=\meun^{\frac{n-2}{2}}\ue(\keun x) \qquad \hbox{ for } x \in k_{1,\epsilon}^{-1}\Omega\setminus\{0\} $$ 
It follows from  \eqref{def: xe} and  \eqref{def:c1}  that $\tueun$ satisfies the hypothesis \eqref{scale lem:hyp 2 on u} of  Lemma \ref{scaling lemma 2} with $\ye=\xeun$, $ \nu_{\eps}=\meun$. Then using lemma  \eqref{scaling lemma 2} we get that  there exists $\tu_{1} \in D^{1.2}(\R^{n}) \cap C^{1}(\R^{n}\setminus \{ 0\})$ weakly satisfying the equation:
$$-\Delta \tu_{1}-\frac{\gamma}{|x|^2}\tu_{1}= b(0)\frac{|\tu_{1}|^{\crits-2} \tu_{1}}{|x|^s}\hbox{ in }\rnp.$$
and
\begin{align*}
 \tu_{1,\eps} \rightharpoonup  & ~\tu_{1} \qquad \hbox{ weakly in } {D}^{1,2}(\R^n) \quad \text{ as } \eps \rightarrow 0 \notag\\
\tu_{1,\eps} \rightarrow  & ~\tu_{1}  \qquad \hbox{in } C^{1}_{loc}(\rnp)  \qquad \text{ as } \eps \rightarrow 0 
\end{align*}
It follows from the definition that $\left| \tu_{1,\eps}\left( \frac{\xeun}{\keun}\right) \right|=1$. From  \eqref{def:c1} we therefore have that $\tu_{1} \not \equiv 0$.  And hence  again from  Lemma \ref{scaling lemma 2} we get that 
$$\int \limits_{\R^{n}} \frac{|\tu_{1} |^{\crits}}{|x|^{s}} \geq \left(\frac{\mu_{\gamma, s,0}(\R^{n})}{b(0)}\right)^{\frac{\crits}{\crits-2}}$$
and there exists $t_{1} \in (0,1]$ such that  $\lim \limits_{\eps\to0}\mu_{1,\eps}^{\pe}=t_{1}$.   
\medskip

\noindent
Next, since $|x|^{\am}\tu_{1}\in C^0(\R^{n})$, we  have
$$\lim \limits_{\delta \to 0}\lim_{\eps\to 0}~\sup_{  B_{\delta k_{1,\eps}}(0)\setminus \{ 0\}}  |x|^{\frac{n-2}{2}}\left|\mu_{1, \eps}^{-\frac{n-2}{2}}\tu_{1}\left(\frac{x}{k_{1,\eps}}\right)\right|^{1-\frac{\pe}{\crits-2}}=0$$
and then using  the definitions  \eqref{def: xe}, \eqref{def:me:ke} it follows that 
$$\lim \limits_{\delta \to 0}\lim_{\eps\to 0}~\sup_{  B_{\delta k_{1,\eps}}(0)\setminus \{ 0\}} |x|^{\frac{n-2}{2}}\left|\ue(x)-\mu_{1,\eps}^{-\frac{n-2}{2}} \tu_{1}\left( \frac{x}{k_{1,\eps} }  \right)\right|^{1-\frac{\pe}{\crits-2}}=0.$$\hfill$\Box$

\noindent{\bf Step \ref{sec:exh}.2:} There exists $C>0$ such that for all $\eps>0$ and all $x\in\Ono$, 
\begin{align}\label{ineq:est:1}
|x|^{\frac{n-2}{2}}|\ue(x)|^{1-\frac{\pe}{\crits-2}} \leq C.
\end{align}

\medskip

\noindent{\it Proof of Step \ref{sec:exh}.2:} We argue by contradiction and let $(\ye)_{\eps>0} \in\Ono$ be such that
\begin{align}
\label{hyp:step32}
\sup_{x\in\Ono}|x|^{\frac{n-2}{2}}|\ue(x)|^{1-\frac{\pe}{\crits-2}}=|\ye|^{\frac{n-2}{2}}|\ue(\ye)|^{1-\frac{\pe}{\crits-2}}\to +\infty ~ \hbox{ as } \eps\to 0.
\end{align}
By the regularity Theorem \ref{th:hopf}  it follows that the sequence $(\ye)_{\eps>0}$ is well-defined and moreover $\lim \limits_{\eps \to 0} y_{\eps} =0$, since $\ue \to u_{0}$ in $C^2_{loc}(\Omegabar\setminus\{0\})$. 
For $\eps >0$ we let
$$ \nu_{\eps}:=|\ue(\ye)|^{-\frac{2}{n-2}},~ \elle:=\nu_{\eps}^{1-\frac{\pe}{\crits-2}} \hbox{ and } \kappa_{\eps}:= \left| y_{\epsilon} \right|^{s/2} \elle^{\frac{2-s}{2}} .$$
Then it follows from (\ref{hyp:step32}) that
\begin{align}\label{lim:infty:34}
\lim_{\eps\to 0}\nu_{\eps}=0,  ~\lim_{\eps\to 0}\frac{|\ye|}{\elle}=+\infty  \hbox{ and } \lim \limits_{\eps\to 0}\frac{\kappa_{\eps}}{|\ye|}=0.
\end{align}

\noindent
Let $R>0$ and let $x\in B_R(0)$ be  such that $\ye+\kappa_{\eps} x\in\Ono$. It follows from the definition (\ref{hyp:step32}) of $\ye$ that for all $\eps>0$
$$|\ye+\kappa_{\eps} x|^{\frac{n-2}{2}}|\ue(\ye+\kappa_{\eps} x)|^{1-\frac{\pe}{\crits-2}}\leq |\ye|^{\frac{n-2}{2}}|\ue(\ye)|^{1-\frac{\pe}{\crits-2}}$$
and then, for all $\eps>0$ 
$$\left(\frac{|\ue(\ye+\kappa_{\eps} x)|}{|\ue(\ye)|}\right)^{1-\frac{\pe}{\crits-2}}\leq\left(\frac{1}{1-\frac{\kappa_{\eps}}{|\ye|}R}\right)^{\frac{n-2}{2}}$$
for all  $x\in B_R(0)$ such that $\ye+\kappa_{\eps} x\in\Ono$. Using  \eqref{lim:infty:34}, we get that there exists $C(R)>0$ such that the hypothesis \eqref{scale lem:hyp 1 on u} of  Lemma \ref{scaling lemma 1}  is satisfied and therefore  one has  $|\ye|=O(\elle)$ when $\eps\to 0$, contradiction to (\ref{lim:infty:34}). This proves (\ref{ineq:est:1}).
\hfill$\Box$
\bigskip

\noindent 
Let $\I \in\mathbb{N}^\star$. We consider the following assertions:
\begin{enumerate}
\item[{\bf(B1)}]
$0< \meun< ...<\mu_{\I,\eps}.$\\
\item[{\bf(B2)}]
 $\lim_{\eps\to 0}\mu_{\eps,\I}=0\hbox{ and }~\lim_{\eps\to 0}\frac{\mu_{\eps,i+1}}{\mei}=+\infty\hbox{ for all } 1 \leq i \leq \I-1$\\
\item[{\bf(B3)}]
For all $ 1 \leq i \leq \I$ there exists $\tui \in {D}^{1,2}(\R^{n}) \cap C^2(\rnp)$  such that $\tui $ weakly solves the equation 
$$-\Delta \tui-\frac{\gamma}{|x|^2}\tui= b(0)\frac{|\tui|^{\crits-2} \tui}{|x|^s} \qquad \hbox{ in } \rnp$$
with 
$$\int \limits_{\R^{n}} \frac{|\tui|^{\crits}}{|x|^{s}} \geq \left(\frac{\mu_{\gamma, s,0}(\R^{n})}{b(0)}\right)^{\frac{\crits}{\crits-2}},$$ 
and
\begin{align*}
\tuei \longrightarrow &~ \tui \qquad \hbox{in } C^{1}_{loc}(\rnp)  \qquad \hbox{ as } \eps \to0, \notag\\
\tuei \rightharpoonup  &~ \tui \qquad \hbox{ weakly in } {D}^{1,2}(\R^n)  \quad \hbox{ as } \eps \to0.
\end{align*}
where for $\eps >0$
$$ \tuei(x):= \mei^{\frac{n-2}{2}}\ue(\kei x) \qquad \hbox{ for } x \in k_{i,\epsilon}^{-1}\Omega\setminus\{0\} $$
with $\kei=\mei^{1-\frac{\pe}{\crits-2}}$.\\
\item[{\bf(B4)}]
For all $1 \leq i \leq \I $, there exists $t_i\in (0,1]$ such that
$\lim_{\eps\to 0}\mei^{\pe}=t_i.$
\end{enumerate}
\medskip

\noindent
We say that ${\mathcal H}_{\I}$ holds if there exists $\I$ sequences  $(\mei)_{\eps>0}$, $i=1,...,\I$ such that points (B1), (B2) (B3) and (B4) holds. Note that it follows from Step \ref{sec:exh}.1 that ${\mathcal H}_{1}$ holds.
Next we  show the following holds:
\medskip

\noindent
{\bf Step \ref{sec:exh}.3} \label{prop:HN}
Let $I\geq 1$. We assume that ${\mathcal H}_{\I}$ holds.
Then either ${\mathcal H}_{\I+1}$ holds or
$$\lim_{R\to +\infty}\lim_{\eps\to 0}\sup_{ \Omega \setminus  B_{0}(R k_{\I,\eps})} |x|^{\frac{n-2}{2}}|\ue(x)-u_0(x)|^{1-\frac{\pe}{\crits-2}}=0.$$
\medskip

\noindent{\it Proof of Step \ref{sec:exh}.3:}
Suppose
$$\lim_{R\to +\infty}\lim_{\eps\to 0}\sup_{ \Omega \setminus  B_{0}(R k_{\I,\eps})} |x|^{\frac{n-2}{2}}|\ue(x)-u_0(x)|^{1-\frac{\pe}{\crits-2}}\neq 0.$$
Then there exists a sequence of points $(\ye)_{\eps>0}\in\Ono$ such that

\bequa\label{hyp:lim:ye:HN}
\lim_{\eps\to 0}\frac{|\ye|}{k_{\I,\eps}}=+\infty\hbox{ and }\lim_{\eps\to 0}|\ye|^{\frac{n-2}{2}}|\ue(\ye)-u_0(\ye)|^{1-\frac{\pe}{\crits-2}}=a>0.
\eequa
Since $\ue \to u_{0}$ in $C^2_{loc}(\Omegabar\setminus\{0\})$ it follows that  $\lim \limits_{\eps \to 0} y_{\eps} =0$. Then by the regularity Theorem \ref{th:hopf} and since $\am < \frac{n-2}{2}$, we get
\begin{align}
\lim_{\eps\to 
0}|\ye|^{\frac{n-2}{2}}|\ue(\ye)|^{1-\frac{\pe}{\crits-2}}=a>0
\end{align}
for some positive constant $a$. In particular, $\lim \limits_{\eps\to 0}|\ue(\ye)|=+\infty$. Let
$$\mu_{\I+1, \eps}:=|\ue(\ye)|^{-\frac{2}{n-2}}\hbox{ and } ~ k_{\I+1,\eps}:=\mu_{\I+1,\eps}^{1-\frac{\pe}{\crits-2}}.$$
As a consequence we have  
\begin{align} \label{hyp:lim:ye:1:HN}
\lim \limits_{\eps\to 0}\mu_{\I +1,\eps}=0 \qquad \hbox{ and } \qquad  \lim \limits_{\eps\to 0}\frac{|\ye|}{k_{\I+1,\eps}}=a>0.
\end{align}
\medskip

\noindent
We rescale and define
$$\tu_{\I+1,\eps}(x):=\mu_{\I+1,\eps}^{\frac{n-2}{2}}\ue(k_{\I+1,\eps}~x) \qquad \hbox{ for } x \in k_{\I+1,\epsilon}^{-1}\Omega\setminus\{0\} $$
It follows from \eqref{ineq:est:1} that for all $\eps >0$
$$|x|^{\frac{n-2}{2}}|\tu_{\I+1,\eps}(x)|^{1-\frac{\pe}{\crits-2}}\leq C \qquad \hbox{ for } x \in k_{\I+1,\epsilon}^{-1}\Omega\setminus\{0\}.$$
so hypothesis \eqref{scale lem:hyp 2 on u} of Lemma \ref{scaling lemma 2} is satisfied. Then using Lemma  \ref{scaling lemma 2} we get that  there exists $\tu_{\I+1} \in D^{1,2}(\R^{n}) \cap C^{1}(\R^{n}\setminus \{ 0\})$ weakly satisfying the equation:
$$-\Delta \tu_{\I+1}-\frac{\gamma}{|x|^2} \tu_{\I+1}= b(0)\frac{|\tu_{\I+1}|^{\crits-2} \tu_{\I+1}}{|x|^s}\hbox{ in }\rnp.$$
and
\begin{align*}
 \tu_{\I+1, \eps} \rightharpoonup  & ~\tu_{\I+1} \qquad \hbox{ weakly in } {D}^{1,2}(\R^n) \quad \text{ as } \eps \rightarrow 0 \notag\\
  \tu_{\I+1, \eps} \rightarrow  & ~\tu_{\I+1} \qquad \hbox{in } C^{1}_{loc}(\rnp)  \qquad \text{ as } \eps \rightarrow 0 
\end{align*}
We denote $\displaystyle \tye:=\frac{\ye}{k_{\I+1,\eps}} $. From \eqref{hyp:lim:ye:1:HN}  it follows that that $\lim \limits_{\eps\to0} |\tye|:= |\tilde{y}_0| =a\neq 0$. Therefore 
$|\tu_{\I+1}(\tilde{y}_0)|=\lim_{\eps\to 0}|\tu_{\I+1,\eps}(\tye)|=1,$
and hence $\tu_{\I+1}\not\equiv 0$. And hence  again from  Lemma \ref{scaling lemma 2} we get
$$\int \limits_{\R^{n}} \frac{|\tu_{\I+1}|^{\crits}}{|x|^{s}} \geq \left(\frac{\mu_{\gamma, s,0}(\R^{n})}{b(0)}\right)^{\frac{\crits}{\crits-2}}$$
and there exists $t_{\I+1}\in (0,1]$ such that 
$\lim \limits_{\eps\to0}\mu_{\I+1,\eps}^{\pe}=t_{\I+1}$. Moreover, it follows from \eqref{hyp:lim:ye:HN} and  \eqref{hyp:lim:ye:1:HN}  that
$$\lim \limits_{\eps\to 0}\frac{\mu_{\I+1,\eps}}{\mu_{\I,\eps}}=+\infty \hbox{ and } \lim_{\eps\to 0}\mu_{\I+1,\eps}=0.$$
Hence the families $(\mei)_{\eps>0}$, $1\leq i \leq \I+1$ satisfy ${\mathcal H}_{\I+1}$.
\hfill$\Box$
\bigskip

 \noindent
The next step is  equivalent to step \ref{sec:exh}.3  at intermediate scales.
\medskip

\noindent
{\bf Step \ref{sec:exh}.4} \label{prop:HN+}
Let $I\geq 1$. We assume that ${\mathcal H}_{\I}$ holds. Then for any $1\leq  i \leq \I-1$ and for any $\delta>0$,  either ${\mathcal H}_{\I+1}$ holds or
$$\lim \limits_{R\to +\infty} \lim \limits_{\eps\to 0} \sup_{B_{\delta k_{i+1,\eps}}(0)\setminus \overline{B}_{R k_{i,\eps}}(0)}|x|^{\frac{n-2}{2}}\left|\ue(x)-\mu_{i+1,\eps}^{-\frac{n-2}{2}}\tu_{i+1}\left(\frac{x}{k_{i+1,\eps}}\right)\right|^{1-\frac{\pe}{\crits-2}}=0.$$
\medskip

\noindent{\it Proof of Step \ref{sec:exh}.4:}
We assume that there exists an $i\leq \I-1$ and   $\delta>0$ such that $$\lim \limits_{R\to +\infty} \lim \limits_{\eps\to 0} \sup_{B_{\delta k_{i+1,\eps}}(0)\setminus \overline{B}_{R k_{i,\eps}}(0)}|x|^{\frac{n-2}{2}}\left|\ue(x)-\mu_{i+1,\eps}^{-\frac{n-2}{2}}\tu_{i+1}\left(\frac{x}{k_{i+1,\eps}}\right)\right|^{1-\frac{\pe}{\crits-2}}>0.$$
It then follows that there exists a sequence  $(\ye)_{\eps>0}\in\Omega$ such that
\beqn
&&\lim_{\eps\to 0}\frac{|\ye|}{k_{i, \eps}}=+\infty,\qquad |\ye|\leq \delta
k_{i+1,\eps}\hbox{ for all }\eps>0\label{hyp:lim:ye:1:bis}\\
&&  |\ye|^{\frac{n-2}{2}}\left|\ue(\ye)-\mu_{i+1,\eps}^{-\frac{n-2}{2}}\tu_{i+1}\left(\frac{\ye}{k_{i+1,\eps}}\right)\right|^{1-\frac{\pe}{\crits-2}}=a>0.\label{hyp:lim:ye:1:ter}
\eeqn
for some positive constant $a$.  Note that $ a< +\infty$ since $$|x|^{\frac{n-2}{2}}\left|\ue(x)-\mu_{i+1,\eps}^{-\frac{n-2}{2}}\tu_{i+1}\left(\frac{x}{k_{i+1,\eps}}\right)\right|^{1-\frac{\pe}{\crits-2}}$$
is uniformly bounded for all $x \in B_{\delta k_{i+1,\eps}}(0)\setminus \overline{B}_{R k_{i,\eps}}(0)$.
\medskip

\noindent
We let  $\tye^{*}\in\R^{n}$ be such that $\ye= k_{i+1,\eps}~\tye^{*}$. It follows from (\ref{hyp:lim:ye:1:bis}) that $|\tye^{*}|\leq \delta$ for all $\eps>0$. . We rewrite (\ref{hyp:lim:ye:1:ter}) as
$$\lim_{\eps\to 0} |\tye^{*}|^{\frac{n-2}{2}}\left|\tu_{i+1, \eps}(\tye^{*})-\tu_{i+1}(\tye^{*})\right|^{1-\frac{\pe}{\crits-2}}=a>0.$$
Then  from point (B3) of ${\mathcal H}_{\I}$ it follows that $\tye^{*}\to 0$ as $\eps\to 0$. And since $|x|^{\am}\tu_{i+1}\in C^0(\R^{n})$, we  get  as $\eps\to 0$
$$|\ye|^{\frac{n-2}{2}}\left|\mu_{i+1, \eps}^{-\frac{n-2}{2}}\tu_{i+1}\left(\frac{\ye}{k_{i+1,\eps}}\right)\right|^{1-\frac{\pe}{\crits-2}}=O\left(\frac{|\ye|}{k_{i+1,\eps}}\right)^{\frac{n-2}{2}-\am}=o(1)$$
Then  (\ref{hyp:lim:ye:1:ter}) becomes 
\begin{align}\label{hyp:lim:ye:bis}
\lim \limits_{\eps\to0}|\ye|^{\frac{n-2}{2}}|\ue(\ye)|^{1-\frac{\pe}{\crits-2}}=a>0.
\end{align}
In particular, $\lim \limits_{\eps\to 0}|\ue(\ye)|=+\infty$. We let
$$\nu_{\eps}:=|\ue(\ye)|^{-\frac{2}{n-2}}\hbox{ and }\elle:=\nu_{\eps}^{1-\frac{\pe}{\crits-2}}.$$
Then we have  
\begin{align} \label{hyp:lim:ye:1:HN+}
\lim \limits_{\eps\to 0}\nu_{\eps}=0 \qquad \hbox{ and } \qquad  \lim \limits_{\eps\to 0}\frac{|\ye|}{\elle}=a>0.
\end{align}
\medskip

\noindent
We rescale and define
$$\tu_{\eps}(x):=\nu_{\eps}^{\frac{n-2}{2}}\ue(\elle~x) \qquad \hbox{ for } x \in \elle^{-1}\Omega\setminus\{0\} $$
It follows from \eqref{ineq:est:1} that for all $\eps >0$
$$|x|^{\frac{n-2}{2}}|\tu_{\eps}(x)|^{1-\frac{\pe}{\crits-2}}\leq C \qquad \hbox{ for } x \in \elle^{-1}\Omega\setminus\{0\}.$$
so hypothesis \eqref{scale lem:hyp 2 on u} of Lemma \ref{scaling lemma 2} is satisfied. Then using lemma  \eqref{scaling lemma 2} we get that  there exists $\tu \in D^{1.2}(\R^{n}) \cap C^{1}(\R^{n}\setminus \{ 0\})$ weakly satisfying the equation:
$$-\Delta \tu-\frac{\gamma}{|x|^2}\tu= b(0)\frac{|\tu|^{\crits-2} \tu}{|x|^s}\hbox{ in }\rnp.$$
and
\begin{align*}
 \tue \rightharpoonup  & ~\tu \qquad \hbox{ weakly in } {D}^{1,2}(\R^n) \quad \text{ as } \eps \rightarrow 0 \notag\\
\tue \rightarrow  & ~\tu  \qquad \hbox{in } C^{1}_{loc}(\rnp)  \qquad \text{ as } \eps \rightarrow 0 
\end{align*}
We denote $\displaystyle \tye:=\frac{\ye}{\elle} $. From \eqref{hyp:lim:ye:bis}  it follows that that $\lim \limits_{\eps\to0} |\tye|:= |\tilde{y}_0| =a\neq 0$. Therefore 
$|\tu(\tilde{y}_0)|=\lim_{\eps\to 0}|\tue(\tye)|=1,$
and hence $\tu \not\equiv 0$. And hence  again from  Lemma \ref{scaling lemma 2} we get
$$\int \limits_{\R^{n}} \frac{|\tu |^{\crits}}{|x|^{s}} \geq \left(\frac{\mu_{\gamma, s,0}(\R^{n})}{b(0)}\right)^{\frac{\crits}{\crits-2}}$$
and there exists $t \in (0,1]$ such that 
$\lim \limits_{\eps\to0}\nu_{\eps}^{\pe}=t$. Moreover it follows from(\ref{hyp:lim:ye:bis}), (\ref{hyp:lim:ye:1:bis})  and since $\lim \limits_{\eps \to 0} \frac{|\ye|}{k_{i+1,\eps}}=0$, that
$$\lim \limits_{\eps\to 0}\frac{\nu_{\eps}}{\mu_{i,\eps}}=+\infty~\hbox{ and }~\lim \limits_{\eps\to 0}\frac{\mu_{i+1, \eps}}{\nu_{\eps}}=+\infty .$$
Hence  the families $(\meun)$,..., $(\mei)$, $(\nu_{\eps})$, $(\mu_{i+1, \eps})$,..., $(\mu_{\I,\eps})$ satisfy ${\mathcal H}_{\I+1}$.\hfill$\Box$
\medskip

\noindent
The last step tells us that family $\{ {\mathcal H}_{\I} \}$ is finite. 
\medskip

\noindent{\bf Step \ref{sec:exh}.5:} Let $N_0=\max\{\I : {\mathcal H}_\I \hbox{ holds } \}$. Then $N_0<+\infty$ and the conclusion of Proposition \ref{prop:exhaust} holds with $N=N_0$.
\medskip

\noindent{\it Proof of Step \ref{sec:exh}.5:}
Indeed, assume that ${\mathcal H}_{\I}$ holds. Since $\mei=o(\mu_{i+1, \eps})$ for all $1 \leq i \leq N-1$,  we get with a change of variable and the
definition of $\tuei$  that for any $R > \delta>0$
\begin{align*}
\int \limits_{\Omega}  \frac{|\ue|^{\crits -\pe}}{|x|^{s}} dx & \geq \sum_{i=1}^{\I} \int \limits_{B_{R \kei}(0)\setminus \overline{B}_{\delta\kei}(0)}  \frac{|\ue|^{\crits-\pe}}{|x|^{s}} dx \notag \\
&\geq  \sum_{i=1}^{\I}  \int \limits_{B_{R}(0)\setminus \overline{B}_{\delta}(0)} \frac{|\tu_{i,\eps}|^{\crits-\pe}}{|x|^{s}} dx. 
\end{align*}
\begin{align}
\hbox{Then  from \eqref{bnd:ue} we have }~\Lambda \geq  \sum_{i=1}^{\I}\int \limits_{B_{R}(0)\setminus \overline{B}_{\delta}(0)} \frac{|\tu_{i,\eps}|^{\crits-\pe}}{|x|^{s}} dx. 
\end{align}
Passing to the limit $\eps\to 0$ and then $\delta \to0$, $R \to +\infty$ we obtain using point (B3) of ${\mathcal H}_{\I}$, that
$$\Lambda \geq  \left(\frac{\mu_{\gamma, s,0}(\R^{n})}{b(0)}\right)^{\frac{\crits}{\crits-2}} \I.$$
 It then follows that $N_0<+\infty$. \hfill $\Box$
\bigskip

\noindent We let families $(\meun)_{\eps>0}$,..., $(\mu_{N_0, \eps})_{\eps>0}$ such that ${\mathcal H}_{N_0}$ holds. We argue
by contradiction and assume that the conclusion of Proposition \ref{prop:exhaust} does not hold with $N=N_0$. Assertions (A1), (A2), (A3), (A4), (A5), (A7) and (A9) holds. Assume that (A6) or (A8) does not hold. It then follows from Steps  $(4.3)$, $(4.4)$ and $(4.5)$ that ${\mathcal H}_{N+1}$ holds. A contradiction with the choice of $N=N_0$ and the proposition is proved.
\hfill$\Box$

\section{Strong Pointwise Estimates}\label{sec:se:1}

\noindent
The objective of this section is the proof of the following strong pointwise control. 

\begin{proposition}\label{prop:fund:est}   Let $\Omega$ be a smooth bounded domain of $\rn$, $n\geq 3$,  such that $0\in  \Omega$ and  assume that  $0< s < 2$ and  $\gamma<\frac{(n-2)^2}{4}$. Let $(\ue)$, $(\pe)$,  $(\he)$ and $b(x)$ be such that $(E_\eps)$,  \eqref{lim:pe}, \eqref{hyp:he},  \eqref{hyp:b}  and \eqref{bnd:ue} holds.  Assume that blow-up occurs, that is
$$\lim_{\eps\to 0}\Vert |x|^{\tau} \ue\Vert_{L^\infty(\Omega)}=+\infty ~\hbox{ where }  ~ \am<\tau<\frac{n-2}{2} .$$ 
Consider the $\mu_{1,\eps},...,\mu_{N,\eps}$  from Proposition  \ref{prop:exhaust}.  Then there exists $C>0$ such that for all $\eps>0$ 
\begin{align}\label{eq:est:global}
|\ue(x)|\leq &~
C \left(~\sum_{i=1}^N\frac{\mei^{\frac{\ap-\am}{2}}}{ \mei^{\ap-\am}|x|^{\am}+|x|^{\ap}}+ \frac{\Vert |x|^{\am}u_0 \vert|_{L^{\infty}(\Omega)}}{|x|^{\am}} \right)
\end{align}
for all $ x \in \Ono$.
\end{proposition}

\noindent The proof of this estimate proceeds in seven steps. 

\medskip

\noindent{\bf Step \ref{sec:se:1}.1:} We claim that for any  $\tau >0$  small and any $R>0$,
 there exists $C(\tau, R)>0$ such that for all $\eps>0$ sufficiently  small
\begin{align}\label{estim:alpha:1}
|\ue(x)|\leq  C(\tau,R)~\left(\frac{\mu_{N,\eps}^{\frac{\ap-\am}{2}-\tau}}{|x|^{\ap-\tau}}+\frac{\Vert |x|^{\am}u_0 \vert|_{L^{\infty}(\Omega)}}{|x|^{\am+\tau}} \right) ~ \hbox{for all }x\in\Omega\setminus \overline{B}_{R k_{N,\eps}}(0).
\end{align}

\smallskip\noindent{\it Proof of  Step \ref{sec:se:1}.1:}  We fix $\gamma'$ such that $\gamma<\gamma'<\frac{(n-2)^2}{4}$. Since the  operator $-\Delta - \frac{\gamma}{|x|^2}-h_{0}(x)$ is coercive, taking $\gamma'$ close to $\gamma$ it follows that the operator $-\Delta-\frac{\gamma'}{|x|^2} -h_{0}(x)$ is  also coercive in  $\Omega$. From Theorem \ref{interior.mass} it get that there exists   $H\in C^\infty(\overline{\Omega}\setminus\{0\})$ such that 
\begin{align} \label{eqn:sup-sol}
 \left\{\begin{array}{ll}
-\Delta H-\frac{\gamma'}{|x|^2}H- h_{0}(x)H=0 &\hbox{ in } \Omega\setminus\{0\}\\
\hfill H>0&\hbox{ in } \Omega\setminus\{0\}\\
\hfill H=0&\hbox{ on }\partial\Omega.
\end{array}\right.
\end{align} 
And we have the following bound on $H$: there exists   $\delta_{1},C_{1}>0$ such that
\begin{align}\label{esti:sup-sol}
\frac{1}{C_{1}}\frac{1}{|x|^{\beta_{+}(\gamma') }} \leq H(x) \leq C_{1}\frac{1}{|x|^{\beta_{+}(\gamma') }} \qquad \hbox{ for all } x  \in B_{2\delta_1}(0).
\end{align}
\medskip

\noindent 
We let $\lambda^{\gamma'}_{1}>0$ be  the first eigenvalue of the coercive operator $-\Delta-\frac{\gamma'}{|x|^2}-h$ on $\Omega$ and we let
$\varphi\in\huno$ be the unique eigenfunction such that
\begin{align}\label{eqn:sub-sol}
\left\{\begin{array}{ll}
-\Delta\varphi-\frac{\gamma'}{|x|^2}\varphi-h_{0}(x) \varphi&~=\lambda^{\gamma'}_{1} \varphi \quad  \hbox{ in } \Omega\\
\hfill \varphi&~>0  \qquad \hbox{ in } \Omega\setminus\{0\}\\
\hfill \varphi&~=0  \qquad \hbox{ on }\partial\Omega.
\end{array}\right.
\end{align}
It follows from  the regularity result, Theorem \ref{th:hopf} that  there exists $C_2,\delta_2>0$ such that
\begin{align}\label{esti:sub-sol}
\frac{1}{C_{2}}\frac{1}{|x|^{\beta_{-}(\gamma') }} \leq \varphi(x) \leq C_{2}\frac{1}{|x|^{\beta_{-}(\gamma') }} \qquad \hbox{ for all } x \in \Omega\cap B_{2\delta_2}(0).
\end{align}

\noindent
We define the operator
$$\mathcal{L}_\eps:=-\Delta-\left(\frac{\gamma}{|x|^2}+\he \right)-b(x)\frac{|\ue|^{\crits-2-p_{\eps}}}{|x|^s}.$$
\medskip

\noindent
{\it Step \ref{sec:se:1}.1.1:} We claim that given any  $\gamma<\gamma'<\frac{(n-2)^2}{4}$  there exist $\delta_0>0$ and $R_0>0$ such that for any  $0<\delta < \delta_0$ and $R>R_0$, we have  for  $\eps>0$ sufficiently small
\begin{align}
\label{ineq:LG}
\mathcal{L}_\eps H(x) &>0,\hbox{ and } \mathcal{L}_\eps \varphi(x)  >0 \qquad \hbox{for all }x\in B_\delta(0)\setminus \overline{B}_{R\keN}(0), \hbox{ if } u_{0} \not\equiv0. \notag\\
& \mathcal{L}_\eps H(x) >0 \qquad \hbox{for all }x\in \Omega\setminus \overline{B}_{R\keN}(0),  \hbox{ if } u_{0} \equiv0. 
\end{align}

\noindent
As one checks for all $\eps >0$ and $ x \neq 0$
\begin{align*}
\frac{\mathcal{L}_\eps H(x)}{ H(x)}= \frac{\gamma' -\gamma}{|x|^2}+(h_{0} -h_{\eps}) -b(x)\frac{|\ue|^{\crits-2-p_{\eps}}}{|x|^s}
\end{align*}
and 
\begin{align*}
\frac{\mathcal{L}_\eps \varphi (x)}{ \varphi(x)}= \frac{\gamma' -\gamma}{|x|^2} +(h_{0} -h_{\eps}) -b(x)\frac{|\ue|^{\crits-2-p_{\eps}}}{|x|^s}+ \lambda^{\gamma'}_{1}.
\end{align*}

\noindent
For $\eps >0$ sufficiently small $\Vert h_{0}- h_{\eps}\Vert_\infty \leq  \frac{\gamma'-\gamma}{4(1+\sup \limits_{\Omega}|x|^{2} )}$ and we choose $0 <\delta_0 < \min\{ 1, \delta_{1}, \delta_{2} \}$ such that 
\begin{align}\label{def:delta0}
&||b||_{L^{\infty} (\Omega)}\delta_{0}^{\left(\crits-2 \right)\left(\frac{n-2}{2}-\am \right)} \Vert |x|^{\am}u_0 \vert|_{L^{\infty}(\Omega)}^{\crits-2} \leq\frac{\gamma'-\gamma}{ 2^{\crits+3}}
\end{align}
It follows from point (A6) of Proposition \ref{prop:exhaust} that, there exists $R_0>0$ such that for any $R>R_0$, we have for all $\eps>0$ sufficiently small
$$|b(x)|^{\frac{1}{\crits-2}} |x|^{\frac{n-2}{2}}|\ue(x)-u_0(x)|^{1-\frac{\pe}{\crits-2}}\leq \left(\frac{\gamma'-\gamma}{2^{\crits+2}}\right)^{\frac{1}{\crits-2}} \qquad \hbox{for all }  x\in \Omega\setminus \overline{B}_{R\keN}(0) $$

\noindent
With this choice of $\delta_0$  and $R_{0}$ we get that for any $0<\delta < \delta_0$ and $R>R_0$, we have  for $\eps>0$ small enough
\begin{align*}
|b(x)|~|x|^{2-s}|\ue(x)|^{\crits-2-\pe}&\leq 2^{\crits-1-\pe}|x|^{2-s}|b(x)||\ue(x)-u_0(x)|^{\crits-2-\pe} \notag\\ 
&+  2^{\crits-1-\pe}|x|^{2-s}|b(x)||u_0(x)|^{\crits-2-\pe} \notag \\
& \leq   \frac{\gamma' -\gamma}{4} \qquad \hbox{ for all } x\in B_\delta(0)\setminus \overline{B}_{R\keN}(0),  \hbox{ if } u_{0} \not\equiv0  \notag\\
\hbox{and }  |b(x)|~|x|^{2-s}|\ue(x)|^{\crits-2-\pe}& \leq \frac{\gamma' -\gamma}{4}    \qquad  \hbox{for all }x\in \Omega\setminus \overline{B}_{R\keN}(0),  \hbox{ if } u_{0} \equiv0. \\
\end{align*}

\noindent
Hence  we obtain that for $\eps>0$ small enough
\begin{align}
\frac{\mathcal{L}_\eps H(x)}{ H(x)} = &~ \frac{\gamma' -\gamma}{|x|^2} +(h_{0}-\he)  -b(x)\frac{|\ue|^{\crits-2-p_{\eps}}}{|x|^s}\notag\\ 
\geq& ~\frac{\gamma' -\gamma}{|x|^2}  +(h_{0}-\he) - \frac{\gamma' -\gamma}{4 |x|^2} \notag\\
>&~0 \qquad \hbox{ for all }  x\in B_\delta(0)\setminus \overline{B}_{R\keN}(0),  \hbox{ if } u_{0} \not\equiv0  \notag\\
\hbox{and }~~\frac{\mathcal{L}_\eps H(x)}{ H(x)}  >&~0  \qquad  \hbox{for all }x\in \Omega\setminus \overline{B}_{R\keN}(0),  \hbox{ if } u_{0} \equiv0.
\end{align} 
And 
\begin{align}
\frac{\mathcal{L}_\eps \varphi (x)}{ \varphi(x)} &~ \geq \frac{\gamma' -\gamma}{|x|^2}  +(h_{0}-\he) -b(x)\frac{|\ue|^{\crits-2-p_{\eps}}}{|x|^s}.  \notag \\
&~\geq ~\frac{\gamma' -\gamma}{|x|^2}  +(h_{0}-\he) - \frac{\gamma' -\gamma}{4 |x|^2} \notag\\
&~>~0 \qquad \hbox{ for all } x\in    B_\delta (0)\setminus \overline{B}_{R\keN}(0)). 
\end{align}
\medskip

\hfill$\Box$
\medskip

\noindent{\it Step \ref{sec:se:1}.1.2:} It follows from point (A4) of Proposition \ref{prop:exhaust} that there exists $C'_1(R)>0$ such that for  all $\eps>0$  small
\begin{align}
|\ue(x)|\leq C'_1(R)  \frac{\meN^{\frac{\beta_{+}(\gamma')-\beta_{-}(\gamma')}{2}}}{|x|^{\beta_{+}(\gamma')}} \qquad \hbox{ for all } x\in \partial B_{R\keN}(0)
\end{align}
By estimate \eqref{esti:sup-sol} on $H$, we then have for some constant $C_1(R)>0$ 
\begin{align}\label{sup:ue:boundary:2}
|\ue(x)|\leq C_1(R) \meN^{\frac{\beta_{+}(\gamma')-\beta_{-}(\gamma')}{2}} H(x)\qquad \hbox{ for all } x\in \partial B_{R\keN}(0).
\end{align}
It follows from point (A1) of Proposition \ref{prop:exhaust}  and the regularity result \eqref{th:hopf}, that there exists $C'_2(\delta)>0$ such that for all $\eps>0$  small
\begin{align}
|\ue(x)|\leq  C'_2(\delta) \frac{\Vert |x|^{\am}u_0 \vert|_{L^{\infty}(\Omega)} }{|x|^{\am}} \qquad \hbox{ for all } x\in \partial B_{\delta}(0),  \hbox{ if } u_{0}\not \equiv0.
\end{align}
And then by the estimate \eqref{esti:sub-sol} on $\varphi$ we then  have for some constant $C_2(\delta)>0$ 
\begin{align}\label{sup:ue:boundary:1}
|\ue(x)|\leq C_2(\delta) \Vert { |x|^{\am}u_0} \Vert_{L^{\infty}(\Omega)}~ \varphi(x) \qquad \hbox{ for all } x\in  \partial B_{\delta}(0)   \hbox{ if } u_{0}\not \equiv0..
\end{align}

\noindent
We now let for $\eps>0$ , 
$$\Psi_{\eps}(x):=  C_1(R)  \meN^{\frac{\beta_{+}(\gamma')-\beta_{-}(\gamma')}{2}} H(x)+C_2(\delta) \Vert |x|^{\am}u_0 \vert|_{L^{\infty}(\Omega)} \varphi(x)   \qquad  \hbox{ for } x  \in \Omega \setminus \{0\}.$$
\noindent
Then  \eqref{sup:ue:boundary:1} and \eqref{sup:ue:boundary:2} implies that for all $\eps>0$  small
\begin{align}\label{control:ue:bord}
|u_{\eps}(x)| \leq \Psi_{\eps}(x) \qquad  \hbox{ for all } x \in \partial(B_\delta(0)\setminus \overline{B}_{R\keN}(0)),   \hbox{ if } u_{0}\not \equiv0.
\end{align}
and if $u_{0}\equiv0$  then 
\begin{align}\label{control:ue:bord:u0}
|u_{\eps}(x)| \leq \Psi_{\eps}(x) \qquad  \hbox{ for all } x \in \partial(\Omega\setminus \overline{B}_{R\keN}(0)).
\end{align}

\noindent
Therefore when $ u_{0}\not \equiv0$ it follows from \eqref{ineq:LG}) and \eqref{control:ue:bord} that for all $\eps>0$ sufficiently small
$$\left\{\begin{array}{ll}
\mathcal{L}_\eps \Psi_\eps \geq 0=\mathcal{L}_\eps \ue & \hbox{ in } B_\delta(0)\setminus \overline{B}_{R\keN}(0) \\
\Psi_\eps \geq \ue& \hbox{ on } \partial(B_\delta(0)\setminus \overline{B}_{R\keN}(0))\\
\mathcal{L}_\eps \Psi_\eps \geq 0=-\mathcal{L}_\eps \ue & \hbox{ in } B_\delta(0)\setminus \overline{B}_{R\keN}(0)\\
\Psi_\eps \geq -\ue& \hbox{ on } \partial(B_\delta(0)\setminus \overline{B}_{R\keN}(0)).
\end{array}\right.$$
and  from \eqref{ineq:LG} and \eqref{control:ue:bord:u0}, in case $u_{0} \equiv 0$,  we have for  $\eps>0$ sufficiently small
$$\left\{\begin{array}{ll}
\mathcal{L}_\eps \Psi_\eps \geq 0=\mathcal{L}_\eps \ue & \hbox{ in } \Omega \setminus \overline{B}_{R\keN}(0) \\
\Psi_\eps \geq \ue& \hbox{ on } \partial(\Omega\setminus \overline{B}_{R\keN}(0))\\
\mathcal{L}_\eps \Psi_\eps \geq 0=-\mathcal{L}_\eps\ue & \hbox{ in } \Omega \setminus \overline{B}_{R\keN}(0)\\
\Psi_\eps \geq -\ue& \hbox{ on } \partial(\Omega \setminus \overline{B}_{R\keN}(0)).
\end{array}\right.$$
Since $\Psi_\eps>0$  and $\mathcal{L}_\eps \Psi_\eps>0$, it follows from the comparison principle of Berestycki-Nirenberg-Varadhan \cite{bnv} that  the operator $\mathcal{L}_\eps$ satisfies the comparison principle on $ B_\delta(0)\setminus \overline{B}_{R\keN}(0)$. Therefore
\begin{align*}
|\ue(x)|\leq&~ \Psi_\eps(x) \qquad  \hbox{ for all } x \in B_\delta(0)\setminus \overline{B}_{R\keN}(0), \notag\\
\hbox{ and }~ |\ue(x)|\leq &~ \Psi_\eps(x) \qquad  \hbox{ for all } x \in \Omega \setminus \overline{B}_{R\keN}(0) ~ \hbox{ if } u_{0} \equiv 0.
\end{align*}
Therefore  when $ u_{0}\not \equiv0$, we have for all $\eps>0$  small
$$|u_{\eps}(x)|\leq  C_1(R)  \meN^{\frac{\beta_{+}(\gamma')-\beta_{-}(\gamma')}{2}} H(x)+C_2(\delta) \Vert |x|^{\am}u_0 \vert|_{L^{\infty}(\Omega)} \varphi(x) $$
for all $x \in B_\delta(0)\setminus \overline{B}_{R\keN}(0)$, for $R$ large and $\delta$ small. 
\medskip

\noindent
Then  when $ u_{0}\not \equiv0$,  using the estimates \eqref{esti:sup-sol} and \eqref{esti:sub-sol} we have or all $\eps>0$  small
\begin{align*}
|u_{\eps}(x)| &~\leq  C_1(R) \frac{ \meN^{\frac{\beta_{+}(\gamma')-\beta_{-}(\gamma')}{2}}   }{|x|^{\beta_{+}(\gamma') }}+
C_2(\delta)  \frac{\Vert {|x|^{\am}u_0} \Vert_{L^{\infty}(\Omega)}}{|x|^{\beta_{-}(\gamma') }} 
\end{align*}
for all $ x \in B_\delta(0)\setminus \overline{B}_{R\keN}(0)$, for $R$ large and $\delta$ small. 
\medskip

\noindent
And if $u_{0} \equiv 0$, then all $\eps>0$  small and $R >0 $ large
$$|u_{\eps}(x)|\leq   C_1(R)  \frac{\meN^{\frac{\beta_{+}(\gamma')-\beta_{-}(\gamma')}{2}}}{|x|^{\beta_{+}(\gamma')}} ~ \hbox{ for all } x \in \Omega \setminus \overline{B}_{R\keN}(0).$$
Taking $\gamma'$ close to $\gamma$, along with  points (A1) and (A4) of Proposition \ref{prop:exhaust} it then follows that estimate \eqref{estim:alpha:1} holds on  $\Omega\setminus \overline{B}_{ Rk_{\eps,N}}(0)$ for all $R>0$.\hfill$\Box$ 
\bigskip

\noindent{\bf Step \ref{sec:se:1}.2:} Let $1 \leq i \leq N-1$. We\ claim that for any $\tau >0$ small  and any $R,\rho>0$, there exists  $C(\tau,R,\rho)>0$ such that  all $\eps>0$.\par
\begin{align}\label{estim:alpha:1:bis}
|\ue(x)|\leq C(\tau,R,\rho)\left(\frac{\mu_{i,\eps}^{\frac{\ap-\am}{2}-\tau}}{|x|^{\ap-\tau}}+\frac{1}{\mu_{i+1,\eps}^{\frac{\ap-\am}{2}-\tau} |x|^{\am +\tau}}  \right)   
\end{align}
for all $x\in B_{\rho k_{i+1,\eps}}(0)\setminus \overline{B}_{R k_{i,\eps}}(0)$.

\smallskip\noindent{\it Proof of  Step \ref{sec:se:1}.2:} We let $i\in \{1,...,N-1\}$. We emulate  the proof of Step \ref{sec:se:1}.1. Fix $\gamma'$ such that $\gamma<\gamma'<\frac{(n-2)^2}{4}$.  Consider the functions $H$ and $\varphi$ defined in  Step \ref{sec:se:1}.1 satisfying \eqref{eqn:sup-sol} and \eqref{eqn:sup-sol} respectively.  
We define the operator
$$\mathcal{L}_\eps:=-\Delta-\left(\frac{\gamma}{|x|^2}+\he \right)-b(x)\frac{|\ue|^{\crits-2-p_{\eps}}}{|x|^s}.$$
\medskip

\noindent{\it Step \ref{sec:se:1}.2.1:} 
We claim that given any  $\gamma<\gamma'<\frac{(n-2)^2}{4}$  there exist $\rho_0>0$ and $R_0>0$ such that for any  $0<\rho < \rho_0$ and $R>R_0$, we have  for  $\eps>0$ sufficiently small
\begin{align}
\label{ineq:LG2}
\mathcal{L}_\eps H(x) >0,\hbox{ and } \mathcal{L}_\eps \varphi(x) >0 \qquad \hbox{for all }x\in B_{\rho k_{i+1,\eps}}(0)  \setminus \overline{B}_{R k_{i, \eps}}(0) 
\end{align}

\noindent
As one checks for all $\eps >0$ and $ x \neq 0$
\begin{align*}
\frac{\mathcal{L}_\eps H(x)}{ H(x)}&~= \frac{\gamma' -\gamma}{|x|^2}  +(h_{0}-\he)-b(x)\frac{|\ue|^{\crits-2-p_{\eps}}}{|x|^s}, \\
\frac{\mathcal{L}_\eps \varphi (x)}{ \varphi(x)}&~\geq \frac{\gamma' -\gamma}{|x|^2}  +(h_{0}-\he) -b(x)\frac{|\ue|^{\crits-2-p_{\eps}}}{|x|^s}.
\end{align*}

\noindent
We choose  $0<\rho_0<1$ such that
\begin{align}\label{def:rho0}
& \rho_0^2\sup \limits_{\Omega}| h_{0}-h_{\eps}|\leq  \frac{\gamma'-\gamma}{4} \quad \hbox{ for all } \eps>0 \hbox{ small and }  \notag \\
& ||b||_{L^{\infty}}  \rho_{0}^{\left(\crits-2 \right)\left(\frac{n-2}{2}-\am \right)} \Vert |x|^{\am} \tu_{i+1} \vert|_{L^\infty(B_2(0)\cap \R^{n})}^{\crits-2} \leq\frac{\gamma'-\gamma}{ 2^{\crits+3}}
\end{align}
It follows from point (A8) of Proposition \ref{prop:exhaust} that there exists $R_0>0$ such that for any $R>R_0$ and any $0<\rho < \rho_0$, we have for all $\eps>0$ sufficiently small
$$|b(x)|^{\frac{1}{\crits-2}} |x|^{\frac{n-2}{2}}\left|\ue(x)-\mu_{i+1,\eps}^{-\frac{n-2}{2}}\tu_{i+1}\left( \frac{x}{k_{i+1,\eps} }  \right)\right|^{1-\frac{\pe}{\crits-2}} \leq  \left(\frac{\gamma'-\gamma}{2^{\crits+2}}\right)^{\frac{1}{\crits-2}} $$
for all $x\in B_{\rho k_{i+1, \eps}}(0)\setminus \overline{B}_{R\kei}(0) $. 
\medskip

\noindent
With this choice of $\rho_0$ and $R_{0}$ we get that for any $0<\rho < \rho_0$ and $R>R_0$, we have  for $\eps>0$ small enough
\begin{align}
|b(x)||x|^{2-s}|\ue(x)|^{\crits-2-\pe} \leq &~ 2^{\crits-1-\pe}|x|^{2-s} |b(x)|\left|\ue(x)-\mu_{i+1,\eps}^{-\frac{n-2}{2}}\tu_{i+1}\left(\frac{x}{k_{i+1, \eps}} \right)\right|^{\crits-2-\pe} \notag\\
&+2^{\crits-1-\pe}\left(\frac{|x|}{k_{i+1,\eps}}\right)^{2-s} |b(x)|\left|\tu_{i+1} \left(\frac{x}{k_{i+1, \eps}}\right) \right|^{\crits-2-\pe} \notag\\
& \leq   \frac{\gamma' -\gamma}{4} \qquad \hbox{ for all } x\in B_{\rho k_{i+1, \eps}}(0)\setminus \overline{B}_{R\kei}(0).
\end{align}
Hence  as in Step \ref{sec:se:1}.1 we have that for $\eps>0$ small enough
\begin{align*}
\frac{\mathcal{L}_\eps H(x)}{ H(x)} >0 ~ \hbox{ and } \frac{\mathcal{L}_\eps \varphi (x)}{ \varphi(x)}>0 ~ \hbox{ for all }x\in B_{\rho k_{i+1,\eps}}(0)  \setminus \overline{B}_{R k_{i, \eps}}(0). 
\end{align*}
\hfill$\Box$
\medskip

\noindent{\it Step \ref{sec:se:1}.2.2:}  let $i\in \{1,...,N-1\}$. It follows from point (A4) of Proposition \ref{prop:exhaust} that there exists $C'_1(R)>0$ such that for  all $\eps>0$  small
\begin{align}
|\ue(x)|\leq C'_1(R)  \frac{\mu_{i,\eps}^{\frac{\beta_{+}(\gamma')-\beta_{-}(\gamma')}{2}}}{|x|^{\beta_{+}(\gamma')}} \qquad \hbox{ for all } x\in \partial B_{R k_{i,\eps}}(0)
\end{align}
And then by the estimate \eqref{esti:sup-sol} on $H$ we have for some constant $C_1(R)>0$ 
\begin{align}\label{sup:ue:boundary:2:bis}
|\ue(x)|\leq C_1(R) \mei^{\frac{\beta_{+}(\gamma')-\beta_{-}(\gamma')}{2}} H(x)\qquad \hbox{ for all } x\in  \partial B_{R\kei}(0).
\end{align}
From point (A4) of Proposition \ref{prop:exhaust} it follows  that there exists $C'_2(\rho)>0$ such that for all $\eps>0$  small
\begin{align}
|\ue(x)|\leq C'_2(\rho)  \frac{1}{\mu_{i+1,\eps}^{\frac{\beta_{+}(\gamma')-\beta_{-}(\gamma')}{2}} |x|^{\beta_{-}(\gamma')}} \qquad \hbox{ for all } x\in \partial B_{\rho k_{i+1, \eps}}(0).
\end{align}
Then by the estimate \eqref{esti:sub-sol} on $\varphi$ we have for some constant $C_2(\delta)>0$ 
\begin{align}\label{sup:ue:boundary:3:bis}
|\ue(x)|\leq C_2(\rho)  \frac{\varphi(x) }{\mu_{i+1,\eps}^{\frac{\beta_{+}(\gamma')-\beta_{-}(\gamma')}{2}} }  \qquad \hbox{ for all } x\in \partial B_{\rho k_{i+1, \eps}}(0).
\end{align}
\noindent
We let for all $\eps>0$ 
$$\tilde{\Psi}_{\eps}(x):=  C_1(R)   \mei^{\frac{\beta_{+}(\gamma')-\beta_{-}(\gamma')}{2}} H(x) +C_2(\rho)  \frac{\varphi(x) }{\mu_{i+1,\eps}^{\frac{\beta_{+}(\gamma')-\beta_{-}(\gamma')}{2}} }   \qquad \hbox{ for } x  \in \Omega \setminus \{0\}$$
Then  \eqref{sup:ue:boundary:2:bis} and \eqref{sup:ue:boundary:3:bis} implies that for all $\eps>0$  small
\begin{align}\label{control:ue:bord:bis}
|u_{\eps}(x)| \leq \tilde{\Psi}_{\eps}(x) \qquad  \hbox{ for all } x \in \partial \left(B_{\rho k_{i+1, \eps}}(0)\setminus \overline{B}_{R\kei}(0 \right).
\end{align}

\noindent
Therefore it follows from \eqref{ineq:LG2} and \eqref{control:ue:bord:bis} that $\eps>0$ sufficiently small
$$\left\{\begin{array}{ll}
\mathcal{L}_\eps \tilde{\Psi}_\eps \geq 0=\mathcal{L}_\eps \ue & \hbox{ in } B_{\rho k_{i+1, \eps}}(0)\setminus \overline{B}_{R\kei}(0) \\
\tilde{\Psi}_\eps \geq \ue& \hbox{ on } \partial \left(B_{\rho k_{i+1, \eps}}(0)\setminus \overline{B}_{R\kei}(0) \right) \\
\mathcal{L}_\eps \tilde{\Psi}_\eps \geq 0=-\mathcal{L}_\eps \ue & \hbox{ in } B_{\rho k_{i+1, \eps}}(0)\setminus \overline{B}_{R\kei}(0) \\
\tilde{\Psi}_\eps \geq -\ue& \hbox{ on } \partial \left(B_{\rho k_{i+1, \eps}}(0)\setminus \overline{B}_{R\kei}(0 )\right).
\end{array}\right.$$
Since $\tilde{\Psi}_\eps>0$  and $\mathcal{L}_\eps \tilde{\Psi}_\eps>0$, it follows from the comparison principle of Berestycki-Nirenberg-Varadhan \cite{bnv} that  the operator $\mathcal{L}_\eps$ satisfies the comparison principle on $ B_{\rho k_{i+1, \eps}}(0)\setminus \overline{B}_{R\kei}(0)$. Therefore
$$|\ue(x)|\leq \tilde{\Psi}_\eps(x) \qquad  \hbox{ for all } x \in B_{\rho k_{i+1, \eps}}(0)\setminus \overline{B}_{R\kei}(0)).$$
So for all $\eps>0$  small
$$|u_{\eps}(x)|\leq    C_1(R)   \mei^{\frac{\beta_{+}(\gamma')-\beta_{-}(\gamma')}{2}} H(x) +C_2(\rho)  \frac{\varphi(x) }{\mu_{i+1,\eps}^{\frac{\beta_{+}(\gamma')-\beta_{-}(\gamma')}{2}} } $$
 for all $ x \in B_{\rho k_{i+1, \eps}}(0)\setminus \overline{B}_{R\kei}(0)$, for $R$ large and $\rho$ small.  Then using the estimates \eqref{esti:sup-sol} and \eqref{esti:sub-sol} we have or all $\eps>0$  small
$$|u_{\eps}(x)|\leq    C_1(R)  \frac{\mu_{i,\eps}^{\frac{\beta_{+}(\gamma')-\beta_{-}(\gamma')}{2}}}{|x|^{\beta_{+}(\gamma')}}+\frac{C_2(\rho)}{\mu_{i+1,\eps}^{\frac{\beta_{+}(\gamma')-\beta_{-}(\gamma')}{2}} |x|^{\beta_{-}(\gamma')}}  ~  \hbox{ for all } x \in B_{\rho k_{i+1, \eps}}(0)\setminus \overline{B}_{R\kei}(0).$$
for $R$ large and $\rho$ small. 

Taking $\gamma'$ close to $\gamma$, along with  point (A4) of Proposition \ref{prop:exhaust} it then follows that estimate \eqref{estim:alpha:1:bis} holds on  $B_{\rho k_{i+1, \eps}}(0)\setminus \overline{B}_{R\kei}(0)$ for all $R, \rho>0$.\hfill$\Box$ 
 \bigskip

\noindent{\bf Step \ref{sec:se:1}.3:}  We\ claim that for any $\tau >0$ small  and any $\rho>0$, there exists  $C(\tau,\rho)>0$ such that  all $\eps>0$.\par
\begin{align}\label{estim:alpha:1:tis}
|\ue(x)|\leq C(\tau,\rho)\frac{1}{\mu_{1,\eps}^{\frac{\ap-\am}{2}-\tau} |x|^{\am +\tau}}   \qquad \hbox{ for all } x\in B_{\rho k_{1,\eps}}(0) \setminus\{ 0\}.
\end{align}

\smallskip\noindent{\it Proof of  Step \ref{sec:se:1}.3:}  Fix $\gamma'$ such that $\gamma<\gamma'<\frac{(n-2)^2}{4}$. Consider the function $\varphi$ defined in  Step \ref{sec:se:1}.1 satisfying  \eqref{eqn:sup-sol}.
\medskip

\noindent{\it Step \ref{sec:se:1}.3.1:} 
We claim that given any  $\gamma<\gamma'<\frac{(n-2)^2}{4}$  there exist $\rho_0>0$  such that for any  $0<\rho < \rho_0$  we have  for  $\eps>0$ sufficiently small
\begin{align}
\label{ineq:LG3}
 \mathcal{L}_\eps \varphi(x) >0 \qquad \hbox{for all } x\in B_{\rho k_{1,\eps}}(0)\setminus\{ 0\}
\end{align}

\noindent
We choose  $0<\rho_0<1$
such that
\begin{align}
& \rho_0^2\sup \limits_{\Omega}| h_{0}-h_{\eps}|\leq  \frac{\gamma'-\gamma}{4} \quad \hbox{ for all } \eps>0 \hbox{ small and }  \notag \\
& ||b||_{L^{\infty}}  \rho_{0}^{\left(\crits-2 \right)\left(\frac{n-2}{2}-\am \right)} \Vert |x|^{\am} \tu_{1} \vert|_{L^\infty(B_2(0)\cap \R^{n})}^{\crits-2} \leq\frac{\gamma'-\gamma}{ 2^{\crits+3}}
\end{align}
It follows from point (A7) of Proposition \ref{prop:exhaust} that for any $0<\rho < \rho_0$, we have for all $\eps>0$ sufficiently small
$$|b(x)|^{\frac{1}{\crits-2}}   |x|^{\frac{n-2}{2}}\left|\ue(x)-\mu_{1,\eps}^{-\frac{n-2}{2}} \tu_{1}\left( \frac{x}{k_{1,\eps} }  \right)\right|^{1-\frac{\pe}{\crits-2}} \leq  \left(\frac{\gamma'-\gamma}{2^{\crits+2}}\right)^{\frac{1}{\crits-2}} $$
for all $x\in B_{\rho k_{1,\eps}}(0)\setminus\{ 0\} $. 
\medskip

\noindent
With this choice of $\rho_0$  we get that for any $0<\rho < \rho_0$  we have  for $\eps>0$ small enough
\begin{align}
|b(x)| |x|^{2-s}|\ue(x)|^{\crits-2-\pe} \leq &~ 2^{\crits-1-\pe}|x|^{2-s} |b(x)|\left|\ue(x)-\mu_{1,\eps}^{-\frac{n-2}{2}}\tu_{1}\left(\frac{x}{k_{1, \eps}} \right)\right|^{\crits-2-\pe} \notag\\
&+2^{\crits-1-\pe}\left(\frac{|x|}{k_{1,\eps}}\right)^{2-s} |b(x)| \left|\tu_{1} \left(\frac{x}{k_{1, \eps}}\right) \right|^{\crits-2-\pe} \notag\\
& \leq   \frac{\gamma' -\gamma}{4} \qquad \hbox{ for all } x\in B_{\rho k_{1,\eps}}(0) \setminus\{ 0\}.
\end{align}
Hence  similarly as in Step \ref{sec:se:1}.1, we obtain that for $\eps>0$ small enough
\begin{align}
\frac{\mathcal{L}_\eps \varphi (x)}{ \varphi(x)}>0  >0 \qquad \hbox{ for all } x\in x\in B_{\rho k_{1,\eps}}(0) \setminus\{ 0\}.
\end{align}
\hfill$\Box$
\medskip

\noindent{\it Step \ref{sec:se:1}.3.2:}     It follows from point (A4) of Proposition \ref{prop:exhaust}  that there exists $C'_2(\rho)>0$ such that for all $\eps>0$  small
\begin{align}
|\ue(x)|\leq C'_2(\rho)  \frac{1}{\mu_{1,\eps}^{\frac{\beta_{+}(\gamma')-\beta_{-}(\gamma')}{2}} |x|^{\beta_{-}(\gamma')}} \qquad \hbox{ for all } x\in \partial B_{\rho k_{1, \eps}}(0)
\end{align}
and then by the estimate \eqref{esti:sub-sol} on $\varphi$ we have for some constant $C_2(\delta)>0$ 
\begin{align}\label{sup:ue:boundary:3:tis}
|\ue(x)|\leq C_2(\rho)  \frac{\varphi(x) }{\mu_{1,\eps}^{\frac{\beta_{+}(\gamma')-\beta_{-}(\gamma')}{2}} }  \qquad \hbox{ for all } x\in  \partial B_{\rho k_{1, \eps}}(0).
\end{align}
\noindent
We let for all $\eps>0$ 
$$\Psi^{0}_{\eps}(x):=  C_2(\rho)  \frac{\varphi(x) }{\mu_{1,\eps}^{\frac{\beta_{+}(\gamma')-\beta_{-}(\gamma')}{2}} }\qquad \hbox{ for } x  \in \Ono. $$
Then   \eqref{sup:ue:boundary:3:tis} implies that for all $\eps>0$  small
\begin{align}\label{control:ue:bord:tis}
|u_{\eps}(x)| \leq \Psi^{0}_{\eps}(x) \qquad  \hbox{ for all } x \in \partial B_{\rho k_{1, \eps}}(0).
\end{align}
\noindent
Therefore it follows from \eqref{ineq:LG3} and \eqref{control:ue:bord:tis} that $\eps>0$ sufficiently small
$$\left\{\begin{array}{ll}
\mathcal{L}_\eps\Psi^{0}_\eps \geq 0=\mathcal{L}_\eps \ue & \hbox{ in }  B_{\rho k_{1,\eps}}(0) \setminus \{ 0\}\\
\Psi^{0}_\eps \geq \ue& \hbox{ on } \partial  B_{\rho k_{1,\eps}}(0) \setminus \{ 0\}\\
\mathcal{L}_\eps \Psi^{0}_\eps \geq 0=-\mathcal{L}_\eps \ue & \hbox{ in }  B_{\rho k_{1,\eps}}(0)\\
\Psi^{0}_\eps \geq -\ue& \hbox{ on } \partial  B_{\rho k_{1,\eps}}(0).
\end{array}\right.$$
Since  the operator $\mathcal{L}_\eps$ satisfies the comparison principle on $ B_{\rho k_{1,\eps}}(0)$. Therefore
$$|\ue(x)|\leq\Psi^{0}_\eps(x) \qquad  \hbox{ for all } x  \in B_{\rho k_{1,\eps}}(0).$$
And so  for all $\eps >0$ small
$$|u_{\eps}(x)|\leq C_2(\rho) \frac{ \varphi(x)}{\mu_{1,\eps}^{\frac{\beta_{+}(\gamma')-\beta_{-}(\gamma')}{2}} }  ~  \hbox{ for all } x  \in B_{\rho k_{1,\eps}}(0)\setminus \{ 0\}.$$
for  $\rho$ small. Using the estimate \eqref{esti:sub-sol} we have or all $\eps>0$  small
$$|u_{\eps}(x)|\leq  \frac{C_2(\rho)}{\mu_{1,\eps}^{\frac{\beta_{+}(\gamma')-\beta_{-}(\gamma')}{2}} |x|^{\beta_{-}(\gamma')}}  ~  \hbox{ for all } x  \in B_{\rho k_{1,\eps}}(0)\setminus\{ 0\}.$$
for  $\rho$ small. It then follows from point (A4) of Proposition \ref{prop:exhaust} that estimate \eqref{estim:alpha:1:tis} holds on  $x  \in B_{\rho k_{1,\eps}}(0) $ for all $\rho>0$.\hfill$\Box$

\medskip
\noindent{\bf Step \ref{sec:se:1}.4:}  Combining the previous three steps,  it follows from \eqref{estim:alpha:1}, \eqref{estim:alpha:1:bis}, \eqref{estim:alpha:1:tis} and Proposition \ref{prop:exhaust} that for any $\tau>0$ small, there exists $C(\tau)>0$ such that for all $\eps >0$  we have 
\begin{align}\label{eq:est:alpha:global}
|\ue(x)|\leq
C(\tau) \left(~\sum_{i=1}^N \frac{\mei^{\frac{\ap-\am}{2}-\tau}}{ \mei^{\left(\ap-\am\right)-2\tau}|x|^{\am+\tau}+|x|^{\ap -\tau}}+ \frac{\Vert |x|^{\am}u_0 \vert|_{L^{\infty}(\Omega)}}{|x|^{\am +\tau}} \right) 
\end{align}
for all $x \in \Omega \setminus\{ 0\}$.
\medskip

\noindent
Next we improve the above estimate and show that one can take  $\tau =0$ in  \eqref{eq:est:alpha:global}.
\bigskip

\noindent
For $\eps$ small, we let $G_{\epsilon}$ be the Green's function for the coercive operator  $-\Delta-\frac{\gamma}{|x|^{2}}-h_{\epsilon}$ on $\Omega$ with Dirichlet boundary condition. Green's representation formula, the pointwise bounds on the Green's function \eqref{est:G:up} [see Ghoussoub-Robert \cite{gr3}] yields for any $ z \in \Omega$
\begin{align}\label{estim:global:z}
\ue(z) =& \int \limits_\Omega G_{\eps}(z,x) b(x) \frac{|\ue (x)|^{\crits-2-\pe} ~\ue(x)}{|x|^s} ~dx \notag \\
|\ue(z)| \leq & \int \limits_\Omega G_{\eps}(z,x) |b(x)| \frac{|\ue (x)|^{\crits-1-\pe} }{|x|^s}~dx \notag\\
\leq&~ C   \int \limits_\Omega  \left(\frac{\max\{|z|,|x|\}}{\min\{|z|,|x|\}}\right)^{\am} \frac{ 1}{|x-z|^{n-2}} \frac{|\ue(x)|^{\crits-1-\pe}}{|x|^s}~ dx.
\end{align}
Using \eqref{eq:est:alpha:global} we then obtain with  $0<\tau < \frac{\crits-2}{\crits-1}\left( \frac{\ap-\am}{2}\right)$ that 
\begin{align}
&\qquad |\ue(z)|  \leq  \notag\\
& C \sum_{i=1}^N   \int \limits_\Omega  \left(\frac{\max\{|z|,|x|\}}{\min\{|z|,|x|\}}\right)^{\am} \frac{ 1}{|x-z|^{n-2}|x|^{s}} \left( \frac{\mei^{\frac{\ap-\am}{2}-\tau}}{ \mei^{\left(\ap-\am\right)-2\tau}|x|^{\am+\tau}+|x|^{\ap -\tau}}\right)^{\crits-1-\pe}dx  \notag\\
&+ C \Vert |x|^{\am}u_0 \vert|^{\crits-1-\pe}_{L^{\infty}(\Omega)}\int \limits_\Omega  \left(\frac{\max\{|z|,|x|\}}{\min\{|z|,|x|\}}\right)^{\am} \frac{ 1}{|x-z|^{n-2}|x|^{s}}  \frac{ 1}{|x|^{(\am +\tau)(\crits -1-\pe)}} ~dx.   
\end{align}
The first term in the above integral was computed for each bubble in Ghoussoub-Robert in \cite{gr3} when $\pe=0$. The proof goes exactly the same with $\pe>0$. The last last term is straightforward to estimate. We then get that there exists a constant $C>0$ such that  for any sequence of points $(\ze)$ in $\Omega \setminus\{ 0\}$ we have 
\begin{align}
|\ue(\ze)| \leq ~C \left(~\sum_{i=1}^N\frac{\mei^{\frac{\ap-\am}{2}}}{ \mei^{\ap-\am}|\ze|^{\am}+|\ze|^{\ap}}+ \frac{\Vert |x|^{\am}u_0 \vert|_{L^{\infty}(\Omega)}}{|\ze|^{\am}} \right).
\end{align}
This completes the proof of Proposition \ref{prop:fund:est}. \hfill $\Box$

\section{Sharp Blow-up rates and proof of Compactness}\label{sec:cpct}
When the expression makes sense, we define
\begin{align}
\mathcal{C}_{n,s}:=\frac{\left(\frac{2-\theta}{2} \right) \frac{1}{t_{N}^{\frac{n-\theta }{\crits -2}}}  \int \limits_{\R^{n}} \frac{\tu_{N}^{2}}{|x|^{\theta}} ~dx }{\frac{1}{\crits} \left( \frac{n-s}{\crits}\right) \left(\sum \limits_{i=1}^{N}  \frac{1}{t_{i}^{ \frac{n-2}{\crits-2}}} \int \limits_{\rn }   \frac{|\tu_{i}|^{\crits }}{|x|^{s}} ~ dx\right)}
\end{align}
The proof of compactness rely on the following two key propositions.
\begin{proposition}\label{prop:rate:sc}
Let $\Omega$ be a smooth bounded domain of $\rn$, $n\geq 3$,  such that $0\in \Omega $ and  assume that  $0< s < 2$,   $\gamma<\frac{(n-2)^2}{4}$. Let $(\ue)$, $(\he)$, $(\pe)$ and $b$ be such that $(E_\eps)$, \eqref{hyp:he}, \eqref{lim:pe},  \eqref{hyp:b} and \eqref{bnd:ue} holds. Assume that blow-up occurs, that is
\begin{equation}\label{blowup}
\lim_{\eps\to 0}\Vert |x|^{\tau} \ue\Vert_{L^\infty(\Omega)}=+\infty ~\hbox{ for some }  ~ \am<\tau<\frac{n-2}{2} .
\end{equation} 
Consider the $\mu_{1,\eps},...,\mu_{N,\eps}$  and $t_{1},...,t_{N}$  from Proposition \ref{prop:exhaust}. Suppose  that
\begin{equation}
\hbox{either }\{\,\ap-\am>4-2\theta\,\}\hbox{ or }\{\ap-\am>2-\theta\hbox{ and }u_0\equiv 0\}.
\end{equation}
Then, we have the following blow-up rate:
\begin{equation}\label{rate:1}
\lim \limits_{\eps \to 0} \frac{\pe}{\mu_{N,\eps}^{2-\theta}}  =-  \frac{\mathcal{C}_{n,s}}{b(0)}\left(K_{h_0} +\frac{{\bf 1}_{\theta=0}}{\crits} \sum \limits_{i,j=1}^{n}\partial_{ij}b(0)\frac{\int \limits_{\rn}X^iX^j\frac{|\tilde{u}_N|^{\crits}}{|X|^s}\, dX}{ \int \limits_{\R^{n}}  \tu_{N}^{2}  ~dx}\right). 
\end{equation}
\end{proposition}

\begin{proposition}[\textit{The positive case}]\label{prop:rate:sc:2}
Let $\Omega$ be a smooth bounded domain of $\rn$, $n\geq 3$,  such that $0\in \Omega $ and  assume that  $0< s < 2$,   $0 \leq \gamma<\frac{(n-2)^2}{4}$.  
Let $(\ue)$, $(\he)$, $(\pe)$ and $b$ be as in Proposition \ref{prop:rate:sc}.
Assume that blow-up occurs as in \eqref{blowup}. 
Consider  $\mu_{1,\eps},...,\mu_{N,\eps}$  and $t_{1},...,t_{N}$  from Proposition \ref{prop:exhaust}. Suppose in addition that
\begin{align*}
\ue>0 \qquad \hbox{for all } \eps >0.
\end{align*}
Then,  we have the following blow-up rates:

\begin{enumerate}
\item
When $\ap-\am > 4-2\theta$, we have 
\begin{align*}
\lim \limits_{\eps \to 0} \frac{\pe}{\mu_{N,\eps}^{2-\theta}}  =-  \frac{\mathcal{C}_{n,s} }{b(0)}\left(K_{h_0}+{\bf 1}_{\theta=0}\frac{(n-2)(n(n-4)-4\gamma)}{4n(2n-2-s)}\frac{\Delta b(0)}{b(0)}\right).  \\
\end{align*}
\item
When $\ap-\am =4-2\theta$, we have 
\begin{align*}
&\lim \limits_{\eps \to 0} \frac{\pe}{\mu_{N,\eps}^{2-\theta}}  =-  \frac{\mathcal{C}_{n,s} }{b(0)}\left(K_{h_0}+{\bf 1}_{\theta=0}\frac{(n-2)(n(n-4)-4\gamma)}{4n(2n-2-s)}\frac{\Delta b(0)}{b(0)}\right)  ~~ \hbox{ if } u_0\equiv 0.\\
&\lim \limits_{\eps \to 0} \frac{\pe}{\mu_{N,\eps}^{2-\theta}}=- \frac{ \mathcal{C}_{n,s}}{b(0)}\left(K_{h_0}+\tilde{K}+{\bf 1}_{\theta=0}\frac{(n-2)(n(n-4)-4\gamma)}{4n(2n-2-s)}\frac{\Delta b(0)}{b(0)}\right) ~~ \hbox{if } u_0 >0, 
\end{align*}
for some $\tilde{K}>0$.\\
\item
When $\ap-\am<4-2\theta$, we have $u_0\equiv 0$ and
\begin{align*}
\lim \limits_{\eps \to 0} \frac{\pe}{\mu_{N,\eps}^{2-\theta}} ~& =-  \frac{\mathcal{C}_{n,s} }{b(0)}\left(K_{h_0}+{\bf 1}_{\theta=0}\frac{(n-2)(n(n-4)-4\gamma)}{4n(2n-2-s)}\frac{\Delta b(0)}{b(0)}\right)\\   &\hbox{ if }\ap-\am>2-\theta.
\end{align*}
\begin{align}\label{asymp:crit:log}
\lim \limits_{\eps \to 0} \frac{\pe}{\mu_{N,\eps}^{2-\theta}\ln\frac{1}{\mu_{N,\eps}}}~=-\mathcal{C}'_{n,s}  \frac{K_{h_0}}{b(0)} \quad \hbox{ if }\ap-\am=2-\theta,
\end{align}
where
$$\mathcal{C}'_{n,s}  :=\frac{\left(\frac{2-\theta}{2} \right)\mathcal{K}^{2} \omega_{n-1} }{\frac{1}{\crits} \left( \frac{n-s}{\crits}\right) \left(\sum \limits_{i=1}^{N}  \frac{1}{t_{i}^{ \frac{n-2}{\crits-2}}} \int \limits_{\rn }   \frac{|\tu_{i}|^{\crits }}{|x|^{s}} ~ dx\right)}~, $$
with $\mathcal{K}$ as defined in \eqref{def:K}. 
\begin{align}\label{est:mass}
\lim_{\eps\to 0}\frac{p_\eps}{\mu_{N,\eps}^{\ap-\am}}=-\chi \cdot m_{\gamma,h_{0}}(\Omega)\quad \hbox{ if }\ap-\am<2-\theta,
\end{align}
where $\chi>0$ is a constant and $m_{\gamma,h_{0}}(\Omega)$ is the  mass of $\Omega$ associated with the operator $-\Delta - \dfrac{\gamma}{|x|^2}-h_{0}(x)$, defined in Theorem \ref{def:mass}.
\end{enumerate}
\end{proposition}

 \medskip\noindent {\it Proof of Theorems \ref{th:cpct:sc} and \ref{th:cpct:sc:pos}:} We argue by contradiction and assume that the family is not pre-compact. Then, up to a subsequence, it blows up. We then apply Propositions \ref{prop:rate:sc} and \ref{prop:rate:sc:2} to get the blow-up rate (that is nonegative). However, the hypothesis of Theorems \ref{th:sc:theta}, \ref{th:cpct:sc} and \ref{th:cpct:sc:pos} yield exactly negative blow-up rates. This is a contradiction, and therefore the family is pre-compact. This proves the three compactness theorems.\qed

\medskip\noindent We now establish Propositions \ref{prop:rate:sc} and \ref{prop:rate:sc:2}. The proof is divided in {\color{red} 14} steps in Sections \ref{loc pohoaev id} to \ref{pf blow-up rates}. These steps are numbered Steps P1, P2, etc. 

\section{Estimates on the localized Pohozaev identity}\label{loc pohoaev id}

\smallskip\noindent We  let $(\ue)$, $(\he)$, $(\pe)$ and $b$ be such that $(E_\eps)$, \eqref{hyp:he}, \eqref{lim:pe},  \eqref{hyp:b} and \eqref{bnd:ue} holds. Assume that blow-up occurs as in \eqref{blowup}. Note that 
$$\gamma<\dfrac{(n-2)^2}{4}-\left(1-\dfrac{\theta}{2}\right)\; \iff \;\ap-\am>2-\theta,$$
In the following, we will use the following consequence of (A9) of Proposition \ref{prop:exhaust}: for all $i=1,...,N$, there exists $c_i>1$ such that
\begin{equation}\label{bnd:k:mu}
c_i^{-1} \mu_{i,\eps}\leq k_{i,\eps}\leq c_i \mu_{i,\eps}. 
\end{equation}

\begin{step}[Pohozaev identity]\label{step:p1} 
For $\delta_0>0$ small, we define 
\begin{equation}\label{def:rho}
\rho_\eps:=\left\{\begin{array}{cc}
\delta_0 &\hbox{ if }u_0\equiv 0,\\
\re:=\sqrt{\mu_{N,\eps}}&\hbox{ if }u_0\not\equiv 0
\end{array}\right.
\end{equation}
and \begin{align}\label{def:Fe}
F_{\eps}(x):=&~( x, \nu) \left( \frac{|\nabla \ue|^2}{2}  -\frac{\gamma}{2}\frac{\ue^{2}}{|x|^{2}} - \frac{h_{\eps}(x)}{2} \ue^{2}  - \frac{b(x)}{\crits -p_{\eps}} \frac{|\ue|^{\crits-p_{\eps} }}{|x|^{s}} \right)  \notag\\
& ~-  \left(x^i\partial_i \ue+\frac{n-2}{2} \ue \right)\partial_\nu \ue.
\end{align}
We then have
\begin{align}\label{PohoId3}
& \int \limits_{B_{ \rho_{\eps} }(0)\setminus B_{ k^{2}_{1,\eps} }(0)} \left( h_{\eps}(x) +\frac{ \left( \nabla h_{\eps}, x \right)}{2} \right)\ue^{2}~ dx   \notag\\
&  +\frac{p_{\eps}}{\crits} \left( \frac{n-s}{\crits -p_{\eps}}\right) \int \limits_{B_{ \rho_{\eps} }(0)\setminus B_{ k^{2}_{1,\eps} }(0)}   b(x)\frac{|\ue|^{\crits-p_{\eps} }}{|x|^{s}} ~ dx \notag \\
&+ \frac{1}{\crits -\pe}  \int \limits_{B_{ \rho_{\eps} }(0)\setminus B_{ k^{2}_{1,\eps} }(0)} \left( x, \nabla b(x) \right) \frac{|\ue|^{\crits-\pe }}{|x|^{s}} ~dx\ \notag\\
 & ~ ~=~ - \int \limits_{ \partial  B_{ \rho_{\eps} }(0)  } F_{\eps}(x)~d\sigma + \int \limits_{\partial  B_{ k^{2}_{1,\eps} }(0) }   F_{\eps}(x)~d\sigma. 
\end{align}
\end{step}

\medskip\noindent{\it Proof of Step P1:} 
We apply the Pohozaev identity \eqref{PohoId} with $y_{0}=0$ and $U_{\eps}= B_{ \rho_{\eps} }(0)\setminus B_{ k^{2}_{1,\eps} }(0)$. This ends the proof of Step \ref{step:p1}. 
\qed
 
\medskip\noindent We will estimate each of the terms in the above integral identities and calculate the limit as $\eps \to 0$.  \medskip

\subsection{Estimates of the $L^{\crits}$ term in the localized Pohozaev identity}

\begin{step}\label{step:p2} Assuming that $b$ satisfies \eqref{hyp:b}, we claim as $\eps \to 0$

\begin{equation}\label{cpct:2*term}
 \int \limits_{B_{ \rho_{\eps} }(0)\setminus B_{ k^{2}_{1,\eps} }(0)} b(x)
 \frac{|\ue|^{\crits-p_{\eps} }}{|x|^{s}} ~ dx =  \sum \limits_{i=1}^{N}  \frac{b(0)}{t_{i}^{ \frac{n-2}{\crits-2}}} \int \limits_{\rn }   \frac{|\tu_{i}|^{\crits }}{|x|^{s}} ~ dx +o(1),
\end{equation}
where $(\rho_\eps)$ is as in \eqref{def:rho} and the $\tu_{i}$'s are as defined in Proposition \ref{prop:exhaust} $(A4)$.
\end{step}

\noindent {\it Proof of Step \ref{step:p2}:} For any $R,\rho >0$ we decompose  the above integral as
\begin{align*}
 \int \limits_{B_{ r_{\eps} }(0)\setminus B_{ k^{2}_{1,\eps} }(0)}      b(x)\frac{|\ue|^{\crits-p_{\eps} }}{|x|^{s}} ~ dx&= 
 \int \limits_{  B_{ r_{\eps} }(0) \setminus \overline{B}_{R k_{N,\eps}}(0)  }    b(x)\frac{|\ue|^{\crits-p_{\eps} }}{|x|^{s}} ~ dx \notag\\
 &+ \sum  \limits_{i=1}^{N}   \int \limits_{  B_{R k_{i,\eps}}(0)  \setminus \overline{B}_{\rho k_{i,\eps}}(0)  }     b(x) \frac{|\ue|^{\crits-p_{\eps} }}{|x|^{s}} ~ dx  \notag\\
&+ \sum  \limits_{i=1}^{N-1}   \int \limits_{   B_{\rho k_{i+1,\eps}}(0)  \setminus \overline{B}_{R k_{i,\eps}}(0) }     b(x) \frac{|\ue|^{\crits-p_{\eps} }}{|x|^{s}} ~ dx  \notag\\
&  +  
 \int \limits_{B_{R k_{1,\eps}}(0)\setminus B_{ k^{2}_{1,\eps} }(0) }   b(x)   \frac{|\ue|^{\crits-p_{\eps} }}{|x|^{s}} ~ dx. 
\end{align*}
We will evaluate each of the above terms and calculate the limit $\lim\limits_{R \to + \infty} \lim\limits_{\rho \to 0} \lim\limits_{\eps \to 0}$. 
\medskip

\noindent
From the estimate \eqref{eq:est:global}, we  get as $\eps \to 0$
\begin{align*}
 &\int \limits_{   B_{ r_{\eps} }(0) \setminus \overline{B}_{R k_{N,\eps}}(0) }   \left| b(x)\frac{|\ue|^{\crits-p_{\eps} }}{|x|^{s}}\right| ~ dx  \notag\\
\leq &~  C~\int \limits_{   B_{ r_{\eps} }(0) \setminus \overline{B}_{R k_{N,\eps}}(0) } \left[  \frac{ \mu_{N,\eps}^{\frac{\ap-\am}{2}(\crits-\pe)}}{|x|^{\ap(\crits-\pe)+s}}+ 
\frac{1}{ |x|^{\am(\crits-\pe)+s}} \right] ~dx\notag\\
\leq &~  C~\int \limits_{   B_{ r_{\eps}}(0)  \setminus \overline{B}_{{R k_{N,\eps}}}(0)  }   \frac{ \mu_{N,\eps}^{\frac{\ap-\am}{2}(\crits-\pe)}}{|x|^{\ap(\crits-\pe)+s}} ~dx \notag \\
+ &~  C~\int \limits_{  B_{ r_{\eps}}(0)  \setminus \overline{B}_{{R k_{N,\eps}}}(0)  }    \frac{1}{ |x|^{\am(\crits-\pe)+s}}~dx \notag \\
\leq &~  C~ \int \limits_{    B_{ \frac{ r_{\eps}}{k_{N,\eps}}}(0)  \setminus \overline{B}_{R }(0)  }    \frac{ 1}{|x|^{n+ \crits \left( \frac{\ap-\am}{2}\right)-\pe\ap}}~dx \notag\\
+ &~  C~\int \limits_{  B_{ 1}(0)  \setminus \overline{B}_{\frac{R k_{N,\eps}}{r_{\eps}}}(0)  }    \frac{ 1 }{ |x|^{n-\crits \left( \frac{\ap-\am}{2}\right)-\pe\am}} ~dx \notag \\
 \leq &C\left( R^{-\crits \left( \frac{\ap-\am}{2}\right)-\pe\ap}  + r_{\eps}^{\crits \left( \frac{\ap-\am}{2}\right)+\pe\am}\right)  .
\end{align*}
Therefore 
\begin{align}{\label{cpct:2*term:1}}
\lim\limits_{R \to + \infty}  \lim\limits_{\eps \to 0}\int \limits_{   B_{ r_{\eps} }(0) \setminus \overline{B}_{R k_{N,\eps}}(0) }  b(x) \frac{|\ue|^{\crits-p_{\eps} }}{|x|^{s}} ~ dx =0.
\end{align}
\medskip

\noindent
It follows from Proposition \ref{prop:exhaust}  with  a change of variable $x=k_{i,\eps}y$, that for any $1 \leq i\leq N$
\begin{align}{\label{cpct:2*term:2}}
\lim\limits_{R \to + \infty} \lim\limits_{\rho \to 0} \lim\limits_{\eps \to 0}  \int \limits_{  B_{R k_{i,\eps}}(0)  \setminus \overline{B}_{\rho k_{i,\eps}}(0)  }     b(x)\frac{|\ue|^{\crits-p_{\eps} }}{|x|^{s}} ~ dx 
= \frac{b(0)}{t_{i}^{ \frac{n-2}{\crits-2}}} \int \limits_{\rn }   \frac{|\tu_{i}|^{\crits }}{|x|^{s}} ~ dx.
\end{align}
\medskip

\noindent
Let $1 \leq i\leq N-1$. In Proposition \ref{prop:fund:est}, we had obtained the following pointwise estimates: For any $R, \rho>0$ and all $\eps>0$ we have 
\begin{align*}
|\ue(x)| \leq  C~ \frac{\mei^{\frac{\ap-\am}{2}} }{|x|^{\ap}}+\frac{ C }{ \mu_{ i+1,\eps}^{\frac{\ap-\am}{2}}|x|^{\am}}
\end{align*}
for all $x \in B_{\rho k_{i+1,\eps}}(0)  \setminus \overline{B}_{R k_{i,\eps}}(0) $.
\medskip

\noindent
Then we have  as $\eps \to 0$
\begin{align*}
 &\int \limits_{   B_{\rho k_{i+1,\eps}}(0)  \setminus \overline{B}_{R k_{i,\eps}}(0)}     \left| b(x)\frac{|\ue|^{\crits-p_{\eps} }}{|x|^{s}}\right|  ~ dx  \notag\\
 \leq &~  C~\int \limits_{   B_{\rho k_{i+1,\eps}}(0)  \setminus \overline{B}_{R k_{i,\eps}}(0)}   \left[   \frac{ \mei^{\frac{\ap-\am}{2}(\crits-\pe)}}{|x|^{\ap(\crits-\pe)+s}}+ \frac{ \mu_{i+1,\eps}^{-\frac{\ap-\am}{2}(\crits-\pe)} }{ |x|^{\am(\crits-\pe)+s}}  \right]~dx \notag\\
  \leq &~  C~\int \limits_{  B_{ \rho k_{i+1,\eps}}(0)  \setminus \overline{B}_{R k_{i,\eps}}(0)  }  \left[   \frac{ \mei^{\frac{\ap-\am}{2}(\crits-\pe)}}{|x|^{\ap(\crits-\pe)+s}}+ \frac{ \mu_{i+1,\eps}^{-\frac{\ap-\am}{2}(\crits-\pe)} }{ |x|^{\am(\crits-\pe)+s}}  \right] \,dx \notag \\
\leq &~  C~ \int \limits_{  B_{ \frac{\rho k_{i+1,\eps}}{k_{i,\eps}}}(0)  \setminus \overline{B}_{R }(0)  }    \frac{ 1}{|x|^{n+ \crits \left( \frac{\ap-\am}{2}\right)-\pe\ap}}\,dx \notag\\
 &+ ~  C~\int \limits_{  B_{2\rho }(0)  \setminus \overline{B}_{ \frac{R k_{i,\eps}}{k_{i+1,\eps}}}(0)  }    \frac{ 1 }{ |x|^{n-\crits \left( \frac{\ap-\am}{2}\right)-\pe\am}} \, dx \notag\\
 \leq &C\left( R^{-\crits \left( \frac{\ap-\am}{2}\right)-\pe\ap}  + \rho^{\crits \left( \frac{\ap-\am}{2}\right)+\pe\am}\right).
\end{align*}
And so 
\begin{align}{\label{cpct:2*term:3}}
\lim\limits_{R \to + \infty} \lim\limits_{\rho \to 0} \lim\limits_{\eps \to 0} \int \limits_{   B_{\rho k_{i+1,\eps}}(0)  \setminus \overline{B}_{R k_{i,\eps}}(0)}    b(x) \frac{|\ue|^{\crits-p_{\eps} }}{|x|^{s}} ~ dx =0.
\end{align}
\medskip

\noindent
Again, from the pointwise  estimates of  Proposition \ref{prop:fund:est}, we have as $\eps \to 0$
\begin{align*}
&\int \limits_{B_{\rho k_{1,\eps}}(0)\setminus B_{ k^{2}_{1,\eps} }(0)}     \left| b(x)\frac{|\ue|^{\crits-p_{\eps} }}{|x|^{s}}\right| ~ dx \notag\\
\leq &~  C~\int \limits_{B_{\rho k_{1,\eps}}(0)\setminus B_{ k^{2}_{1,\eps} }(0)}   \frac{ \mu_{1,\eps}^{-\frac{\ap-\am}{2}(\crits-\pe)} }{ |x|^{\am(\crits-\pe)+s}}  ~dx \notag\\
\leq& ~  C~\int \limits_{   B_{\rho }(0)  \setminus B_{ k_{1,\eps}^{2}}(0)  }    \frac{ 1 }{ |x|^{n-\crits \left( \frac{\ap-\am}{2}\right)-\pe\am}} ~dx \notag\\
 \leq &C   ~\rho^{\crits \left( \frac{\ap-\am}{2}\right)+\pe\am}.
\end{align*}
Therefore
\begin{align}{\label{cpct:2*term:4}}
 \lim\limits_{\rho \to 0} \lim\limits_{\eps \to 0}   \int \limits_{ B_{\rho k_{1,\eps}}(0)\setminus B_{ k^{2}_{1,\eps} }(0) }   b(x)  \frac{|\ue|^{\crits-p_{\eps} }}{|x|^{s}}~dx=0.
\end{align}
\medskip

\noindent
Combining  \eqref{cpct:2*term:1}, \eqref{cpct:2*term:2}, \eqref{cpct:2*term:3} and \eqref{cpct:2*term:4} we obtain \eqref{cpct:2*term}.

\medskip\noindent We now prove \eqref{cpct:2*term} under the assumption that $u_0\equiv 0$. We decompose the integral as
\begin{align*}
 \int \limits_{   B_{\delta_0 }(0)\setminus B_{ k^{2}_{1,\eps} }(0)}      b(x)\frac{|\ue|^{\crits-p_{\eps} }}{|x|^{s}} ~ dx&=
 \int \limits_{B_{ \delta_0 }(0) \setminus \overline{B}_{r_\eps}(0)}    b(x)\frac{|\ue|^{\crits-p_{\eps} }}{|x|^{s}} ~ dx\notag\\
 & \qquad +  \int \limits_{   B_{r_\eps }(0)\setminus B_{ k^{2}_{1,\eps} }(0)}     b(x) \frac{|\ue|^{\crits-p_{\eps} }}{|x|^{s}} ~ dx,  
\end{align*}
with $r_\eps:=\sqrt{\mu_{N,\eps}}$. From the estimate \eqref{eq:est:global} and $u_0\equiv 0$, we  get as $\eps \to 0$
\begin{equation*}
 \int \limits_{  B_{ \delta_0 }(0) \setminus \overline{B}_{r_\eps}(0) }  \left| b(x) \frac{|\ue|^{\crits-p_{\eps} }}{|x|^{s}} \right| ~ dx 
\leq  C~\int \limits_{  B_{ \delta_0 }(0) \setminus \overline{B}_{r_{\eps}}(0) } \left[  \frac{ \mu_{N,\eps}^{\frac{\ap-\am}{2}(\crits-\pe)}}{|x|^{\ap(\crits-\pe)+s}} \right] ~dx
\end{equation*}
Since $\ap\crits+s>n$, we then get that
$$ \int \limits_{  B_{ \delta_0 }(0) \setminus \overline{B}_{r_\eps}(0)} \left| b(x)   \frac{|\ue|^{\crits-p_{\eps} }}{|x|^{s}}  \right|~ dx \leq C\left(\frac{\mu_{N,\eps}}{r_\eps}\right)^{\frac{\crits}{2}(\ap-\am)}=o(1)$$
as $\eps\to 0$. Then with \eqref{cpct:2*term} we get \eqref{cpct:2*term} and this proves \eqref{cpct:2*term}.
\smallskip\noindent This ends the proof of Step \ref{step:p2}.\qed
\medskip

\begin{step}\label{step:p:grad} Assuming that $b$ satisfies \eqref{hyp:b} and taking $(\rho_\eps)$ as in \eqref{def:rho}, we obtain as $\eps \to 0$
\begin{itemize}
\item If $\frac{\crits}{2}(\ap-\am)>2$, then, we have that
\begin{align}\label{eq:gradb}
 \int \limits_{ B_{\rho_\eps }(0)\setminus B_{ k^{2}_{1,\eps} }(0)}(x,\nabla b(x))
 \frac{|\ue|^{\crits-p_{\eps} }}{|x|^{s}} ~ dx =&~  \mu_{N,\eps}^2\left(\frac{\partial_{ij}b(0)}{2 t_N^{\frac{n}{\crits-2}}}\sum \limits_{i,j=1}^{n}\int \limits_{\rn}X^iX^j\frac{|\tilde{u}_N|^{\crits}}{|X|^s}\, dX\right) \notag\\ &~~+o(\mu_{N,\eps}^{2}).
\end{align}

\item If $\frac{\crits}{2}(\ap-\am)=2$ then, we have that
\begin{align}\label{eq:gradb:log}
 \int \limits_{ B_{\rho_\eps }(0)\setminus B_{ k^{2}_{1,\eps} }(0)}(x,\nabla b(x))
 \frac{|\ue|^{\crits-p_{\eps} }}{|x|^{s}} ~ dx =  O\left(\mu_{N,\eps}^2\ln\frac{1}{\mu_{N,\eps}}\right).
\end{align}

\item If $\frac{\crits}{2}(\ap-\am)<2$, then
\begin{align}\label{eq:gradb.3}
 \int \limits_{ B_{r_\eps }(0)\setminus B_{ k^{2}_{1,\eps} }(0)}(x,\nabla b(x))
 \frac{|\ue|^{\crits-p_{\eps} }}{|x|^{s}} ~ dx = o\left( \mu_{N,\eps}^{\ap-\am}\right)
\end{align}

\end{itemize}
\end{step}

\noindent {\it Proof of Step \ref{step:p:grad}:} The proof is very similar to Step \ref{step:p3} and we only sketch it. For convenience, for any $i=1,...,N$, we define
\begin{equation}
B_{i,\eps}(x):=\frac{\mei^{\frac{\ap-\am}{2}}}{\mei^{\ap-\am}|x|^{\am}+|x|^{\ap}}.
\end{equation}
We have the following useful estimates:
\begin{lemma}\label{est:Bp} For any $i\in \{1,...,N\}$ and $(\re)$ such that $c\mei\leq \re\leq \delta_0$, we have that for $k\geq 0$
\begin{eqnarray*}
&&\int \limits_{B_{\re}(0)}\frac{|x|^kB_{i,\eps}^{\crits-p_\eps}}{|x|^s}\, dx\\
&&\leq C\left\{\begin{array}{cl}
\mei^k &\hbox{ if }\frac{\crits}{2}(\ap-\am)>k;\\
\mei^k\left(1+\left|\ln\frac{\re}{\mei}\right|\right) &\hbox{ if }\frac{\crits}{2}(\ap-\am)=k;\\
\mei^{\frac{\crits}{2}(\ap-\am)}\re^{k-\frac{\crits}{2}(\ap-\am)} &\hbox{ if }\frac{\crits}{2}(\ap-\am)<k.
\end{array}\right.
\end{eqnarray*}
and
\begin{eqnarray*}
&&\int \limits_{B_{\re}(0)}|x|^k\left(1+\left|\ln\frac{1}{|x|}\right|\right)\frac{B_{i,\eps}^{\crits-p_\eps}}{|x|^s}\, dx\\
&&\leq C\left\{\begin{array}{cl}
\mei^k \ln\frac{1}{\mei}&\hbox{ if }\frac{\crits}{2}(\ap-\am)>k;\\
\mei^k\ln\frac{1}{\mei}\left(1+\left|\ln\frac{\re}{\mei}\right|\right) &\hbox{ if }\frac{\crits}{2}(\ap-\am)=k;\\
\mei^{\frac{\crits}{2}(\ap-\am)}\re^{k-\frac{\crits}{2}(\ap-\am)} \left(1+\left|\ln\frac{1}{|\re|}\right|\right)&\hbox{ if }\frac{\crits}{2}(\ap-\am)<k.
\end{array}\right.
\end{eqnarray*}
\end{lemma}
The proof follows from straightforward estimates and the change of variable $x=\mei z$. Note that  we have used that $\mei^{p_\eps}\to t_i\neq 0$ as $\eps\to 0$. 

\medskip

\noindent
We start with the case $\dfrac{\crits}{2}(\ap-\am)>2$. With the same change of variable, we get that
\begin{equation}
\int \limits_{\left(B_{\re}(0)\setminus B_{R k_{N,\eps}}(0)\right)\cup B_{R^{-1} k_{N,\eps}}(0)}\frac{|x|^2B_{N,\eps}^{\crit-p_\eps}}{|x|^s}\, dx\leq \epsilon(R)\mu_{N,\eps}^2,
\end{equation}
where $\lim \limits_{R\to \infty}\epsilon(R)=0$. With the change of variables $x=k_{N,\eps}y$, the definition (A4) of $\tilde{u}_N$ and $k_{i,\eps}$ and the asymptotic (A9) in Proposition \ref{prop:exhaust}, we get that
\begin{eqnarray}
&&\int \limits_{B_{R k_{N,\eps}}(0)\setminus B_{R^{-1} k_{N,\eps}}(0)}(x,\nabla b(x))\frac{|\ue|^{\crits-\pe}}{|x|^s}\, dx\\
&&=\mu_{N,\eps}^{2-\frac{n}{\crits-2}p_\eps}\int \limits_{B_{R }(0)\setminus B_{R^{-1}}(0)}\left(x,\frac{\nabla b(k_{N,\eps}x)}{k_{\eps,N}}\right))\frac{|\tilde{u}_{N,\eps}|^{\crits-\pe}}{|y|^s}\, dy.\label{ineq:18}
\end{eqnarray}
Since $\nabla b(0)=0$, we have that
\begin{equation}\label{ineq:19}
(x,\nabla b(x))=\partial_{ij}b(0)x^ix^j+o(|x|^2)\hbox{ as }x\to 0.
\end{equation}
Finally, with $r_\eps:=\sqrt{\mu_{\eps,N}}$, we have that
$$\int \limits_{B_{r_\eps}(0)}\frac{|x|^2}{|x|^s}|x|^{-\am(\crits-p_\eps)}\, dx=O\left(\mu_{N,\eps}^{1+\frac{\crits}{4}(\ap-\am)}\right).$$
Now, for $R>0$, with the upper-bound \eqref{eq:est:global}, we get that
\begin{eqnarray}
&&\left|~\int \limits_{\left(B_{r_\eps}(0)\setminus B_{R k_{N,\eps}}(0)\right)\cup B_{R^{-1} k_{N,\eps}}(0)}(x,\nabla b(x))\frac{|\ue|^{\crits}}{|x|^s}\, dx \right| \notag \\
&&\leq C \int \limits_{B_{r_\eps}(0)}\frac{|x|^2}{|x|^s}|x|^{-\am(\crits-p_\eps)}\, dx+C\sum_{i<N}\int \limits_{B_{\re}(0)}\frac{|x|^2B_{i,\eps}^{\crit-p_\eps}}{|x|^s}\, dx \notag\\
&&+C\int \limits_{\left(B_{\re}(0)\setminus B_{R k_{N,\eps}}(0)\right)\cup B_{R^{-1} k_{N,\eps}}(0)}\frac{|x|^2B_{N,\eps}^{\crit-p_\eps}}{|x|^s}\, dx.\label{ineq:20}
\end{eqnarray}
Using \eqref{ineq:19}, letting $\eps\to 0$ and then $R\to +\infty$ in \eqref{ineq:18} and \eqref{ineq:20} in the estimates above, we get \eqref{eq:gradb}. This is similar for \eqref{eq:gradb}. The estimates for the case $\dfrac{\crits}{2}(\ap-\am)\leq 2$ are obtained via a rough upper-bound. This ends Step \ref{step:p:grad}.\qed
\medskip

\subsection{Estimates of the $L^2$ term in the localized Pohozaev identity}

\begin{step}\label{step:p3} We have, as $\eps\to 0$,

\begin{itemize}
\item When  $\ap -\am >2-\theta$,
\begin{align}\label{cpct:L2term:sc}
 \int \limits_{  B_{ r_{\eps} }(0)\setminus B_{ k^{2}_{1,\eps} }(0)}  \left( \he(x) +\frac{ \left( \nabla \he, x \right)}{2} \right)\ue^{2}~ dx=   \mu_{N,\eps}^{2-\theta}  \left[ ~  \left(\frac{2-\theta}{2} \right) \frac{K_{h_0} }{t_{N}^{\frac{n-\theta }{\crits -2}}}  \int \limits_{\R^{n}} \frac{\tu_{N}^{2}}{|x|^{\theta}} ~dx  + o(1)\right].
\end{align}
where  $K_{h_0}$ is as in \eqref{hyp:h0}.\\
\item When $\ap -\am =2-\theta$,
\begin{align}\label{cpct:L2term:sc:log}
 \int \limits_{  B_{ r_{\eps} }(0)\setminus B_{ k^{2}_{1,\eps} }(0)}  \left( \he(x) +\frac{ \left( \nabla \he, x \right)}{2} \right)\ue^{2}~ dx=  O\left( \mu_{N,\eps}^{2-\theta}\ln\frac{1}{\mu_{\eps,N}}\right).
\end{align}
\item When $\ap-\am<2-\theta$,
\begin{align}\label{cpct:L2term:sc:bis}
 \int \limits_{  B_{ r_{\eps} }(0)\setminus B_{ k^{2}_{1,\eps} }(0)}  \left( \he(x) +\frac{ \left( \nabla \he, x \right)}{2} \right)\ue^{2}~ dx=   O\left(\mu_{N,\eps}^{\frac{2-\theta+(\ap-\am)}{2}}  \right).
\end{align}
\end{itemize}

\end{step}
\noindent{\it Proof of Step \ref{step:p3}:} Due to the hypothesis \eqref{hyp:h0} and \eqref{hyp:he}, we have that 
\begin{equation}
 \left| \he(x) +\frac{ \left( \nabla \he, x \right)}{2} \right|\leq C|x|^{-\theta}\hbox{ for all }x\in \Omega\setminus\{0\}\hbox{ and }\eps>0.
 \end{equation}
\noindent{\it Case 1: $\ap-\am<2-\theta$.}  We estimate roughly the integral with \eqref{eq:est:global}
\begin{eqnarray*}
&&\left|~ \int \limits_{  B_{ r_{\eps} }(0)\setminus B_{ k^{2}_{1,\eps} }(0)} |x|^{-\theta}\ue^{2}~ dx\right|\leq C\int \limits_{B_{ r_{\eps} }(0)\setminus B_{ k^{2}_{1,\eps} }(0)}  |x|^{-\theta} \left(~\sum_{i=1}^N\frac{\mei^{\ap-\am}}{ |x|^{2\ap}}+ \frac{1}{|x|^{2\am}} \right)\, dx\\
&&\leq C\mu_{N,\eps}^{\ap-\am}\int \limits_{  B_{ r_{\eps} }(0)\setminus B_{ k^{2}_{1,\eps} }(0)}  |x|^{-\theta-2\ap} \, dx +C\int \limits_{  B_{ r_{\eps}}(0)\setminus B_{ k^{2}_{1,\eps} }(0)}  |x|^{-\theta-2\am} \, dx\\
&&\leq C\mu_{N,\eps}^{\ap-\am}\int \limits_{0}^{r_\eps}r^{n-\theta-2\ap-1} \, dx +C\int \limits_0^{r_{\eps}}  r^{n-\theta-2\am-1} \, dx\\
&&\leq C\mu_{N,\eps}^{\ap-\am}\int \limits_{0}^{r_\eps}r^{2-\theta-(\ap-\am)-1} \, dr +C\int \limits_0^{r_{\eps}}  r^{\ap-\ap+2-\theta-1} \, dr\\
&&\leq   C \mu_{N,\eps}^{\frac{2-\theta+(\ap-\am)}{2}} 
\end{eqnarray*}
since $\ap-\am<2-\theta$ and $r_\eps=\sqrt{\mu_{N,\eps}}$. This proves \eqref{cpct:L2term:sc:bis}.

\medskip\noindent{\it Case 2: $\ap-\am=2-\theta$.} Here again since  $r_\eps=\sqrt{\mu_{N,\eps}}$, we estimate roughly the integral with \eqref{eq:est:global} to get
\begin{eqnarray*}
&&\left|~\int \limits_{B_{ r_{\eps} }(0)\setminus B_{ k^{2}_{1,\eps} }(0)} |x|^{-\theta}\ue^{2}~ dx~\right| \notag\\
&&\leq C\int \limits_{B_{ r_{\eps} }(0)\setminus B_{ k^{2}_{1,\eps} }(0)} |x|^{-\theta} \left(~\sum_{i=1}^N\frac{\mei^{\ap-\am}}{ \mei^{2(\ap-\am)}|x|^{2\am}+|x|^{2\ap}}+ \frac{1}{|x|^{2\am}} \right)\, dx\\
&&\leq C\sum_{i=1}^N\mei^{2-\theta}\int \limits_{ B_{ \mei^{-1}r_{\eps} }(0) }  \frac{1}{|x|^\theta\left(|x|^{2\am}+|x|^{2\ap}\right)} \, dx +
C\int \limits_{  B_{ r_{\eps} }(0)  }  |x|^{-\theta-2\am} \, dx \\
&&\leq C\sum_{i=1}^N\mei^{2-\theta}\ln\frac{r_\eps}{\mei}+C\mu_{N,\eps}^{2-\theta}\leq  C\sum_{i=1}^N\mei^{2-\theta}\ln\frac{1}{\mei}\leq C\mu_{N,\eps}^{2-\theta}\ln\frac{1}{\mu_{N,\eps}}.
\end{eqnarray*}
Since $\mei=o(\mu_{N,\eps})$ for $i<N$, this proves \eqref{cpct:L2term:sc:log}.

\medskip\noindent\noindent{\it Case 3: $\ap-\am>2-\theta$.} It follows from point \eqref{prop:bdd:blowupsol} of Theorem \ref{th:hopf} that the for all $1 \leq i \leq N$, the  function $\dfrac{\tu_{i}}{|x|^{\theta/2}} \in L^{2}(\rn)$ in this case.  For any $R,\rho >0$ we decompose  the  integral as
\begin{align*}
& \int \limits_{   B_{r_{\eps} }(0)\setminus B_{ k^{2}_{1,\eps} }(0)}\left( \he(x) +\frac{ \left( \nabla \he, x \right)}{2} \right) \ue^{2}   ~ dx= \notag\\
& \int \limits_{  B_{ r_{\eps} }(0) \setminus \overline{B}_{R k_{N,\eps}}(0) }  \left( \he(x) +\frac{ \left( \nabla \he, x \right)}{2} \right)\ue^{2}~ dx + \sum  \limits_{i=1}^{N}   \int \limits_{ B_{R k_{i,\eps}}(0)  \setminus \overline{B}_{\rho k_{i,\eps}}(0) }    \left( \he(x) +\frac{ \left( \nabla \he, x \right)}{2} \right)\ue^{2}~ dx   \notag\\
&+ \sum  \limits_{i=1}^{N-1}   \int \limits_{   B_{\rho k_{i+1,\eps}}(0)  \setminus \overline{B}_{R k_{i,\eps}}(0)  }   \left( \he(x) +\frac{ \left( \nabla \he, x \right)}{2} \right)\ue^{2}~ dx    +  
 \int \limits_{   B_{\rho k_{1,\eps}}(0)\setminus B_{ k^{2}_{1,\eps} }(0)}    \left( \he(x) +\frac{ \left( \nabla \he, x \right)}{2} \right)\ue^{2}~ dx. 
\end{align*}

\noindent
From the estimate \eqref{eq:est:global}, we  get as $\eps \to 0$
\begin{align*}
 &\mu_{N,\eps}^{-(2-\theta)} \int \limits_{   B_{ r_{\eps} }(0) \setminus \overline{B}_{R k_{N,\eps}}(0) } \left| \he(x) +\frac{ \left( \nabla \he, x \right)}{2} \right| \ue^{2} ~ dx  \notag\\
 \leq &~  C~ \mu_{N,\eps}^{-(2-\theta)}\int \limits_{    B_{ r_{\eps} }(0) \setminus \overline{B}_{R k_{N,\eps}}(0)  } \frac{1}{|x|^{\theta}} \left[ \frac{ \mu_{N,\eps}^{\ap-\am}}{|x|^{2\ap}} +   \frac{ 1}{|x|^{2\am}} \right]~dx \notag\\
\leq &~  C~\int \limits_{     B_{ \frac{r_{\eps}}{k_{N,\eps}}}(0)  \setminus \overline{B}_{R }(0)  } \frac{1}{|x|^{\theta}}   \frac{ 1}{|x|^{n+ (\ap-\am-2)} } ~dx \notag\\
+ &~  C~\int \limits_{     B_{1}(0)  \setminus \overline{B}_{\frac{R k_{N,\eps}}{r_{\eps}}}(0)  }   \frac{1}{|x|^{\theta}} \frac{  \mu_{N,\eps}^{ \frac{\ap-\am-(2-\theta)}{2}} }{|x|^{n- (\ap-\am+2)} } ~dx \notag\\
 &\leq C ~\left(R^{- (\ap-\am-(2-\theta))}+  \mu_{N,\eps}^{ \frac{\ap-\am-(2-\theta)}{2}} \right) .
\end{align*}
Therefore when  $\ap -\am >2-\theta$
\begin{align}{\label{cpct:L2term:1}}
\lim\limits_{R \to + \infty}  \lim\limits_{\eps \to 0} ~ \mu_{N,\eps}^{-(2-\theta)} \int \limits_{  B_{ r_{\eps} }(0) \setminus \overline{B}_{R k_{N,\eps}}(0)  }   \left( \he(x) +\frac{ \left( \nabla \he, x \right)}{2} \right) \ue^{2} ~ dx =0.
\end{align}

\noindent
Since in this case $\dfrac{\tui}{|x|^{\theta/2}} \in L^{2}(\rn)$ for any $1 \leq i\leq N$, it follows from Proposition \ref{prop:exhaust}  with  a change of variable $x=k_{N,\eps}y$, that 
\begin{align}{\label{cpct:L2term:2}}
\lim\limits_{R \to + \infty} \lim\limits_{\rho \to 0} \lim\limits_{\eps \to 0}~  \left(\mu_{i,\eps}^{-(2-\theta)} \int \limits_{  B_{R k_{i,\eps}}(0)  \setminus \overline{B}_{\rho k_{i,\eps}}(0)  }      \left( \he(x) +\frac{ \left( \nabla \he, x \right)}{2} \right)\ue^{2}~ dx \right)\notag\\
= \left(\frac{2-\theta}{2} \right) \frac{K_{h_0}}{t_{i}^{ \frac{n-\theta}{\crits-2}}} \int \limits_{\rn } \frac{\tu_{i}^{2 }}{|x|^{\theta}} ~ dx.
\end{align}

\noindent
Let $1 \leq i\leq N-1$. Using the pointwise estimates \eqref{eq:est:global}, for any $R, \rho>0$ and all $\eps>0$ we have  as $\eps \to 0$
\begin{align*}
 &\mu_{i+1,\eps}^{-(2-\theta)} \int \limits_{   B_{\rho k_{i+1,\eps}}(0)  \setminus \overline{B}_{R k_{i,\eps}}(0)  }  \left| \he(x) +\frac{ \left( \nabla \he, x \right)}{2} \right|  \ue^{2} ~ dx  \notag\\
  \leq &~  C~ \mu_{i+1,\eps}^{-(2-\theta)} \int \limits_{   B_{\rho k_{i+1,\eps}}(0)  \setminus \overline{B}_{R k_{i,\eps}}(0)  } \frac{1}{|x|^{\theta}}  \left[   \frac{ \mei^{\ap-\am}}{|x|^{2\ap}}+ \frac{ \mu_{i+1,\eps}^{-(\ap-\am) } }{ |x|^{2\am}}  \right]~dx \notag \\
\leq &~  C~\mu_{i+1,\eps}^{-(2-\theta)}  \int \limits_{      B_{ \frac{\rho k_{i+1,\eps}}{k_{i,\eps}}}(0)  \setminus \overline{B}_{R }(0)  }  \frac{1}{|x|^{\theta}}   \frac{ \mei^{2-\theta}}{|x|^{n+ (\ap-\am-2)} }~dx \notag\\
 &+ ~  C~ \mu_{i+1,\eps}^{-(2-\theta)} \int \limits_{   B_{\rho }(0)  \setminus \overline{B}_{ \frac{R k_{i,\eps}}{k_{i+1,\eps}}}(0)  }  \frac{1}{|x|^{\theta}}   \frac{ \mu_{i+1,\eps}^{2-\theta} }{|x|^{n- (\ap-\am+2)} } ~dx \notag\\
 &\leq C~\mu_{i+1,\eps}^{-(2-\theta)}\left( \mu_{i,\eps}^{2-\theta} R^{- (\ap-\am-(2-\theta))}  +  \mu_{i+1,\eps}^{2-\theta} \rho^{\ap-\am+(2-\theta)} \right).
\end{align*}
And so  as $\ap -\am >2$
\begin{align}{\label{cpct:L2term:3}}
\lim\limits_{R \to + \infty} \lim\limits_{\rho \to 0} \lim\limits_{\eps \to 0} ~ \mu_{i+1,\eps}^{-(2-\theta)} \int \limits_{    B_{\rho k_{i+1,\eps}}(0)  \setminus \overline{B}_{R k_{i,\eps}}(0)  } \left( \he(x) +\frac{ \left( \nabla \he, x \right)}{2} \right)\ue^{2} ~ dx =0.
\end{align}
\medskip

\noindent
Similarly  from the pointwise  estimates \eqref{eq:est:global} of  Theorem \ref{prop:fund:est}, we have as $\eps \to 0$
\begin{align*}
 &\mu_{1,\eps}^{-(2-\theta)} \int \limits_{B_{\rho k_{1,\eps}}(0)\setminus B_{ k^{2}_{1,\eps} }(0)}  \left| \he(x) +\frac{ \left( \nabla \he, x \right)}{2} \right|\ue^{2} ~ dx  \notag\\
 \leq &~  C~\mu_{1,\eps}^{-(2-\theta)}\int \limits_{B_{\rho k_{1,\eps}}(0)\setminus B_{ k^{2}_{1,\eps} }(0)}  \frac{1}{|x|^{\theta}} \frac{ \mu_{1,\eps}^{-(\ap-\am)}}{|x|^{2\am}} ~dx \notag \\
\leq& ~  C~\int \limits_{ B_{\rho }(0)   }   \frac{1}{|x|^{\theta}}  \frac{ 1}{|x|^{n- (\ap-\am+2)} } ~dx \notag\\
 &\leq C   ~\rho^{\ap-\am+(2-\theta)}.
\end{align*}
Therefore
\begin{align}{\label{cpct:L2term:4}}
 \lim\limits_{\rho \to 0} \lim\limits_{\eps \to 0} \mu_{1,\eps}^{-(2-\theta)} \int \limits_{B_{\rho k_{1,\eps}}(0)\setminus B_{ k^{2}_{1,\eps} }(0)} \left( \he(x) +\frac{ \left( \nabla \he, x \right)}{2} \right)\ue^{2} ~ dx =0.
\end{align}
\medskip

\noindent
From \eqref{cpct:L2term:1}, \eqref{cpct:L2term:2}, \eqref{cpct:L2term:3}, \eqref{cpct:L2term:4}  and Proposition \ref{prop:exhaust} we  then obtain \eqref{cpct:L2term:sc}.  This ends Step \ref{step:p3}.\qed \bigskip

\begin{step}\label{step:p3bis} Suppose that $u_0\equiv 0$. We claim that, as $\eps\to 0$
\begin{itemize}
\item When  $\ap -\am >2-\theta$,
\begin{align}\label{cpct:L2term:sc:2}
 \int \limits_{  B_{ \delta_0}(0)\setminus B_{ k^{2}_{1,\eps} }(0)}  \left( \he(x) +\frac{ \left( \nabla \he, x \right)}{2} \right)\ue^{2}~ dx=   \mu_{N,\eps}^{2-\theta}  \left[ ~  \left(\frac{2-\theta}{2} \right) \frac{K_{h_0} }{t_{N}^{\frac{n-\theta }{\crits -2}}}  \int \limits_{\R^{n}} \frac{\tu_{N}^{2}}{|x|^{\theta}} ~dx  + o(1)\right].
\end{align}
where  $K_{h_0}$ is as in \eqref{hyp:h0};\\
\item When $\ap -\am =2-\theta$,
\begin{align}\label{cpct:L2term:sc:log:2}
 \int \limits_{  B_{ \delta_0 }(0)\setminus B_{ k^{2}_{1,\eps} }(0)}  \left( \he(x) +\frac{ \left( \nabla \he, x \right)}{2} \right)\ue^{2}~ dx=  O\left( \mu_{N,\eps}^{2-\theta}\ln\frac{1}{\mu_{\eps,N}}\right).
\end{align}\\
\item When $\ap-\am<2-\theta$,
\begin{align}\label{cpct:L2term:sc:bis:2}
 \int \limits_{  B_{ \delta_0 }(0)\setminus B_{ k^{2}_{1,\eps} }(0)}  \left( \he(x) +\frac{ \left( \nabla \he, x \right)}{2} \right)\ue^{2}~ dx=   O\left(\mu_{N,\eps}^{\ap-\am}  \right).
\end{align}
\end{itemize}
\end{step}
\noindent{\it Proof of Step \ref{step:p3bis}:} It follows from \eqref{eq:est:global} that
\begin{eqnarray*}
&&\int \limits_{B_{\delta_0}(0)\setminus B_{r_\eps}(0)}\frac{\ue^2}{|x|^\theta}\, dx\leq C\sum_{i=1}^N\int \limits_{B_{\delta_0}(0)\setminus B_{r_\eps}(0)}\frac{1}{|x|^\theta}\frac{\mei^{\ap-\am}}{ \mei^{2(\ap-\am)}|x|^{2\am}+|x|^{2\ap}}\, dx\\
&&\leq C\sum_{i=1}^N\int \limits_{B_{\delta_0}(0)\setminus B_{r_\eps}(0)}\frac{\mei^{\ap-\am}}{ |x|^{2\ap+\theta}}\, dx\leq C\mu_{N,\eps}^{\ap-\am}\int \limits_{r_\eps}^{\delta_0} r^{2-\theta-(\ap-\am)-1}\, dr.
\end{eqnarray*}
Combining this estimate with \eqref{cpct:L2term:sc} and \eqref{cpct:L2term:sc:log} yield \eqref{cpct:L2term:sc:2} and \eqref{cpct:L2term:sc:log:2}: this proves Step \ref{step:p3bis} for $\ap-\am\geq 2-\theta$. For the case $\ap-\am<2-\theta$, the same estimate as above yields
\begin{eqnarray*}
&&\int \limits_{B_{\delta_0}(0)}\frac{\ue^2}{|x|^\theta}\, dx\leq C\mu_{\eps,N}^{\ap-\am}\int \limits_{0}^{\delta_0} r^{2-\theta-(\ap-\am)-1}\, dr\leq C\mu_{N,\eps}^{\ap-\am},
\end{eqnarray*}
which gives us \eqref{cpct:L2term:sc:bis:2}. This ends Step \ref{step:p3bis}.\qed \bigskip

 \subsection{Proof of the sharp blow-up rates when $u_0\equiv 0$ and $\ap-\am = 2-\theta$}
\begin{step}\label{prop:11} We  let $(\ue)$, $(\he)$, $(\pe)$ and $b$ be such that $(E_\eps)$, \eqref{hyp:he}, \eqref{lim:pe},  \eqref{hyp:b}  and \eqref{bnd:ue} holds. Assume that blow-up occurs as in \eqref{blowup}. Suppose $\ue>0$ for all $\eps >0$ and $u_0\equiv 0$.  We fix a family of parameters $(\lambda_\eps)_{\eps>0}\in (0,+\infty)$ such that
\begin{equation}\label{hyp:lambda}
\lim_{\eps\to 0}\lambda_\eps=0\hbox{ and }\lim_{\eps\to 0}\frac{\mu_{N,\eps}}{\lambda_\eps}=0.
\end{equation}
Then, for all $x\in\rnp$, $x\neq 0$, we have that
\begin{equation*}
\lim_{\eps\to 0}\frac{\lambda_\eps^{\ap}}{\mu_{N,\eps}^{\frac{\ap-\am}{2}}}\ue(\lambda_\eps x)=\frac{\mathcal{K}}{|x|^{\ap}},
\end{equation*}
where \begin{equation}\label{def:K}
\mathcal{K}:=\frac{L_{\gamma,\Omega}}{t_N^{\frac{\ap}{\crits-2}}}\int \limits_{ \rn} \frac{\tilde{u}_{N}(z)^{\crits-1}}{|z|^{s+\am}}\,dz>0
\end{equation}
and $L_{\gamma,\Omega}$ is defined in \eqref{est:G:4}. Moreover, this limit holds in $C^2_{loc}(\rnp)$.
\end{step}
\noindent{\it Proof of Step \ref{prop:11}:} We define
$$w_\eps(x):=\frac{\lambda_\eps^{\ap}}{\mu_{N,\eps}^{\frac{\ap-\am}{2}}}\ue(\lambda_\eps x)$$
for all $x\in \rnp\cap B_{\delta\lambda_\eps^{-1}}(0)$. From   $(E_{\eps})$ it follows that for all $\eps >0$, we have that
\begin{equation*}
\left\{\begin{array}{cl}
-\Delta  w_\eps - \frac{\gamma}{\left| x\right|^{2}}  w_\eps -\lambda_{\eps}^2~\he \circ (\lambda_{\eps} x) w_\eps =\Xi_{\eps} b(\lambda_{\eps}x)\frac{w_\eps^{\crits-1-\pe}}{\left|x\right|^{s}}& \hbox{ in }\rnp\cap B_{\delta\lambda_\eps^{-1}}(0)\\
 w_\eps>0 & \hbox{ in }\rnp\cap B_{\delta\lambda_\eps^{-1}}(0).
\end{array}\right.
\end{equation*}
With
$$\Xi_{\eps}:=\left(\frac{\mu_{N,\eps}^{\frac{\ap-\am}{2}}}{\lambda_\eps^{\ap}}\right)^{\crit-2-p_\eps}\lambda_\eps^{2-s} .$$
Since $\mu_{N,\eps}^{p_\eps}\to t_N>0$ (see (A9) of Proposition \ref{prop:exhaust}) and $$\ap(\crits-2)-(2-s)=(\crits-2)\frac{\ap-\am}{2},$$ then using the hypothesis \eqref{hyp:lambda}, we get that
$$\left(\frac{\mu_{N,\eps}^{\frac{\ap-\am}{2}}}{\lambda_\eps^{\ap}}\right)^{\crit-2-p_\eps}\lambda_\eps^{2-s}\leq C
\left(\frac{\mu_{N,\eps}}{\lambda_\eps}\right)^{\ap(\crits-2-p_\eps)-(2-s)}=o(1)$$
as $\eps\to 0$. Since $u_0\equiv 0$, it follows from the pointwise control \eqref{eq:est:global} that there exists $C>0$ such that $0<w_\eps(x)\leq C |x|^{-\ap}$ for all $x\in \rnp\cap B_{\delta\lambda_\eps^{-1}}(0)$. It then follows from standard elliptic theory that there exists $w\in C^2(\rnp)$ such that
\begin{equation}\label{lim:w:c2}
\lim_{\eps\to 0}w_\eps=w\hbox{ in }C^2_{loc}(\rnp)
\end{equation}
with
\begin{equation*}
\left\{\begin{array}{ll}
-\Delta w - \frac{\gamma}{|x|^2}  w=0&\hbox{ in }\rnp\\
0\leq w(x)\leq C |x|^{-\ap}&\hbox{ in }\rnp.
\end{array}\right.
\end{equation*}
From Proposition 11 in Ghoussoub-Robert \cite{gr3} it follows that there exists $\Lambda\geq 0$ such that $w(x)= \Lambda  |x|^{-\ap}$ for all $x\in\rnp$. We are left with proving that $\Lambda=\mathcal{K}$ defined in \eqref{def:K}. We fix $x\in \rnp$. Green's representation formula yields
\begin{eqnarray}
w_\eps(x)&=&\int \limits_\Omega\frac{\lambda_\eps^{\ap}}{\mu_{N,\eps}^{\frac{\ap-\am}{2}}}G_\eps( \lambda_\eps x,y)\frac{\ue(y)^{\crits-1-p_\eps}}{|y|^s}\,dy\nonumber\\
&=& \int \limits_{B_{R k_{N,\eps}}(0)\setminus B_{\delta k_{N,\eps}}(0)}+\int \limits_{\Omega\setminus(B_{R k_{N,\eps}}(0)\setminus B_{\delta k_{N,\eps}}(0))}\label{est:w:0}
\end{eqnarray}
Here $G_{\epsilon}$ is  the Green's function for the coercive operator  $-\Delta-\frac{\gamma}{|x|^{2}}-h_{\epsilon}$, for $\eps $ small,   on $\Omega$ with Dirichlet boundary condition. 
\medskip

\noindent{\it Step \ref{prop:11}.1:} We estimate the first term of the right-hand-side. A change of variable yields
\begin{eqnarray*}
&&\int \limits_{ B_{Rk_{N,\eps}}(0)\setminus B_{\delta k_{N,\eps}}(0)}\frac{\lambda_\eps^{\ap}}{\mu_{N,\eps}^{\frac{\ap-\am}{2}}}G_\eps( \lambda_\eps x,y)\frac{\ue(y)^{\crits-1-p_\eps}}{|y|^s}\,dy\\
&&=\Xi_\eps^{(1)}\int \limits_{B_{R}(0)\setminus B_{\delta}(0)}G_\eps(\lambda_\eps x,k_{N,\eps}z)\frac{\tilde{u}_{N,\eps}(z)^{\crits-1-p_\eps}}{|z|^s}(1+o(1))\,dz
\end{eqnarray*}
with
$$\Xi_\eps^{(1)}:=\frac{\lambda_\eps^{\ap}}{\mu_{N,\eps}^{\frac{\ap-\am}{2}}}k_{N,\eps}^{n-s}\mu_{N,\eps}^{-\frac{n-2}{2}(\crits-1-p_\eps)}$$
It follows from \eqref{est:G:4} that for any $z\in \rnp$, we have that
$$G_\eps( \lambda_\eps x,k_{N,\eps}z)=(L_{\gamma,\Omega}+o(1))\frac{1}{\lambda_\eps^{\ap}|x|^{\ap}}\cdot \frac{1}{k_{N,\eps}^{\am}|z|^{\am}},$$
and that the convergence is uniform with repect to $z\in  B_{R}(0)\setminus B_{\delta}(0)$. Plugging this estimate in the above equality, using that $k_{N,\eps}=\mu_{N,\eps}^{1-p_\eps/(\crits-2)}$, $\mu_{N,\eps}^{p_\eps}\to t_N>0$ and the convergence of $\tilde{u}_{N,\eps}$ to $\tilde{u}_N$ (see Proposition \ref{prop:exhaust}), we get that
\begin{eqnarray*}
&&\int \limits \limits_{ B_{R k_{N,\eps}}(0)\setminus B_{\delta k_{N,\eps}}(0)}\frac{\lambda_\eps^{\ap}}{\mu_{N,\eps}^{\frac{\ap-\am}{2}}}G_\eps( \lambda_\eps x,y)\frac{\ue(y)^{\crits-1-p_\eps}}{|y|^s}\,dy\\
&&=L_{\gamma,\Omega}\frac{1}{|x|^{\ap}}t_N^{-\frac{\ap}{\crits-2}}\int \limits_{ B_{R}(0)\setminus B_{\delta}(0)}\frac{1}{|z|^{\am}}\frac{\tilde{u}_{N}(z)^{\crits-1}}{|z|^s}\,dz+o(1)
\end{eqnarray*}
as $\eps\to 0$. Therefore,
\begin{eqnarray}
&&\lim \limits \limits_{R\to +\infty, \delta\to 0}\lim_{\eps\to 0}\int \limits_{ B_{R k_{N,\eps}}(0)\setminus B_{\delta k_{N,\eps}}(0)}\frac{\lambda_\eps^{\ap}}{\mu_{N,\eps}^{\frac{\ap-\am}{2}}}G_\eps( \lambda_\eps x,y)\frac{\ue(y)^{\crits-1}}{|y|^s}\,dy\nonumber\\
&&=\frac{\mathcal{K}}{|x|^{\ap}}\label{est:w:1}
\end{eqnarray}
where $\mathcal{K}$ is as in \eqref{def:K}. 

\medskip\noindent{\it Step  \ref{prop:11}.2:} With the control \eqref{est:G:up} on the Green's function and the pointwise control \eqref{eq:est:global} on $\ue$, we get that
\begin{eqnarray}
&&\int \limits_{\Omega\setminus (B_{R k_{N,\eps}}(0)\setminus B_{\delta k_{N,\eps}}(0)}\frac{\lambda_\eps^{\ap}}{\mu_{N,\eps}^{\frac{\ap-\am}{2}}}G_\eps( \lambda_\eps x,y)\frac{\ue(y)^{\crits-1-\pe}}{|y|^s}\,dy\nonumber\\
&&\leq \sum_{i=1}^{N-1}A_{i,\eps} +B_\eps(R)+C_\eps(\delta)\label{est:w:0:bis}
\end{eqnarray}
where
\begin{equation*}
A_{i,\eps}:=C\frac{\lambda_\eps^{\ap}}{\mu_{N,\eps}^{\frac{\ap-\am}{2}}}\int \limits_{B_{R_0}(0)}\frac{\ell_\epsilon (x, y)^{\am}}{|\lambda_\eps x-y|^{n-2}|y|^s}\left(\frac{\mu_{i,\eps}^{\frac{\ap-\am}{2}}}{\mu_{i,\eps}^{\ap-\am}|y|^{\am}+|y|^{\ap}}\right)^{\crits-1-\pe}\, dy
\end{equation*}
\begin{equation*}
B_\eps(R):=C\frac{\lambda_\eps^{\ap}}{\mu_{N,\eps}^{\frac{\ap-\am}{2}}}\int \limits_{B_{R_0}(0)\setminus B_{R k_{N,\eps}}(0)}\frac{\ell_\epsilon (x, y)^{\am} }{|\lambda_\eps x-y|^{n-2}}\frac{\mu_{N,\eps}^{\frac{\ap-\am}{2}(\crits-1-\pe)}}{|y|^{\ap(\crits-1-\pe)+s}}\, dy
\end{equation*}

\begin{equation*}
C_\eps(\delta):=C(x)\frac{\lambda_\eps^{\ap+2-n+\am}}{\mu_{N,\eps}^{\frac{\ap-\am}{2}\cdot(\crits-\pe)}}\int \limits_{ B_{\delta k_{N,\eps}}(0)}\frac{dy}{ |y|^{\am(\crits-1-\pe)+s+\am}}
\end{equation*}
where $\ell_\epsilon (x, y):=\frac{\max\{\lambda_\eps|x|,|y|\}}{\min\{\lambda_\eps|x|,|y|\}}$, 

\medskip\noindent{\it Step \ref{prop:11}.3.} We first estimate $C_\eps(\delta)$. Since $n>s+\crits\am$ (this is a consequence of $\am<(n-2)/2$), straightforward computations yield
$$C_\eps(\delta)\leq C(x)\delta^{\frac{\crits}{2}(\ap-\am)},$$
and therefore
\begin{equation}\label{est:w:2}
\lim_{\delta\to 0}\lim_{\eps\to 0}C_\eps(\delta)=0.
\end{equation}

\medskip\noindent{\it Step \ref{prop:11}.4.} We estimate $B_\eps(R)$. We split the integral as
$$B_\eps(R)=\int \limits_{R k_{\eps,N}<|y|<\frac{\lambda_\eps|x|}{2}}I_\eps(y)\, dy+\int \limits_{\frac{\lambda_\eps|x|}{2}<|y|<2\lambda_\eps|x|}I_\eps(y)\, dy+\int \limits_{|y|>2\lambda_\eps|x|}I_\eps(y)\, dy$$
where $I_\eps(y)$ is the integrand. Since 
$$n-(s+\ap(\crits-1)+\am)=-\frac{\crits-2}{2}(\ap-\am)<0,$$ straightforward computations yield
\begin{eqnarray*}
&&\int \limits_{R k_{N,\eps}<|y|<\frac{\lambda_\eps|x|}{2}}I_\eps(y)\, dy\\
&&\leq C(x)\frac{\lambda_\eps^{\ap+\am+2-n}}{\mu_{N,\eps}^{\frac{\ap-\am}{2}}} \int \limits_{R k_{N,\eps}<|y|<\frac{\lambda_\eps|x|}{2}}\frac{\mu_{N,\eps}^{\frac{\ap-\am}{2}(\crits-1-\pe)}}{|y|^{\ap(\crits-1-\pe)+s+\am-1}}\, dy\\
&&\leq C(x)R^{-\frac{\crits-2}{2}(\ap-\am)},
\end{eqnarray*}
For the next term, a change of variable yields
\begin{eqnarray*}
&&\int \limits_{\frac{\lambda_\eps|x|}{2}<|y|<2\lambda_\eps|x|}I_\eps(y)\, dy\\
&&\leq C(x)\frac{\lambda_\eps^{\ap}}{\mu_{N,\eps}^{\frac{\ap-\am}{2}}} \int \limits_{\frac{\lambda_\eps|x|}{2}<|y|<2\lambda_\eps|x|}|\lambda_\eps x-y|^{2-n}\frac{\mu_{N,\eps}^{\frac{\ap-\am}{2}(\crits-1)}}{|y|^{\ap(\crits-1)+s}}\, dy\\
&&\leq C(x) \left(\frac{\mu_{N,\eps}}{\lambda_\eps}\right)^{\frac{\crits-2}{2}(\ap-\am)}\int \limits_{\frac{|x|}{2}<|z|<2|x|}|x-z|^{2-n}\, dz=o(1)
\end{eqnarray*}
as $\eps\to 0$. Finally, since $\ap+\am=n-2$ and $n-s-\ap\crits=\frac{\crits}{2}(\ap-\am)$, we estimate the last term
\begin{eqnarray*}
&&\int \limits_{|y|>2\lambda_\eps|x|}I_\eps(y)\, dy\\
&&\leq C(x) \mu_{N,\eps}^{\frac{\crits-2}{2}(\ap-\am)}\lambda_\eps^{\ap-\am}\int \limits_{|y|>2\lambda_\eps|x|}\frac{|y|^{\am+1-n-s}\, dy}{|y|^{\ap(\crits-1)}}\\
&&\leq C(x) \left(\frac{\mu_{N,\eps}}{\lambda_\eps}\right)^{\frac{\crits-2}{2}(\ap-\am)}=o(1)
\end{eqnarray*}
as $\eps\to 0$. All these inequalities yield
\begin{equation}\label{est:w:3}
\lim_{R\to +\infty}\lim_{\eps\to 0}B_\eps(R)=0.
\end{equation}

\medskip\noindent{\it Step \ref{prop:11}.5.} We fix $i\in \{1,...,N-1\}$ and estimate $A_{i,\eps}$. As above, we split the integral as
$$A_{i,\eps}=\int \limits_{|y|<\frac{\lambda_\eps|x|}{2}}J_{i,\eps}(y)\, dy+\int \limits_{\frac{\lambda_\eps|x|}{2}<|y|<2\lambda_\eps|x|}J_{i,\eps}(y)\, dy+\int \limits_{|y|>2\lambda_\eps|x|}J_{i,\eps}(y)\, dy,$$
where $J_{i,\eps}$ is the integrand. Since $\mu_{i,\eps}\leq \mu_{N,\eps}$, as one checks, the second and the third integral of the right-hand-side are controled from above respectively by $\int \limits_{\frac{\lambda_\eps|x|}{2}<|y|<2\lambda_\eps|x|}I_\eps(y)\, dy$ and $\int \limits_{|y|>2\lambda_\eps|x|}I_\eps(y)\, dy$ that have been computed just above and go to $0$ as $\eps\to 0$. We are then left with the first term. With a change of variables, we have that
\begin{eqnarray*}
&&\int \limits_{|y|<\frac{\lambda_\eps|x|}{2}}J_{i,\eps}(y)\, dy\\
&&\leq  C(x)\frac{\lambda_\eps^{\ap+\am+2-n}}{\mu_{N,\eps}^{\frac{\ap-\am}{2}}} \\
&&\times\int \limits_{|y|<\frac{\lambda_\eps|x|}{2}}
\left(\frac{\mu_{i,\eps}^{\frac{\ap-\am}{2}}}{\mu_{i,\eps}^{\ap-\am}|y|^{\am}+|y|^{\ap}}\right)^{\crits-1-\pe}
\frac{dy}{|y|^{s+\am}}\\
&&\leq C(x) \frac{\mu_{i,\eps}^{n-s-\am-\frac{n-2}{2}(\crits-1-\pe)}}{\mu_{N,\eps}^{\frac{\ap-\am}{2}}}\\
&&\times\int \limits_{|z|<\frac{\lambda_\eps|x|}{2\mu_{i,\eps}}}\frac{1}{|z|^{\am+s}}\left(\frac{1}{|z|^{\am}+|z|^{\ap}}\right)^{\crits-1-\pe}\, dz\\
&&\leq C(x)\left(\frac{\mu_{i,\eps}}{\mu_{N,\eps}}\right)^{\frac{\ap-\am}{2}}
\end{eqnarray*}
since $n>s+\crits\am$ and $n<\am+s+(\crits-1)\ap$. Since $\mu_{i,\eps}=o(\mu_{N,\eps})$ as $\eps\to 0$, we get that
\begin{equation}\label{est:w:4}
\lim_{\eps\to 0}A_{i,\eps}=0.
\end{equation}

\medskip\noindent{\it Step \ref{prop:11}.6:} Plugging \eqref{est:w:1}, \eqref{est:w:2}, \eqref{est:w:3} and \eqref{est:w:4} into \eqref{est:w:0} and \eqref{est:w:0:bis} yields $\lim \limits_{\eps\to 0}w_\eps(x)=\dfrac{\mathcal{K}}{|x|^{\ap}}$ for all $x\in \rnp$. With \eqref{lim:w:c2}, we then get that $\Lambda={\mathcal{K}}$. This proves Step \ref{prop:11}.\\

\noindent Now we get the optimal asymptotic when $\ap-\am=2-\theta$:

\begin{step}\label{step:p13} We  let $(\ue)$, $(\he)$, $(\pe)$ and $b$ be such that $(E_\eps)$, \eqref{hyp:he}, \eqref{lim:pe},  \eqref{hyp:b}  and \eqref{bnd:ue} holds. Assume that blow-up occurs as in \eqref{blowup} and $\ap-\am=2-\theta$. Suppose $\ue>0$ for all $\eps >0$ and $u_0\equiv 0$. Then

\begin{equation}
 \int \limits_{B_{ \delta_0}(0)\setminus B_{ k^{2}_{1,\eps} }(0)}  \left( \he(x) +\frac{ \left( \nabla \he, x \right)}{2} \right)\ue^{2}~ dx= \left(\frac{2-\theta}{2}\mathcal{K}^2\omega_{n-1}K_{h_0}+o(1)\right)\mu_{N,\eps}^{2-\theta}\ln\frac{1}{\mu_{N,\eps}}.\label{eq:log}
\end{equation}
\end{step}
\noindent{\it Proof of Step \ref{step:p13}:} Note that it follows from \eqref{hyp:he} that
\begin{equation}\label{eq:grad}
\left( \he(x) +\frac{ \left( \nabla \he, x \right)}{2} \right)=\frac{(2-\theta)K_{h_0}+o(1)}{2|x|^\theta}\hbox{ as }x\to 0\hbox{ and }\eps\to 0.
\end{equation}

\medskip\noindent We define $\theta_\eps:=\dfrac{1}{\sqrt{|\ln\mu_{N,\eps}|}}$, $\alpha_\eps:=\mu_{N,\eps}^{\theta_\eps}$ and $\beta_\eps:=\mu_{N,\eps}^{1-\theta_\eps}$. As one checks, we have that
\begin{equation}
\left. \begin{array}{ccc}\label{ppty:a:b}
\mu_{N,\eps}=o(\beta_\eps) & \beta_\eps=o(\alpha_\eps) & \alpha_\eps=o(1) \\
\ln \dfrac{\alpha_\eps}{\beta_\eps}\simeq \ln\dfrac{1}{\mu_{N,\eps}} & \ln \dfrac{\beta_\eps}{\mu_{N,\eps}}=o\left(\ln\dfrac{1}{\mu_{N,\eps}} \right) & \ln\alpha_\eps=o(\ln\mu_{N,\eps})
\end{array}\right\}
\end{equation}
as $\eps\to 0$. It then follows from \eqref{eq:est:global} and the properties \eqref{ppty:a:b} that
\begin{equation}
\left. \begin{array}{l}
\ds \int\limits_{   B_{ \delta_0}(0)\setminus B_{\alpha_\eps}(0)  }\frac{\ue^2}{|x|^\theta}\, dx =O\left(\mu_{N,\eps}^{2-\theta} \ln\frac{1}{\alpha_\eps}\right)=o\left(\mu_{N,\eps}^{2-\theta} \ln\frac{1}{\mu_{N,\eps}}\right);\\
\ds \int \limits_{ B_{ \beta_\eps}(0)\setminus B_{ k^{2}_{1,\eps} }(0)}\frac{\ue^2}{|x|^\theta}\, dx  =o\left(\mu_{N,\eps}^{2-\theta} \ln\frac{1}{\mu_{N,\eps}}\right).
\end{array}\right\}\label{est:25}
\end{equation}
Since $\mu_{N,\eps}=o(\beta_\eps)$ and $\alpha_\eps=o(1)$ as $\eps\to 0$, it follows from Proposition \ref{prop:11} that
\begin{equation}\label{asymp:ue}
\lim_{\eps\to 0}\sup_{x\in  B_{\alpha_\eps}(0)\setminus B_{\beta_\eps}(0)}\left|\frac{|x|^{2\ap}\ue^2(x)}{\mu_{N,\eps}^{\ap-\am}}-\mathcal{K}^2\right|=0
\end{equation}
and therefore
\begin{eqnarray*}
&&\int \limits_{   B_{ \alpha_\eps}(0)\setminus B_{\beta_\eps}(0)  }\frac{\ue^2}{|x|^\theta}\, dx =(\mathcal{K}^2+o(1))\mu_{N,\eps}^{\ap-\am}\int \limits_{   B_{ \alpha_\eps}(0)\setminus B_{\beta_\eps}(0)  }\frac{dx}{|x|^{\theta+2\ap}}\\
&&=(\mathcal{K}^2\omega_{n-1}+o(1))\mu_{N,\eps}^{2-\theta}\ln\frac{1}{\mu_{N,\eps}}.
\end{eqnarray*}
Plugging this last estimate into \eqref{eq:grad} yields \eqref{eq:log}. This ends the proof of Step \ref{step:p13}.
\hfill \qed

\subsection{Estimates of the boundary terms}
\begin{step}\label{step:p5} We let $(\ue)$, $(\he)$, $(\pe)$ and $b(x)$ be such that $(E_\eps)$, \eqref{hyp:he}, \eqref{lim:pe},  \eqref{hyp:b}  and \eqref{bnd:ue} holds. We assume that blow-up occurs as in \eqref{blowup}. For any $\eps>0$, we define
\begin{align*}
\tve(x):= r_{\eps}^{\am}  \ue( r_{\eps} x) \qquad \hbox{for } x \in B_{2\delta_0 r_{\eps}^{-1}}(0)\setminus\{0\} ,
\end{align*}
with $r_{\eps}:= \sqrt{\mu_{N,\eps}}$. We claim that there exists  $\tv \in C^{1} ( \rn) $ such that
$$\lim \limits_{\eps \to 0}  \tve(x)= \tv~\hbox{ in } C^{1}_{loc} ( \rnp) $$
where $\tv$ is a solution of 
\begin{equation}\label{blowup:eqn:cpct}
\left\{ \begin{array}{llll}
 -\Delta \tv -\dfrac{\gamma}{|x|^{2}} \tv&=&0  & \hbox{ in } \rnp
 \end{array}\right.
\end{equation}
\end{step}
\noindent {\it Proof of Step \ref{step:p5}:} From $(E_{\eps})$ it follows that for all $\eps >0$, the rescaled functions $\tve $  weakly satisfies the equation 
\begin{align}{\label{eqn:cpct:1}}
-\Delta  \tve - \dfrac{\gamma}{|x|^2}  \tve -r_{\eps}^2~\he \circ(r_{\eps} x)~ \tve =r_{\eps}^{\vartheta+\pe \am} b(r_{\eps} x)\dfrac{| \tve |^{\crits-2-\pe} \tve}{|x|^s}.
\end{align}
with  $\vartheta:=(\crits-2)\dfrac{\ap-\am}{2}>0$ .
\medskip

\noindent
Using the pointwise estimates  \eqref{eq:est:global}  we obtain the bound, that as $\eps \to 0$ we have for $ x \in \rn$
\begin{align*}
|\tv_{\eps}(x)|  \leq & ~C~r_{\eps}^{\am} \sum_{i=1}^N\frac{\mei^{\frac{\ap-\am}{2}}}{ \mei^{\ap-\am}| x|^{\am}+|x|^{\ap}}  \notag\\
&~+~C~ r_{\eps}^{\am}\frac{\Vert |x|^{\am}u_0 \Vert|_{L^{\infty}(\Omega)}}{|x|^{\am}}  \notag\\
 \leq & ~C~\sum_{i=1}^N \frac{\left( \frac{\mei}{\mu_{N,\eps}} \right)^{\frac{\ap-\am}{2}}  }{ \left( \frac{\mei}{ \sqrt{\mu_{N,\eps}}} \right)^{\ap-\am} \left| \frac{x}{r_{\eps}} \right|^{\am}+ \left| \frac{x}{r_{\eps}} \right|^{\ap}}  \notag\\
&~+~C~  \frac{\Vert |x|^{\am}u_0 \Vert|_{L^{\infty}(\Omega)}   }{ \left| \frac{x}{r_{\eps}} \right|^{\am}}  \notag\\
 \leq & ~C \left(~\sum_{i=1}^N\frac{ \left( \frac{\mei}{\mu_{N,\eps}} \right)^{\frac{\ap-\am}{2}} }{ \left( \frac{\mei}{ \sqrt{\mu_{N,\eps}}} \right)^{\ap-\am}| x|^{\am}+|x|^{\ap}}+ \frac{\Vert |x|^{\am}u_0 \vert|_{L^{\infty}(\Omega)} }{|x|^{\am}} \right)\\
 & \leq~C ~ \left(\frac{1}{ |x|^{\ap}} + \frac{\Vert |x|^{\am}u_0 \vert|_{L^{\infty}(\Omega)}}{ |x|^{\am} }\right).
\end{align*} 
Then  passing to limits in the  equation \eqref{eqn:cpct:1}, standard elliptic theory  yields the existence of $\tilde{v}\in C^2(\rnp)$ such that $\tv_{\eps} \to \tv$ in $C^2_{loc}(\rnp)$ and  $\tv$ satisfies the equation: 
\begin{align*}
\left\{ \begin{array}{llll}
-\Delta \tv -\dfrac{\gamma}{|x|^{2}} \tv&=&0  \ \ & \hbox{ in } \rn\setminus\{0\}
\end{array}\right.
\end{align*}
and we have the following bound on $\tv$
$$|\tv(x)|  \leq C ~ \left(\frac{1}{ |x|^{\ap}} + \frac{\Vert |x|^{\am}u_0 \Vert|_{L^{\infty}(\Omega)}}{ |x|^{\am} } \right) \qquad \hbox{ for all  } x\in\rnp.$$ 
This ends the proof of Step \ref{step:p5}.\qed
\bigskip

\begin{step}\label{step:p6} We let $(\ue)$, $(\he)$, $(\pe)$ and $b(x)$ be such that $(E_\eps)$, \eqref{hyp:he}, \eqref{lim:pe},  \eqref{hyp:b}  and \eqref{bnd:ue} holds. We assume that blow-up occurs as in \eqref{blowup}. We claim that, as $\eps \to 0$,
\begin{align}\label{cpct:bdry:term1and2}
\int \limits_{ \partial B_{ r_{\eps} }(0) }  F_{\eps}(x)~d\sigma =&~\mu_{N,\eps}^{\frac{\ap-\am}{2}} \left(  \mathcal{F}_0 +o(1) \right) \notag\\
\hbox{ and }~\int \limits_{\partial  B_{ k^{2}_{1,\eps} }(0) }  F_{\eps}(x)~d\sigma=& ~~ o \left( \mu_{N,\eps}^{\frac{\ap-\am}{2}}  \right),
\end{align}
with
\begin{align}\label{def:F1}
\mathcal{F}_0 := \int \limits_{\partial  B_{1 }(0) } ( x, \nu) \left( \frac{|\nabla \tv|^2}{2}  -\frac{\gamma}{2}\frac{\tv^{2}}{|x|^{2}} \right) -  \left(x^i\partial_i \tv+\frac{n-2}{2} \tv \right)\partial_\nu \tv ~d \sigma 
\end{align}
\end{step}

\noindent{\it Proof of Step \ref{step:p6}:} We keep the notations of Step \ref{step:p5}.
With a change of variable and the definition of $\tve$, and $\vartheta:=(\crits-2)\dfrac{\ap-\am}{2}>0$, we get
\begin{align*}
& \int \limits_{ \partial B_{ r_{\eps} }(0) } F_{\eps}(x)~d\sigma= \notag\\
& ~ r_{\eps}^{\ap-\am} \int \limits_{\partial  B_{1  }(0) } ( x, \nu) \left( \dfrac{|\nabla \tve|^2}{2}  -\frac{\gamma}{2}\dfrac{\tve^{2}}{|x|^{2}} \right) -  \left(x^i\partial_i \tve+\dfrac{n-2}{2} \tve \right)\partial_{\nu} \tve ~d \sigma  \notag\\
& ~-r_{\eps}^{\ap-\am} \int \limits_{\partial  B_{1 }(0) } ( x, \nu) \left(r_{\eps}^{2}\frac{h_{\eps}(r_{\eps} x)}{2} \tve^{2}  -b(r_{\eps}x)\dfrac{ r_{\eps}^{\vartheta+\am\pe} }{\crits -p_{\eps} } 
\dfrac{|\tve|^{\crits-p_{\eps} }}{|x|^{s}}\right) d\sigma. 
\end{align*}
Then from the convergence result of Step \ref{step:p5}, we  then get \eqref{cpct:bdry:term1and2}. This ends Step \ref{step:p6}.\qed

\begin{step}\label{step:p7} We let $(\ue)$, $(\he)$, $(\pe)$ and $b(x)$ be such that $(E_\eps)$, \eqref{hyp:he}, \eqref{lim:pe},  \eqref{hyp:b}  and \eqref{bnd:ue} holds. Assume that blow-up occurs as in \eqref{blowup}. Suppose $u_0\equiv 0$. We define 
\begin{equation}\label{def:bue}
\bar{u}_\eps:=\dfrac{\ue}{\mu_{N,\eps}^{\frac{\ap-\am}{2}}}.
\end{equation}
We claim that  there exists $\bar{u}\in C^2(\overline{\Omega}\setminus\{0\})$ such that
\begin{equation}\label{cv:bue}
\lim_{\eps\to 0}\bar{u}_\eps=\bar{u}\hbox{ in }C^2_{loc}(\overline{\Omega}\setminus\{0\})\hbox{ with }\left\{\begin{array}{ll}
-\Delta \bar{u}-\left(\dfrac{\gamma}{|x|^2}+h_0\right)\bar{u}=0 &\hbox{ in }\Omega\setminus\{0\}\\
\bar{u}=0&\hbox{ in }\partial\Omega
\end{array}\right.
\end{equation}
\end{step}
\noindent{\it Proof of Step \ref{step:p7}:} Since $u_0\equiv 0$, it follows from \eqref{eq:est:global} that there exists $C>0$ such that 
\begin{equation}\label{control:bue}
|\bar{u}_\eps(x)|\leq C |x|^{-\ap}\hbox{ for all }x\in\Omega\setminus\{0\}\hbox{ and }\eps>0.
\end{equation}
Moreover, equation  $(E_\eps)$ rewrites
$$-\Delta \bar{u}_\eps-\left(\dfrac{\gamma}{|x|^2}+h_\eps\right)\bar{u}_\eps=\mu_{N,\eps}^{\frac{\ap-\am}{2}(\crits-2-p_\eps)} b(x)\dfrac{|\bar{u}_\eps|^{\crits-2-p_\eps}\bar{u}_\eps}{|x|^s}\;\hbox{ in }\Omega,
$$
and $\bar{u}_\eps=0$ on $\partial\Omega$. It then follows from standard elliptic theory that the claim holds. This ends Step \ref{step:p7}.\qed

\begin{step}\label{step:p8} We let $(\ue)$, $(\he)$, $(\pe)$ and $b(x)$ be such that $(E_\eps)$, \eqref{hyp:he}, \eqref{lim:pe},  \eqref{hyp:b}  and \eqref{bnd:ue} holds. Assume that blow-up occurs as in \eqref{blowup}. Suppose that $u_0\equiv 0$. We claim that
\begin{align}\label{cpct:bdry:term1and2:lala}
\int \limits_{\partial B_{ \delta_0 }(0) }  F_{\eps}(x)~d\sigma&=\left( \mathcal{F}_{\delta_{0}}+o(1)\right)\mu_{N,\eps}^{\ap-\am} \notag\\
\hbox{ and }~\int \limits_{\partial  B_{ k^{2}_{1,\eps} }(0) }  F_{\eps}(x)~d\sigma&= ~~ o \left( \mu_{N,\eps}^{\ap-\am}  \right),
\end{align}
where
\begin{equation}\label{def:Fd}
 \mathcal{F}_{\delta_{0}} := \int \limits_{  \partial B_{ \delta_0 }(0) } ( x, \nu) \left( \frac{|\nabla \bar{u}|^2}{2}  -\left(\frac{\gamma}{|x|^2}+h_0\right)\frac{\bar{u}^{2}}{2} \right) -  \left(x^i\partial_i \bar{u}+\frac{n-2}{2} \bar{u} \right)\partial_{\nu} \bar{u} ~d \sigma.
\end{equation}
\end{step}
\noindent{\it Proof of Step \ref{step:p8}:} 
With a change of variable, the definition of $\bar{u}_\eps$ and the convergence \eqref{cv:bue}, we get
\begin{eqnarray}
&& \int \limits_{ \partial B_{ \delta_0 }(0) } F_{\eps}(x)~d\sigma\\
&&=  \mu_{N,\eps}^{\ap-\am} \int \limits_{ \partial B_{ \delta_0 }(0) } ( x, \nu) \left( \frac{|\nabla \bar{u}_\eps|^2}{2}  -\left(\frac{\gamma}{|x|^{2}}+h_\eps\right)\frac{\bar{u}_\eps^{2}}{2}\right)~d \sigma  \notag \\
 && - \frac{\mu_{N,\eps}^{\frac{\ap-\am}{2}(\crit-\pe)}}{\crit- \pe} \int \limits_{ \partial B_{\delta_0 }(0)} ( x, \nu) b(x) \frac{|\bar{u}_\eps|^{\crit-\eps}\bar{u}_\eps}{|x|^2}~d \sigma  \notag \\
 &&- \int \limits_{ \partial B_{\delta_0 }(0)}  \left(x^i\partial_i \bar{u}_\eps+\frac{n-2}{2} \bar{u}_\eps \right)\partial_{\nu} \bar{u}_\eps ~d \sigma  \notag\\
&&=\mu_{N,\eps}^{\ap-\am} \left(   \mathcal{F}_{\delta_{0}}+o(1) \right).\label{est:Fe}
\end{eqnarray}
where $\mathcal{F}_{\delta_{0}}$ is as above. And directly from the bound \eqref{eq:est:global} we obtain
\begin{align*}
\int \limits_{\partial B_{ k^{2}_{1,\eps} }(0)}  F_{\eps}(x)~d\sigma= ~~ o \left( \mu_{N,\eps}^{\ap-\am}  \right)
\end{align*}
This ends Step \ref{step:p8}.\qed

\begin{step}\label{step:p9} We let $(\ue)$, $(\he)$, $(\pe)$ and $b(x)$ be such that $(E_\eps)$, \eqref{hyp:he}, \eqref{lim:pe},  \eqref{hyp:b}  and \eqref{bnd:ue} holds. Assume that blow-up occurs as in \eqref{blowup}. Suppose that $\ue>0$ for all $\eps>0$. We have $\mathcal{F}_0\geq 0$ and
$${\mathcal F}_0>0\;\iff \; u_0>0.$$
where $\mathcal{F}_0$ is as in \eqref{def:F1}.
\end{step}
\noindent{\it Proof of Step \ref{step:p9}:} We let $\tv$ be defined as in Step \ref{step:p5}. It then follows  from Step \ref{step:p5} that $\tv \geq 0$  and $\tv$ satisfies \eqref{blowup:eqn:cpct} and we have the following bound on $\tv$
\begin{equation}
|\tv(x)|  \leq C ~ \left(\frac{1}{ |x|^{\ap}} + \frac{\Vert |x|^{\am}u_0 \vert|_{L^{\infty}(\Omega)}}{ |x|^{\am} } \right) \qquad \hbox{ for all  } x\in\rnp.\label{est:tv}
\end{equation} 
Given $\alpha\in\rr$, we define $v_\alpha(x):=|x|^{-\alpha}$ for all $x\in \rnp$. Since $\tv\geq 0$, it follows from Proposition 6.4 in Ghoussoub-Robert \cite{gr3} that  there exists $A,B\geq 0$ such that
\begin{equation}\label{eq:tv:A:B}
\tv:=A v_{\ap}+B v_{\am}.
\end{equation}

\noindent{\it Step \ref{step:p9}.1:} We claim that $B=0$ when $u_0\equiv 0$.

\smallskip\noindent This is a direct consequence of controling \eqref{eq:tv:A:B} with \eqref{est:tv} when $u_0\equiv 0$ and letting $|x|\to\infty$.

\medskip\noindent{\it Step \ref{step:p9}.2:} We claim that $B>0$ when $u_0>0$.

\smallskip\noindent We prove the claim. We fix $x\in \rnp$.  Green's representation formula yields
\begin{eqnarray*}
\tv_\eps(x)&=& \int \limits_\Omega r_{\eps}^{\am} G_\eps(r_{\eps} x, y)\frac{\ue^{\crits-1-\pe}(y)}{|y|^s}\, dy.
\end{eqnarray*} 
We fix $\mathcal{D} \subset\subset  \Omega \setminus \{0\}$. Then there exists $\tilde{c}(\mathcal{D})>0$ such that $|y|\geq \tilde{c}(\mathcal{D})$ for all $y\in \mathcal{D}$. Moreover, the control \eqref{est:G:low} of the Green's function yields
\begin{eqnarray*}
\tv_\eps(x)&\geq & C \int \limits_{\mathcal{D}} r_{\eps}^{\am} \frac{1}{r_\eps^{\am}|x|^{\am}}|\tilde{c}(\mathcal{D})-r_\eps |x||^{-(n-2)}\frac{\ue^{\crits-1-\pe}(y)}{|y|^s}\, dy,
\end{eqnarray*}
and then, passing to the limit $\eps\to 0$,  we get that
\begin{eqnarray*}
\tv(x)&\geq &  \frac{C}{|x|^{\am}} \int \limits_{\mathcal{D}} \frac{u_0^{\crits-1}(y)}{|y|^s}\, dy,
\end{eqnarray*}
for all $x\in \rnp$. As one checks, this yields $B\geq C \int \limits_{\mathcal{D}} \dfrac{u_0^{\crits-1}(y)}{|y|^s}\, dy >0$ when $u_0>0$. This ends Step \ref{step:p9}.2.

\medskip\noindent{\it Step \ref{step:p9}.3:} We claim that $A>0$.

\smallskip\noindent The proof is similar to Step \ref{step:p9}.2. We fix $x\in \rnp $ and $\tilde{\mathcal{D}} \subset\subset \rnp$. Green's representation formula and the pointwise control \eqref{est:G:low} yields
\begin{eqnarray*}
\tv_\eps(x)&\geq & \int \limits_{k_{N,\eps}\tilde{\mathcal{D}} } r_{\eps}^{\am} G_\eps(r_{\eps} x, y)\frac{\ue^{\crits-1}(y)}{|y|^s}\, dy\\
&\geq & \tilde{c}(\tilde{\mathcal{D}} ) \int \limits_{\tilde{\mathcal{D}} } r_{\eps}^{\am} G_\eps(r_{\eps} x, k_{N,\eps} y) k_{N,\eps}^n \frac{\ue(k_{N,\eps} y)^{\crits-1}}{|k_{N,\eps} y|^s}\, dy\\
&\geq &  \tilde{c}(\tilde{\mathcal{D}} ) \int \limits_{\tilde{\mathcal{D}} } r_{\eps}^{\am} \left(\frac{r_{\eps}|x|}{k_{N,\eps}|y|}\right)^{\am}|r_\eps x-k_{N,\eps} y|^{2-n} k_{N,\eps}^{\frac{n-2}{2}}\frac{\tilde{u}_{\eps,N}(y)^{\crits-1}}{|y|^s}\, dy
\end{eqnarray*}
Since $r_\eps:=\sqrt{\mu_{N,\eps}}$, letting $\eps\to 0$, we get with the convergence (A4) of Proposition \ref{prop:exhaust} that
\begin{eqnarray*}
\tv_\eps(x)&\geq & \tilde{c}(\tilde{\mathcal{D}} ) \int \limits_{\tilde{\mathcal{D}} } r_{\eps}^{2\am-(n-2)} |x|^{\am} \left|x-\frac{k_{N,\eps}}{r_\eps} y \right|^{2-n} k_{N,\eps}^{\frac{n-2}{2}-\am}\frac{\tilde{u}_{\eps,N}(y)^{\crits-1}}{|y|^s}\, dy\\
&\geq & \frac{ \tilde{c}(\tilde{\mathcal{D}} )}{|x|^{\ap}}\int \limits_{\tilde{\mathcal{D}} } \frac{\tilde{u}_N(y)^{\crits-1}}{|y|^s}\, dy
\end{eqnarray*}
for all $x\in\rnp$. Therefore, as one checks, $A\geq c(\omega)\int \limits_{\tilde{\mathcal{D}} } \dfrac{\tilde{u}_N(y)^{\crits-1}}{|y|^s}\, dy>0$. This ends Step \ref{step:p9}.3.

\medskip\noindent{\it Step \ref{step:p9}.4:} We claim that
\begin{equation}\label{id:F:A:B}
\mathcal{F}_0=2\omega_{n-1}\left(\frac{(n-2)^2}{4}-\gamma\right)\cdot AB.
\end{equation}
We prove the claim. Definition \eqref{def:F1} reads \begin{align}\label{def:F0:2}
\mathcal{F}_0 := \int \limits_{ \partial  B_{1 }(0) } ( x, \nu) \left( \frac{|\nabla \tv|^2}{2}  -\frac{\gamma}{2}\frac{\tv^{2}}{|x|^{2}} \right) -  \left(x^i\partial_i \tv+\frac{n-2}{2} \tv \right)\partial_\nu \tv ~d \sigma 
\end{align}
For simplicity, we define the bilinear form:
\begin{eqnarray*}
{\mathcal H}_\delta(u,v)&=&\int \limits_{\partial  B_{\delta}(0) } \left[( x, \nu) \left( (\nabla u,\nabla v)  -\gamma \frac{uv}{|x|^{2}} \right) -  \left(x^i\partial_i u+\frac{n-2}{2} u \right)\partial_\nu v\right.\\
&&\left.-\left(x^i\partial_i v+\frac{n-2}{2} v \right)\partial_\nu u\right] ~d \sigma
\end{eqnarray*}
As one checks, using \eqref{eq:tv:A:B}
\begin{eqnarray}
\mathcal{F}_0&=&\frac{1}{2}{\mathcal H}_1(A v_{\ap}+B v_{\am},A v_{\ap}+B v_{\am})\nonumber\\
&=& \frac{A^2}{2}{\mathcal H}_1(v_{\ap},v_{\ap})+AB{\mathcal H}_1(v_{\ap},v_{\am})\nonumber\\
&&+\frac{B^2}{2}{\mathcal H}_1(v_{\am},v_{\am})\nonumber
\end{eqnarray}
In full generality, we compute ${\mathcal H}_\delta(v_\alpha,v_\beta)$ for all $\alpha,\beta\in\rr$ and all $\delta>0$. Consequently, straightforward computations yield
$$  \left(x^i\partial_i v_\alpha+\frac{n-2}{2} v_\alpha \right)\partial_\nu v_\beta=-\beta\left(\frac{n-2}{2}-\alpha\right)\frac{v_\alpha v_\beta}{|x|}$$
and
$$(x,\nu)\left((\nabla v_\alpha,\nabla v_\beta)-\frac{\gamma}{|x|^2}v_\alpha v_\beta\right)=(\alpha\beta-\gamma)\frac{v_\alpha v_\beta}{|x|}$$
and then
\begin{eqnarray*}
{\mathcal H}_\delta(v_\alpha,v_\beta) &=& \int \limits_{ \partial B_\delta(0)}\left(\alpha\beta-\gamma   +\beta\left(\frac{n-2}{2}-\alpha\right)+\alpha\left(\frac{n-2}{2}-\beta\right) \right)\frac{v_\alpha v_\beta}{|x|}\, d\sigma
\end{eqnarray*}
Plugging all these identities together yields
$${\mathcal H}_\delta(v_\alpha,v_\beta) =\left(-\alpha\beta-\gamma   +(\alpha+\beta)\frac{n-2}{2} \right)\delta^{-\alpha-\beta-1+n-1}\omega_{n-1}$$
Since $\ap,\am$ are solutions to $X^2-(n-2)X+\gamma=0$, we get that 
$${\mathcal H}_\delta(v_{\am},v_{\am})= {\mathcal H}_\delta(v_{\ap},v_{\ap})= 0.$$
 Since $\ap+\am=n-2$ and $\ap\am=\gamma$, we get that
$${\mathcal H}_\delta(v_{\am},v_{\ap})=2\omega_{n-1}\left(\frac{(n-2)^2}{4}-\gamma\right).$$
Plugging all these results together yields \eqref{id:F:A:B}. This ends Step \ref{step:p9}.4.

\smallskip\noindent These substeps end the proof of Step \ref{step:p9}.\qed

\begin{step}\label{step:p10} We let ($\ue)$, $(\he)$, $(\pe)$ and $b(x)$ be such that $(E_\eps)$, \eqref{hyp:he}, \eqref{lim:pe},  \eqref{hyp:b}  and \eqref{bnd:ue} holds. Assume that blow-up occurs as in \eqref{blowup}. If $\ue>0$ for all $\eps>0$ then $u_0\equiv 0$, when  $\ap-\am<4-2\theta$.
\end{step}
\noindent{\it Proof of Step \ref{step:p10}:} Since $\dfrac{\ap-\am}{2}<2-\theta$, plugging \eqref{cpct:2*term}, \eqref{cpct:L2term:sc}, \eqref{cpct:L2term:sc:log}, \eqref{cpct:L2term:sc:bis} and \eqref{cpct:bdry:term1and2} and \eqref{eq:gradb:log} and \eqref{eq:gradb:log}  into the Pohozaev identity \eqref{PohoId3}, we get as $\eps \to 0$
\begin{align}\label{id:poho:36}
  & \frac{p_{\eps}}{\crits} \left( \dfrac{n-s}{\crits}\right)  \left(\sum \limits_{i=1}^{N}  \frac{b(0)}{t_{i}^{ \frac{n-2}{\crits-2}}} \int \limits_{\rn }   \frac{|\tu_{i}|^{\crits }}{|x|^{s}} ~ dx +o(1)\right)
 =-\left(\mathcal{F}_0+o(1)\right)\mu_{N,\eps}^{\frac{\ap-\am}{2}},
\end{align}
where $\mathcal{F}_0$ is as in \eqref{def:F0:2}. Therefore $\mathcal{F}_0\leq 0$. Since $\ue>0$, it then follows from \eqref{id:F:A:B} of Step \ref{step:p9} that $u_0\equiv 0$. This proves Step \ref{step:p10}.\qed \bigskip

\section{Proof of the sharp blow-up rates}\label{pf blow-up rates}
We now prove the sharp blow-up rates claimed in  Propositions \ref{prop:rate:sc} and \ref{prop:rate:sc:2}. 

\subsection{The case $\ap-\am>2-\theta$.} It follows from the Pohozaev identity \eqref{PohoId3}, \eqref{cpct:2*term},  \eqref{eq:gradb}, \eqref{cpct:L2term:sc} and Step \ref{step:p6} when $u_0\not\equiv 0$, and the Pohozaev identity \eqref{PohoId3}, \eqref{cpct:2*term:2},  \eqref{eq:gradb},  \eqref{cpct:L2term:sc:2} and Step \ref{step:p8} when $u_0\equiv 0$ that
\begin{eqnarray}\label{id:poho:final}
&&  \mu_{N,\eps}^{2-\theta}  \left[ ~  \left(\frac{2-\theta}{2} \right) \frac{K_{h_0} }{t_{N}^{\frac{n-\theta }{\crits -2}}}  \int \limits_{\R^{n}} \frac{\tu_{N}^{2}}{|x|^{\theta}} ~dx  + o(1)\right].  \notag\\
 && +\frac{p_{\eps}}{\crits} \left( \frac{n-s}{\crits}\right) \left(\sum \limits_{i=1}^{N}  \frac{b(0)}{t_{i}^{ \frac{n-2}{\crits-2}}} \int \limits_{\rn }   \frac{|\tu_{i}|^{\crits }}{|x|^{s}} ~ dx +o(1)\right) \notag\\
&&+\frac{{\bf 1}_{\theta=0}}{\crits} \mu_{\eps,N}^2\left(\frac{\partial_{ij}b(0)}{t_N^{\frac{n}{\crits-2}}}\int \limits_{\rn}X^iX^j\frac{|\tilde{u}_N|^{\crits}}{|X|^s}\, dX\right) +o(\mu_{\eps,N}^{2}).\\
 &&=\left\{\begin{array}{ll}
  -\mu_{N,\eps}^{\frac{\ap-\am}{2}} \left(  \mathcal{F}_0 +o(1) \right) & \notag\\ 
  -\mu_{N,\eps}^{\ap-\am} \left(  \mathcal{F}_{\delta_0} +o(1) \right) &\hbox{if }u_0\equiv 0.
  \end{array}\right.
\end{eqnarray} 
We will use  the following lemma:
\begin{lemma}\label{lem:positive} Suppose we have  $\ap-\am>2$, $\gamma\geq 0$ and $\tilde{u}_N>0$. Then
\begin{align}\label{b:positive}
\frac{\partial_{ij}b(0)}{\crits}\frac{\int \limits_{\rn}X^iX^j\frac{|\tilde{u}_N|^{\crits}}{|X|^s}\, dX}{\int \limits_{\R^{n}} \frac{\tu_{N}^{2}}{|x|^{\theta}} ~dx }=\frac{\Delta b(0)}{b(0)}\frac{(n-2)[(\ap-\am)^2-4]}{4n(2n-2-s)}.
\end{align}
\end{lemma}
This lemma will be proved at the end of the present section.

\medskip\noindent{\it Case 1: $\ap-\am>2(2-\theta)$  or $\{u_0\equiv 0$ and $\ap-\am> 2-\theta\}$.} Then \eqref{id:poho:final} yields \eqref{rate:1} and  Proposition \ref{prop:rate:sc} is proved. When $\tue>0$, then $\tilde{u}_N\geq 0$ and $\tilde{u}_N\not\equiv 0$: it then follows from the strong comparison principle that $\tilde{u}_N>0$ and then Lemma \ref{lem:positive} holds. This yields Proposition \ref{prop:rate:sc:2} for $\ap-\am>2(2-\theta)$  or $\{u_0\equiv 0$ and $\ap-\am> 2-\theta\}$.

\medskip\noindent{\it Case 2: $\ue>0$ and $\ap-\am=2(2-\theta)$.} It then follows from \eqref{id:poho:final} that 
 \begin{eqnarray}
&&\lim_{\eps\to 0} \frac{p_{\eps}}{\mu_{N,\eps}^{2-\theta} }=-  \dfrac{ \left(\frac{2-\theta}{2} \right) \frac{K_{h_0} }{t_{N}^{\frac{n-\theta }{\crits -2}}}  \int \limits_{\R^{n}} \frac{\tu_{N}^{2}}{|x|^{\theta}} ~dx  +\frac{{\bf 1}_{\theta=0}}{\crits}\frac{\partial_{ij}b(0)}{ t_N^{\frac{n}{\crits-2}}}\int \limits_{\rn}X^iX^j\frac{|\tilde{u}_N|^{\crits}}{|X|^s}\, dX+ \mathcal{F}_0}{\frac{1}{\crits} \left( \frac{n-s}{\crits}\right) \left(\sum \limits_{i=1}^{N}  \frac{b(0)}{t_{i}^{ \frac{n-2}{\crits-2}}} \int \limits_{\rn }   \frac{|\tu_{i}|^{\crits }}{|x|^{s}} ~ dx\right)}.\notag 
\end{eqnarray} 
Then Step \ref{step:p9} and \eqref{b:positive} yield Proposition \ref{prop:rate:sc:2} when $\ap-\am=2(2-\theta)$.

\medskip\noindent{\it Case 3: $\ue>0$ and $2-\theta<\ap-\am<2(2-\theta)$.} It follows from Step \ref{step:p10} that $u_0\equiv 0$. Then  \eqref{id:poho:final} and \eqref{b:positive} yields Proposition \ref{prop:rate:sc} for $2-\theta<\ap-\am<2(2-\theta)$.

\subsection{The case $\ue>0$ and $\ap-\am=2-\theta$.} The Pohozaev identity \eqref{PohoId3}, \eqref{cpct:2*term:2},  \eqref{eq:gradb},  Step\ref{step:p13} and Step \ref{step:p8} yields  Proposition \ref{prop:rate:sc} for $\ap-\am=2-\theta$.

\subsection{The case $\ue>0$ and $\ap-\am< 2-\theta$}
\begin{step}\label{step:p11} We let $(\ue)$, $(\he)$, $(\pe)$ and $b(x)$ be such that $(E_\eps)$, \eqref{hyp:he}, \eqref{lim:pe},  \eqref{hyp:b}  and \eqref{bnd:ue} holds. We assume that blow-up occurs as in \eqref{blowup}. Suppose  $\ue>0$ for all $\eps >0$ and $\ap-\am<2-\theta$. Then \eqref{est:mass} holds, that is
\begin{equation}
\lim_{\eps\to 0}\frac{p_\eps}{\mu_{N,\eps}^{\ap-\am}}=-\dfrac{2\omega_{n-1}\crits^2 \left(\frac{(n-2)^2}{4}-\gamma\right)A^2}{(n-s)\sum \limits_{i=1}^{N}  \frac{b(0)}{t_{i}^{ \frac{n-2}{\crits-2}}} \int \limits_{\rn }   \frac{|\tu_{i}|^{\crits }}{|x|^{s}} ~ dx}\cdot m_{\gamma,h}(\Omega)
\label{est:mass:bis}
\end{equation}
for some $A>0$, where $m_{\gamma,h}(\Omega)$ is the mass.
\end{step}

\noindent{\it Proof of Step \ref{step:p11}:} It follows from Step \ref{step:p10} that $u_0\equiv 0$. 

\smallskip\noindent{\it Step \ref{step:p11}.1:} We now claim that
\begin{equation}
\frac{p_{\eps}}{\crits} \left( \frac{n-s}{\crits }\right)  \left(\sum \limits_{i=1}^{N}  \frac{b(0)}{t_{i}^{ \frac{n-2}{\crits-2}}} \int \limits_{\rn }   \frac{|\tu_{i}|^{\crits }}{|x|^{s}} ~ dx +o(1)\right) = \mu_{N,\eps}^{\ap-\am}\left(M_{\delta_0}+o(1)\right)\nonumber\end{equation}
where
\begin{eqnarray}\label{est:inf:1}
M_{\delta_0}&:=&-\int \limits_{B_{\delta_0 }(0) } \left( h_0(x) +\frac{ \left( \nabla h_0, x \right)}{2} \right)\bar{u}^{2}~ dx- \mathcal{F}_{\delta_0}
\end{eqnarray}
with $\mathcal{F}_{\delta_{0}}$ is as in \eqref{def:Fd} and $\bar{u}$ is as in \eqref{cv:bue}.

\smallskip\noindent Indeed  the convergence \eqref{def:bue}, \eqref{cv:bue}, \eqref{control:bue}  and $\ap-\am<2-\theta$ gives
\begin{eqnarray}\label{est:l2}
&&\int \limits_{B_{\delta_0 }(0) } \left( h_{\eps}(x) +\frac{ \left( \nabla h_{\eps}, x \right)}{2} \right)\ue^{2}~ dx\\
&&=\mu_{N,\eps}^{\ap-\am}\left(\int \limits_{  B_{\delta_0 }(0)  } \left( h_0(x) +\frac{ \left( \nabla h_0, x \right)}{2} \right)\bar{u}^{2}~ dx+o(1)\right)\nonumber
\end{eqnarray}
\smallskip\noindent Then from the Pohozaev identity \eqref{PohoId3}, \eqref{cpct:2*term:2},  \eqref{eq:gradb},  \eqref{cpct:L2term:sc:2} and Step \ref{step:p8} we obtain  \eqref{est:inf:1}, which  proves the claim and ends Step \ref{step:p11}.1.

\medskip\noindent Fix $\delta<\delta'$. Taking $U:= B_{\delta'}(0)\setminus B_\delta(0)$, $b\equiv0$ and $u=\bar{u}$ in \eqref{PohoId}, and using \eqref{cv:bue}, we get that $M_\delta$ is independent of the choice of $\delta>0$ small enough.

\medskip\noindent{\it Step \ref{step:p11}.2:} We claim that $\bar{u}>0$.

\medskip\noindent We prove this claim. Since $\bar{u}\geq 0$ is a solution to \eqref{cv:bue}, it is enough to prove that $\bar{u}\not\equiv 0$. We argue as in the proof of Step \ref{step:p9}. We fix $x\in \Omega\setminus\{0\}$. Green's identity and $\ue>0$  gives
\begin{eqnarray*}
&&\quad \bue(x)=\mu_{N,\eps}^{-(\ap-\am)/2}\int \limits_\Omega G_\eps(x,y)\frac{\ue(y)^{\crits-1-p_\eps}}{|y|^s}\, dy\\
&&\geq  \mu_{N,\eps}^{-(\ap-\am)/2}\int \limits_{\mathcal{D}_\eps} G_\eps(x,y)\frac{\ue(y)^{\crits-1-p_\eps}}{|y|^s}\, dy\\
&&\geq  C \mu_{N,\eps}^{n-s-(\ap-\am)/2}\int \limits_{\mathcal{D}} G_\eps(x,\mu_{N,\eps}y)\frac{\ue(\mu_{N,\eps}y)^{\crits-1-p_\eps}}{|y|^s}\, dy,
\end{eqnarray*}
where $\mathcal{D}_\epsilon:=  B_{2\mu_{N,\eps}}(0)\setminus B_{\mu_{N,\eps}}(0)$, $\mathcal{D}:= B_2(0)\setminus B_1(0)$. With the pointwise control \eqref{est:G:up}, we get
\begin{eqnarray*}
&&\quad \bue(x)\geq\\
&&  C(x) \int \limits_\mathcal{D} \left(\frac{|x|}{|y|}\right)^{\am}|x-\mu_{N,\eps}y|^{2-n}
\frac{u_{\eps,N}(y)^{\crits-1-p_\eps}}{|y|^s}\, dy
\end{eqnarray*}
where $u_{\eps,N}$ is as in Proposition \ref{prop:exhaust}. Letting $\eps\to 0$ and using the convergence (A4) of Proposition \ref{prop:exhaust}, we get that
$$\bar{u}(x)\geq C(x)\frac{1}{|x|^{\ap}} \hbox{ for all }x\in \Omega.$$
And then $\bar{u}>0$ in $\Omega$. This proves the claim and Step \ref{step:p11}.2.

\medskip\noindent We fix $r_0>0$ and $\eta\in C^\infty(\rn)$ such that $\eta(x)=1$ in $B_{r_0}(0)$ and $\eta(x)=0$ in $\rn\setminus B_{2r_0}(0)$. It then follows from \cites{gr4,gr5} that, for $r_0>0$ small enough, there exists $A>0$ and  $\beta\in H_{0}^1(\Omega)$ such that
$$\bar{u}(x)=A\left(\frac{\eta(x)}{|x|^{\ap}}+\beta(x)\right)\hbox{ for all }x\in \Omega$$
with
$$\beta(x)=m_{\gamma,h}(\Omega)\frac{\eta(x)  }{|x|^{\am}}+o\left(\frac{\eta(x)  }{|x|^{\am}}\right)$$
as $\eps\to 0$. Here, $m_{\gamma,h_{0}}(\Omega)$ is the  mass of $\Omega$ associated with the operator $-\Delta - \dfrac{\gamma}{|x|^2}-h_{0}(x)$, defined in Theorem \ref{def:mass}.. 

\medskip\noindent{\it Step \ref{step:p11}.3:}  We claim that
\begin{equation}\label{est:Md}
\lim_{\delta\to 0}M_\delta= - 2\omega_{n-1}\left(\frac{(n-2)^2}{4}-\gamma\right)A^2\cdot m_{\gamma,h}(\Omega)
\end{equation}
We prove the claim. Since $\bar{u}$ is a solution to \eqref{cv:bue}, it follows from standard elliptic theory that there exists $C>0$ such that $\bar{u}(x)+|x||\nabla\bar{u}(x)|\leq C|x|^{-\ap}$ for all $x\in B_{2\delta_0}(0)$. Therefore, since $\ap-\am<2-\theta$, we get that 
$$\lim_{\delta\to 0}\left(~\int \limits_{ B_\delta(0)}\bar{u}^2\, dx+\int \limits_{\partial B_\delta(0)}(x,\nu)\bar{u}^2\, d\sigma\right)=0.$$
And therefore,
\begin{equation*}
M_\delta=-\frac{A^2}{2}\bar{\mathcal H}_\delta(\bar{v}_{\ap}+\bar{v}_{\am},\bar{v}_{\ap}+\bar{v}_{\am} )+o(1)
\end{equation*}
as $\delta\to 0$, where
\begin{eqnarray*}
\bar{\mathcal H}_\delta(u,v)&:=&\int \limits_{ \partial B_{\delta_0 }(0)} \left[( x, \nu) \left( (\nabla u,\nabla v)  -\frac{\gamma}{|x|^2}uv \right) -  \left(x^i\partial_i u+\frac{n-2}{2} u \right)\partial_{\nu} v\right.\\
&&\left.- \left(x^i\partial_i v+\frac{n-2}{2} v \right)\partial_{\nu} u\right]\,d \sigma
\end{eqnarray*}
with
$$\bar{v}_{\ap}(x):=\frac{\eta(x) }{|x|^{\ap}}\hbox{ and }\bar{v}_{\am}(x)=\beta(x)\hbox{ for all }x\in \Omega.$$
We then get that
\begin{eqnarray*}
M_\delta&=&-\frac{A^2}{2}\bar{\mathcal H}_\delta(\bar{v}_{\ap},\bar{v}_{\ap} )-A^2\bar{\mathcal H}_\delta(\bar{v}_{\ap},\bar{v}_{\am} )\\
&&-\frac{A^2}{2}\bar{\mathcal H}_\delta(\bar{v}_{\am},\bar{v}_{\am} )+o(1)
\end{eqnarray*}
as $\delta\to 0$. Since $\ap-\am<2-\theta$ and $\ap+\am=n-2$, we get with a change of variable that as $\delta\to 0$, 
\begin{eqnarray*}
\bar{\mathcal H}_\delta(\bar{v}_{\ap},\bar{v}_{\ap} )&=&{\mathcal H}_\delta(v_{\ap},v_{\ap} )+O(\delta^{2-(\ap-\am)})\\
\bar{\mathcal H}_\delta(\bar{v}_{\ap},\bar{v}_{\am} )&=&m_{\gamma,h}(\Omega)\cdot {\mathcal H}_\delta(v_{\ap},v_{\am} )+o(\delta^{2})\\
\bar{\mathcal H}_\delta(\bar{v}_{\am},\bar{v}_{\am} )&=&O(\delta^{2+ (\ap-\am)}).
\end{eqnarray*}
Using the computations performed in the proof of Step \ref{step:p9}, we then get \eqref{est:Md}. This proves the claim and ends Step \ref{step:p11}.3.

\medskip\noindent{\it End of the proof of Step \ref{step:p11}:} Since $M_\delta$ is independent of $\delta$ small, we then get that $M_{\delta_0}=-2\omega_{n-1}\left(\dfrac{(n-2)^2}{4}-\gamma\right)A^2m_{\gamma,h}(\Omega)$. Putting this estimate in \eqref{est:inf:1}, we then get \eqref{est:mass:bis}. This end Step \ref{step:p11}.\qed

\medskip\noindent{\bf Proof of Lemma \ref{lem:positive}.} We assume that $\tilde{u}_N>0$. We define 
$$U(x):=\left(\frac{1}{|x|^{\frac{2-s}{n-2}\am}+|x|^{\frac{2-s}{n-2}\ap}}\right)^{\frac{n-2}{2-s}}$$ 
for all $x\in\rnp$. As one checks, we have that
\begin{eqnarray*}
-\Delta U-\frac{\gamma}{|x|^2}U&=&\frac{4(n-s)}{n-2}\left(\frac{(n-2)^2}{4}-\gamma\right)\frac{U^{\crits-1}}{|x|^s}\\
&=&\frac{(n-s)(\ap-\am)^2}{n-2}\frac{U^{\crits-1}}{|x|^s}.
\end{eqnarray*}
We define 
$$\mu^{\crits-2}=\frac{(n-s)(\ap-\am)^2}{(n-2)b(0)}.$$
so that $\mu U$ is a solution to \eqref{eq:tui}. 

\medskip\noindent We claim that there exists $\nu>0$ such that $\tilde{u}_N=\mu \nu^{-\frac{n-2}{2}}U(\nu^{-1}\cdot)$.

\smallskip\noindent We prove the claim. It follows from the convergence (A4) and the pointwise control \eqref{eq:est:global} that
$$\tilde{u}_N(x)\leq \frac{C}{|x|^{\am}+|x|^{\ap}}\hbox{ for }x\in\rnp.$$
Since $\tilde{u}_N$ solves \eqref{eq:tui}, one has that
$$-\hbox{div}(|x|^{-2\am}\nabla (|x|^{\am}\tilde{u}_N))=|x|^{-\am\crits-s}(|x|^{\am}\tilde{u}_N)^{\crits-1}\hbox{ in }\rnp$$
with $|x|^{\am}\tilde{u}_N\in C^2(\rnp)$ and $|x|^{\am}\tilde{u}_N\leq C(1+|x|^{\ap-\am})^{-1}$ in $\rnp$. Since $s>0$, it then follows from Chou-Chu \cite{chouchu} that $\tilde{u}_N$ is radially symmetrical when $\gamma\geq 0$.

\smallskip\noindent We  define $\varphi(t):=e^{-\frac{n-2}{2}t}\tilde{u}_N(e^{-t})$ for $t\in\rr$. We  get that $-\varphi"+\left(\frac{(n-2)^2}{4}-\gamma\right)\varphi=\varphi^{\crits-1}$ in $\rr$, $\varphi>0$ and $\lim \limits_{t\to\pm\infty}\varphi(t)=0$. The ODE yields the existence of $\mathcal{S}\in\rr$ such that 
\begin{equation}\label{int:prem}
\frac{(\varphi')^2}{2}+\left(\frac{(n-2)^2}{4}-\gamma\right)\frac{\varphi^2}{2}-\frac{\varphi^{\crits}}{\crits}=\mathcal{S}.
\end{equation}
Letting $t\to +\infty$ yields $\mathcal{S}=0$. Since $\lim_{t\to\pm\infty}\varphi(t)=0$, there exists $t_0\in\rr$ such that $\varphi'(t_0)=0$, and \eqref{int:prem} yields a unique possible value for $\varphi(t_0)$. Therefore, the theory of ODEs yields uniqueness of $\varphi$ up to translation. Since $\lambda U$ is also such a positive solution, we then get that there exists $\nu>0$ such that $\tilde{u}_N=\mu \nu^{-\frac{n-2}{2}}U(\nu^{-1}\cdot)$. This proves the claim.

\medskip\noindent With a change of variables and symmetry, we get that
\begin{eqnarray*}
&&\frac{\partial_{ij}b(0)}{ \crits}\frac{\int \limits_{\rn}X^iX^j\frac{|\tilde{u}_N|^{\crits}}{|X|^s}\, dX}{\int \limits_{\R^{n}} \tu_{N}^{2} ~dx }=\frac{\Delta b(0)}{n~ \crits}\frac{\int \limits_{\rn}|X|^2\frac{|\tilde{u}_N|^{\crits}}{|X|^s}\, dX}{\int  \limits_{\R^{n}}  \tu_{N}^{2}  ~dx }\\
&&=\frac{\Delta b(0)\mu^{\crits-2}(n-2)}{2n(n-s)}\frac{\int \limits_{\rn}|X|^2\frac{U^{\crits}}{|X|^s}\, dX}{\int \limits_{\R^{n}} U^2  ~dx }=\frac{\Delta b(0)}{b(0)}\frac{(\ap-\am)^2}{2n}\frac{\int \limits_{\rn}|X|^2\frac{U^{\crits}}{|X|^s}\, dX}{\int \limits_{\R^{n}} U^2~dx }
\end{eqnarray*}
We estimate this ratio as in Jaber \cite{jaber}. Passing to radial coordinates and taking the change of variable $t=r^{\frac{2-s}{n-2}(\ap-\am)}$, we get that
\begin{eqnarray*}
\frac{\int \limits_{\rn}|X|^2\frac{U^{\crits}}{|X|^s}\, dX}{\int \limits_{\R^{n}}  U^2  ~dx }=\frac{\int \limits_0^\infty\frac{r^{2-s+n-1}}{\left(r^{\frac{2-s}{n-2}\am}+r^{\frac{2-s}{n-2}\ap}\right)^{\frac{2(n-s)}{2-s}}}\, dr}{\int \limits_0^\infty\frac{r^{n-1}}{\left(r^{\frac{2-s}{n-2}\am}+r^{\frac{2-s}{n-2}\ap}\right)^{\frac{2(n-2)}{2-s}}}\, dr}=\frac{I_{\mathcal{Q}+2}^{\mathcal{P}+1}}{I_{\mathcal{Q}}^P}
\end{eqnarray*}
where $$I_{\mathcal{Q}}^\mathcal{P}=\int \limits_0^\infty\frac{t^\mathcal{P}}{(1+t)^\mathcal{Q}}\, dt\,;\, \mathcal{Q}=\frac{2(n-2)}{2-s}\hbox{ and }\mathcal{P}=\frac{(n-2)(n-2\am)}{(2-s)(\ap-\am)}-1$$
Integrarting by parts, see Jaber \cite{jaber}, we have that
$$I^{\mathcal{P}+1}_{\mathcal{Q}+1}=\frac{\mathcal{P}+1}{\mathcal{Q}}I_{\mathcal{Q}}^{\mathcal{P}}\hbox{ and }I_{\mathcal{Q}+1}^{\mathcal{P}}=\frac{\mathcal{Q}-\mathcal{P}-1}{\mathcal{Q}}I_{\mathcal{Q}}^{\mathcal{P}}.$$
We then get that
\begin{eqnarray*}
\frac{\int \limits_{\rn}|X|^2\frac{U^{\crits}}{|X|^s}\, dX}{\int \limits_{\R^{n}}  U^2  ~dx }=\frac{(\mathcal{P}+1)(\mathcal{Q}-\mathcal{P}-1)}{\mathcal{Q}(\mathcal{Q}+1)}=\frac{(n-2)(\ap-\am+2)(\ap-\am-2)}{2(\ap-\am)^2(2n-2-s)},
\end{eqnarray*}
so that
\begin{eqnarray*}
&&\frac{\partial_{ij}b(0)}{ \crits}\frac{\int \limits_{\rn}X^iX^j\frac{|\tilde{u}_N|^{\crits}}{|X|^s}\, dX}{\int \limits_{\R^{n}}  \tu_{N}^{2}  ~dx }=\frac{\Delta b(0)}{b(0)}\frac{(n-2)((\ap-\am)^2-4)}{4n(2n-2-s)}\hbox{ when }\tilde{u}_N>0.
\end{eqnarray*}
This completes Lemma \ref{lem:positive}.
\bigskip

\section{Blow-up rates when $b\not\in C^2(\overline{\Omega})$}
The function $b$ that arises from the hyperbolic model (see Lemma \ref{Lemma:relation between Hyp and Euc solution on bounded domain}) is not in $C^2$ when $n=3,4$, and therefore, the asymptotic rates of Propositions \ref{prop:rate:sc} and \ref{prop:rate:sc:2} do not apply. However, due to the behavior of $(x,\nabla b(x))$ in \eqref{est:nablab:3} and  \eqref{est:nablab:4}, we are in position to get optimal rates.

\subsection{The case of behavior like \eqref{est:nablab:4}} When the expressions make sense, we define
\begin{align}
\mathcal{C}_{n,s}:=\frac{\left(\frac{2-\theta}{2} \right) \frac{1}{t_{N}^{\frac{n-\theta }{\crits -2}}}  \int \limits_{\R^{n}} \frac{\tu_{N}^{2}}{|x|^{\theta}} ~dx }{\frac{n-s}{\crits^2}  \sum \limits_{i=1}^{N}  \frac{1}{t_{i}^{ \frac{n-2}{\crits-2}}} \int \limits_{\rn }   \frac{|\tu_{i}|^{\crits }}{|x|^{s}} ~ dx }\hbox{ and }\mathcal{D}_{n,s}:=\frac{\left(\frac{1}{\crits} \right) \frac{1}{t_{N}^{\frac{n }{\crits-2}}}  \int \limits_{\R^{n}} |x|^2\frac{|\tu_{N}|^{\crits}}{|x|^{s}} ~dx }{\frac{n-s}{\crits^2}  \sum \limits_{i=1}^{N}  \frac{1}{t_{i}^{ \frac{n-2}{\crits-2}}} \int \limits_{\rn }   \frac{|\tu_{i}|^{\crits }}{|x|^{s}} ~ dx }
\end{align}
The proof of compactness rely on the following two key propositions.
\begin{proposition}[$b\in C^1\setminus C^2$]\label{prop:rate:sc:4}
Let $\Omega$ be a smooth bounded domain of $\rn$, $n\geq 3$,  such that $0\in \Omega $ and  assume that  $0< s < 2$,   $\gamma<\frac{(n-2)^2}{4}$. Let $(\ue)$, $(\he)$, $(\pe)$ and $b$ be such that $(E_\eps)$, \eqref{hyp:he}, \eqref{lim:pe},  and \eqref{bnd:ue} holds. We assume that $b\in C^1(\overline{\Omega})$ and that there exists $C_1\in\rr$ such that
\begin{equation*}
(x,\nabla b(x))=C_1|x|^2\ln\frac{1}{|x|}+O(|x|^2)\hbox{ as }x\to 0.
\end{equation*}
Assume that blow-up occurs, that is
\begin{equation*}
\lim_{\eps\to 0}\Vert |x|^{\tau} \ue\Vert_{L^\infty(\Omega)}=+\infty ~\hbox{ for some }  ~ \am<\tau<\frac{n-2}{2} .
\end{equation*} 
Consider the $\mu_{1,\eps},...,\mu_{N,\eps}$  and $t_{1},...,t_{N}$  from Proposition \ref{prop:exhaust}. Then, we have the following blow-up rates:
\begin{equation*}\label{rate:1:4}
\lim \limits_{\eps \to 0} \frac{\pe}{\mu_{N,\eps}^{2-\theta}}  =-  \frac{\mathcal{C}_{n,s}}{b(0)}\cdot K_{h_0}\hbox{ if }\theta>0\hbox{ and }\left\{\begin{array}{c}
\hbox{either }\ap-\am>  4-2\theta\\
\hbox{or }\{\ap-\am>2-\theta\hbox{ and }u_0\equiv 0\}
\end{array}\right. 
\end{equation*}
\begin{equation*}\label{rate:1:4:bis}
\lim \limits_{\eps \to 0} \frac{\pe}{\mu_{N,\eps}^{2}\ln\frac{1}{\mu_{N,\eps}}}  =-  \frac{\mathcal{D}_{n,s}}{b(0)} \cdot C_1 \hbox{ if } \theta=0\hbox{ and }\left\{\begin{array}{c}
\hbox{either }\ap-\am\geq 4\\
\hbox{or }\{\ap-\am>2\hbox{ and }u_0\equiv 0\}
\end{array}\right. 
. 
\end{equation*}

\end{proposition}

\begin{proposition}[$b\in C^1\setminus C^2$, \textit{The positive case}]\label{prop:rate:sc:2:4}
Let $\Omega$ be a smooth bounded domain of $\rn$, $n\geq 3$,  such that $0\in \Omega $ and  assume that  $0< s < 2$,   $0 \leq \gamma<\frac{(n-2)^2}{4}$.  
Let $(\ue)$, $(\he)$, $(\pe)$ and $b$ be as in Proposition \ref{prop:rate:sc:4}. In particular
\begin{equation*}
(x,\nabla b(x))=C_1|x|^2\ln\frac{1}{|x|}+O(|x|^2)\hbox{ as }x\to 0.
\end{equation*}
Assume that blow-up occurs as in \eqref{blowup}. 
Consider  $\mu_{1,\eps},...,\mu_{N,\eps}$  and $t_{1},...,t_{N}$  from Proposition \ref{prop:exhaust}. Suppose in addition that
\begin{align*}
\ue>0 \qquad \hbox{for all } \eps >0.
\end{align*}
Then,  we have the following blow-up rates:

\begin{enumerate}
\item
When $\ap-\am =4-2\theta$ and $\theta>0$, we have 
\begin{align*}
&\lim \limits_{\eps \to 0} \frac{\pe}{\mu_{N,\eps}^{2-\theta}}  =-  \frac{\mathcal{C}_{n,s} }{b(0)}\left(K_{h_0}+\tilde{K}\right), 
\end{align*}
for some $\tilde{K}\geq 0$ such that $\tilde{K}> 0$ iff $u_0>0$. For $\theta=0$, see Proposition \ref{prop:rate:sc:4}.\\
\item
When $2-\theta<\ap-\am<4-2\theta$, we have $u_0\equiv 0$ and  Proposition \ref{prop:rate:sc:4} applies.

\item
When $\ap-\am=2-\theta$, we have $u_0\equiv 0$ and  
\begin{align}\label{asymp:crit:log:4}
\lim \limits_{\eps \to 0} \frac{\pe}{\mu_{N,\eps}^{2-\theta}\ln\frac{1}{\mu_{N,\eps}}}~=-\frac{\mathcal{C}'_{n,s}}{b(0)} \cdot K_{h_0}- {\bf 1}_{\theta=0}\frac{\mathcal{D}_{n,s}}{b(0)} \cdot C_1
\end{align}
where
$$\mathcal{C}'_{n,s}  :=\frac{\left(\frac{2-\theta}{2} \right)\mathcal{K}^{2} \omega_{n-1} }{ \frac{n-s}{\crits^2} \sum \limits_{i=1}^{N}  \frac{1}{t_{i}^{ \frac{n-2}{\crits-2}}} \int \limits_{\rn }   \frac{|\tu_{i}|^{\crits }}{|x|^{s}} ~ dx }~, $$
with $\mathcal{K}$ as defined in \eqref{def:K}. 

\item When $\ap-\am<2-\theta$, then $u_0\equiv 0$ and

\begin{align}\label{est:mass:4}
\lim_{\eps\to 0}\frac{p_\eps}{\mu_{N,\eps}^{\ap-\am}}=-\chi \cdot m_{\gamma,h_{0}}(\Omega)\quad \hbox{ if }\ap-\am<2-\theta,
\end{align}
where $\chi>0$ is a constant and $m_{\gamma,h_{0}}(\Omega)$ is the  mass of $\Omega$ associated with the operator $-\Delta - \dfrac{\gamma}{|x|^2}-h_{0}(x)$, defined in Theorem \ref{def:mass}.
%

\end{enumerate}
\end{proposition}
The proof of the propositions goes as the proof of Propositions \ref{prop:rate:sc} and \ref{prop:rate:sc:2}, with the use of the pointwise control of Lemma \ref{est:Bp}.

\subsection{The case of behavior like \eqref{est:nablab:3}} In this case, when writing the Pohozaev identity there are several terms to compare. For the sake of simplicity, we only deal here with the case $\theta=1$ that corresponds to the hyperbolic case when $n=3$. When the expression makes sense, we define
\begin{align}
\mathcal{E}_{n,s}:=\frac{\left(\frac{1}{2} \right) \frac{1}{t_{N}^{\frac{n-1 }{\crits -2}}}  \int \limits_{\R^{n}} \frac{\tu_{N}^{2}}{|x|} ~dx }{\frac{n-s}{\crits^2}  \sum \limits_{i=1}^{N}  \frac{1}{t_{i}^{ \frac{n-2}{\crits-2}}} \int \limits_{\rn }   \frac{|\tu_{i}|^{\crits }}{|x|^{s}} ~ dx }
\end{align}
The proof of compactness rely on the following two key propositions.
\begin{proposition}[$b\in C^0\setminus C^1$]\label{prop:rate:sc:3}
Let $\Omega$ be a smooth bounded domain of $\rn$, $n\geq 3$,  such that $0\in \Omega $ and  assume that  $0< s < 2$,   $\gamma<\frac{(n-2)^2}{4}$. Let $(\ue)$, $(\he)$, $(\pe)$ and $b$ be such that $(E_\eps)$, \eqref{hyp:he}, \eqref{lim:pe},  and \eqref{bnd:ue} holds. We assume that $\theta=1$, that $b\in C^{0,1}(\overline{\Omega})$ and that there exists $C_2\in\rr$ such that
\begin{equation*}
(x,\nabla b(x))=C_2|x|+O(|x|^2)\hbox{ as }x\to 0.
\end{equation*}
Assume that blow-up occurs, that is
\begin{equation*}
\lim_{\eps\to 0}\Vert |x|^{\tau} \ue\Vert_{L^\infty(\Omega)}=+\infty ~\hbox{ for some }  ~ \am<\tau<\frac{n-2}{2} .
\end{equation*} 
Consider the $\mu_{1,\eps},...,\mu_{N,\eps}$  and $t_{1},...,t_{N}$  from Proposition \ref{prop:exhaust}. Assume that
$$\{\ap-\am>2\}\hbox{ or }\{\ap-\am>1\hbox{ and }u_0\equiv 0\}.$$
Then, we have the following blow-up rates:
\begin{equation*}\label{rate:1:3}
\lim \limits_{\eps \to 0} \frac{\pe}{\mu_{N,\eps}}  =-  \frac{\mathcal{E}_{n,s}}{b(0)}\cdot \left(K_{h_0}+\frac{2}{\crits}\frac{\int \limits_{\R^{n}} |x|\frac{|\tu_{N}|^{\crits}}{|x|^s} ~dx}{\int \limits_{\R^{n}} \frac{\tu_{N}^{2}}{|x|} ~dx}\cdot C_2\right)
\end{equation*}
\end{proposition}

\begin{proposition}[$b\in C^0\setminus C^1$, \textit{The positive case for $\theta=1$}]\label{prop:rate:sc:2:3}
Let $\Omega$ be a smooth bounded domain of $\rn$, $n\geq 3$,  such that $0\in \Omega $ and  assume that  $0< s < 2$,   $0 \leq \gamma<\frac{(n-2)^2}{4}$.  
Let $(\ue)$, $(\he)$, $(\pe)$ and $b$ be as in Proposition \ref{prop:rate:sc:3}. In particular $\theta=1$ and
\begin{equation*}
(x,\nabla b(x))=C_2|x|+O(|x|^2)\hbox{ as }x\to 0.
\end{equation*}
Assume that blow-up occurs as in \eqref{blowup}. 
Consider  $\mu_{1,\eps},...,\mu_{N,\eps}$  and $t_{1},...,t_{N}$  from Proposition \ref{prop:exhaust}. Suppose in addition that
\begin{align*}
\ue>0  \hbox{ for all } \eps >0\hbox{ and }C_2>0
\end{align*}
Then,  we have the following blow-up rates:

\begin{enumerate}
\item
When $\ap-\am =2$, we have that

\begin{equation*} 
\lim \limits_{\eps \to 0} \frac{\pe}{\mu_{N,\eps}}  =-  \frac{\mathcal{E}_{n,s}}{b(0)}\cdot \left(K_{h_0}+\frac{(n-2)^2}{2(2n+2-s)}\left((n-2)^2-4\gamma-1\right)
\cdot \frac{C_2}{b(0)}+\tilde{K}\right)
\end{equation*}
for some $\tilde{K}\geq 0$ such that $\tilde{K}> 0$ iff $u_0>0$. 
\item
When $1<\ap-\am<2$, we have $u_0\equiv 0$ and Proposition \ref{prop:rate:sc:3} applies.

\item
When $\ap-\am=1$, we have $u_0\equiv 0$ and  
\begin{align*}
\lim \limits_{\eps \to 0} \frac{\pe}{\mu_{N,\eps}\ln\frac{1}{\mu_{N,\eps}}}~=-\frac{\mathcal{E}'_{n,s}}{b(0)} \cdot K_{h_0}
\end{align*}
where
$$\mathcal{E}'_{n,s}  :=\frac{\left(\frac{1}{2} \right)\mathcal{K}^{2} \omega_{n-1} }{ \frac{n-s}{\crits^2} \sum \limits_{i=1}^{N}  \frac{1}{t_{i}^{ \frac{n-2}{\crits-2}}} \int \limits_{\rn }   \frac{|\tu_{i}|^{\crits }}{|x|^{s}} ~ dx }~, $$
with $\mathcal{K}$ as defined in \eqref{def:K}. 

\item When $\ap-\am<1$, then $u_0\equiv 0$ and

\begin{align}\label{est:mass:3}
\lim_{\eps\to 0}\frac{p_\eps}{\mu_{N,\eps}^{\ap-\am}}=-\chi \cdot m_{\gamma,h_{0}}(\Omega)\end{align}
where $\chi>0$ is a constant and $m_{\gamma,h_{0}}(\Omega)$ is the  mass of $\Omega$ associated with the operator $-\Delta - \dfrac{\gamma}{|x|^2}-h_{0}(x)$, defined in Theorem \ref{def:mass}.
%

\end{enumerate}
\end{proposition}
The proof of the propositions goes as the proof of Propositions \ref{prop:rate:sc} and \ref{prop:rate:sc:2},  with the use of the pointwise control of Lemma \ref{est:Bp}. We have also used that when $\tu_N>0$, then
$$\frac{\int \limits_{\R^{n}} |x|\frac{|\tu_{N}|^{\crits}}{|x|^s} ~dx}{\int \limits_{\R^{n}} \frac{\tu_{N}^{2}}{|x|} ~dx}=\frac{(n-2)(n-s)}{2b(0)(2n+2-s)}\left((\ap-\am)^2-1\right),$$
which is proved like Lemma \ref{lem:positive}.

\medskip\noindent For the proof of compactness in dimension $3$, the following proposition will be useful: its proof is similar to the other ones:

\begin{proposition}[$b\in C^0\setminus C^1$, \textit{The positive case with $\theta=\gamma=0$}]\label{prop:rate:sc:3:3}
Let $\Omega$ be a smooth bounded domain of $\rn$, $n= 3$,  such that $0\in \Omega $ and  assume that  $0< s < 2$,   $\gamma=\theta=0$.  
Let $(\ue)$, $(\he)$, $(\pe)$ and $b$ be as in Proposition \ref{prop:rate:sc:3}. Assume that 
\begin{equation*}
(x,\nabla b(x))=C_2|x|+O(|x|^2)\hbox{ as }x\to 0.
\end{equation*}
Assume that blow-up occurs as in \eqref{blowup}. 
Consider  $\mu_{1,\eps},...,\mu_{N,\eps}$  and $t_{1},...,t_{N}$  from Proposition \ref{prop:exhaust}. Suppose in addition that $\ue>0$ for all $\eps>0$ and that $C_2>0$. Then, $\ap-\am=1$, $u_0\equiv 0$ and 
\begin{align}\label{est:mass:3:bis}
\lim_{\eps\to 0}\frac{p_\eps}{\mu_{N,\eps}^{\ap-\am}}=-\chi_1 \cdot m_{0,h_{0}}(\Omega)-\chi_2\cdot C_2\end{align}
where $\chi_1,\chi_2>0$ are constants and $m_{0,h_{0}}(\Omega)$ is the  mass of $\Omega$ associated with the operator $-\Delta -h_{0}(x)$, defined in Theorem \ref{def:mass}.
\end{proposition}

 \section{Proof of Multiplicity}\label{sec:proof:th}
\noindent {\bf Proof of Theorem \ref{th:cpct:sc}:} We fix $\gamma<(n-2)^2/4$. Let  $h_0\in\ctheta$ be such that $-\Delta - \frac{\gamma}{|x|^2}-h_0(x)$ is coercive and let $b$ satisfy \eqref{hyp:b}. For each  $2<p\leq \crits$, we consider the $C^2$-functional
\begin{equation*}
I_{p, \gamma}(u)=\dfrac{1}{2}\int \limits_\Omega\left(\vert \nabla u\vert^2\, dx-{\gamma} \dfrac{|u|^2}{|x|^2}-h_{0}u^2\right)\, dx-
\frac{1}{p}\int \limits_{\Omega} b(x)\frac {|u|^p}{|x|^s}\, dx
\end{equation*}
on $\huno$, whose critical points are the weak solutions of
\begin{equation}
\label{main}
\left\{ \begin{array}{lll}
-\Delta u-\dfrac{\gamma}{|x|^2}u-h_{0} u&= b(x)\dfrac{|u|^{p-2}u}{|x|^s} &{\rm on} \ \Omega \\
  \hfill   u  &=   0  &{\rm on} \ \partial\Omega.
\end{array} \right.
\end{equation}
For a fixed $u\in \huno$, $u\not\equiv 0$, we have that
$$I_{p, \gamma}(\lambda u)=\dfrac{\lambda^2}{2}\int \limits_\Omega\vert \nabla u\vert^2\, dx-\gamma \dfrac{\lambda^2}{2}\int \limits_{\Omega} \dfrac {|u|^2}{|x|^2}\, dx-\dfrac{\lambda^2}{2}\int \limits_\Omega h_{0}u^2\, dx-\dfrac{\lambda^p}{p}\int \limits_{\Omega} b(x) \frac {|u|^p}{|x|^s}\, dx.$$
Then, by coercivity, we have  $\lim \limits_{\lambda \to \infty}I_{p, \gamma}(\lambda u)=-\infty$, which means that for each finite dimensional subspace $E_k \subset E:=\huno$, there  exists $R_k>0$ such that
\begin{equation}\label{neg}
\sup \{I_{p, \gamma}(u); u\in E_k, \Vert u\Vert_{H_1^2} >R_k\} <0
\end{equation}
as $p\to\crits$. Let $(E_k)^\infty_{k=1}$ be an increasing sequence of subspaces of
$\huno$ such that
$\dim E_k=k $ and $\overline{\cup^\infty_{k=1} E_k}=E:=\huno$ and define
the min-max values:
$$c_{p,k} = {\ds \inf_{g \in {\bf H}_k} \sup_{x\in E_k} I_{p, \gamma}(g(x))},$$
where
$$ {\bf H}_k=\{g \in C(E,E); \hbox{ $g$ is odd and $g(v)=v$
   for $\|v\| > R_k$ for some $R_k>0 \}$}.$$
\begin{proposition} \label{critical.values} With the above notation and assuming  $n\ge 3$, we have:
  \begin{enumerate}
  \item For each $k\in \nn$, $c_{p,k}>0$ and  $\lim\limits_{p\to \crits} c_{p,k}=c_{\crits,k}:=c_k.$
   \item If $2<p<{\crits}$, there exists for each $k$, functions  $u_{p,k}\in
\huno$ such that $I'_{p, \gamma}(u_{p,k})=0$, and $I_{p, \gamma}(u_{p,k})=c_{p,k}$.
   \item  For each  $2<p < {\crits}$, we have  
$ c_{p,k}\ge D_{n,p} k^{\frac{p+1}{p-1}\frac{2}{n}}$ where $D_{n,p}>0$ is 
such that
$\lim \limits_{p\to  {\crits}}D_{n,p}=0$.
\item $\lim\limits_{k\to \infty} c_k=\lim\limits_{k\to \infty} 
c_{\crits,k}=+\infty$.
\end{enumerate}
\end{proposition}
\noindent {\bf Proof:} (1) Coercivity yields the existence of $\Lambda_0>0$ such that
\begin{equation}\label{coer}
\int \limits_\Omega\left(|\nabla u|^2-\frac{\gamma}{|x|^2}u^2-h_{0}u^2\right)\, dx\geq \Lambda_0\int \limits_\Omega|\nabla u|^2\, dx\hbox{ for all }u\in \huno.
\end{equation}
With coercivity, the Hardy and the Hardy-Sobolev inequality, there exists $C>0$ and $\alpha>0$ such that
\[
I_{p, \gamma}(u)\geq \frac{\Lambda_0}{2}\|\nabla u\|_2^2-C\|\nabla u\|_2^p =\|\nabla 
u\|_2^2\left(\frac{\Lambda_0}{2}-C\|\nabla u\|_2^{p-2}\right)\geq \alpha >0
\]
for all $u\in\huno$ such that  $\Vert \nabla u\Vert_2= \rho>0$ is small enough.  
Then the sphere $S_\rho=\{u\in E; \|u\|_{\huno}= \rho\}$ intersects every  image $g(E_k)$  by an odd continuous function $g$. It follows that
\[
c_{p,k}\geq \inf \{I_{p, \gamma}(u); u\in S_\rho \} \geq \alpha >0.
\]
In view of (\ref{neg}), it follows that for each $g\in {\bf H}_k$, we have 
that
\[
\sup\limits_{x\in
E_k}I_{p_i, \gamma}(g(x))=\sup\limits_{x\in D_k}I_{p, \gamma}(g(x))
\]
  where $D_k$ denotes the
ball in $E_k$ of radius $R_k$.
  Consider now a sequence $p_i \to \crits$ and note
first that for each $u\in E$,  we have that $I_{p_i, \gamma}(u) \to I_{\crits, \gamma}
(u)$.   Since $g(D_k)$ is compact and the family of
functionals $(I_{p, \gamma})_p$ is equicontinuous, it follows that
   $\sup\limits_{x\in E_k}I_{p, \gamma}(g(x))\to \sup\limits_{x\in E_k}I_{\crits, \gamma
}(g(x))$, from which follows that
   $\limsup\limits_{i\in \nn}c_{p_i,k}\leq \sup\limits_{x\in E_k}I_{\crits, \gamma}
(g(x))$. Since this holds for any $g\in {\bf H}_k$, it follows that
   \[
    \limsup\limits_{i\in \nn}c_{p_i,k}\leq c_{\crits, k}=c_k.
   \]
On the other hand, the function $f(r)=\frac{1}{p}r^p-\frac{1}{\crits}r^{\crits}$ 
attains its maximum on $[0, +\infty)$ at $r=1$ and therefore
$f(r) \leq \frac{1}{p}-\frac{1}{\crits}$ for all $r>0$. It follows
  \begin{eqnarray*}   I_{\crits, \gamma}(u) &= &I_{p, \gamma}(u)+\int \limits_\Omega \frac{b(x)}{|x|^s}
\left(\frac{1}{p}|u(x)|^p-\frac{1}{\crits}|u(x)|^{\crits }\right)\, dx\\
&\leq&
I_{p, \gamma}(u)+\int \limits_\Omega \frac{b(x)}{|x|^s} \left(\frac{1}{p}-\frac{1}{\crits}\right)\, dx
   \end{eqnarray*}
from which follows that $c_k\leq  \liminf\limits_{i\in \nn}c_{p_i,k}$, and 
claim (1) is proved. \\
If now $p< {\crits}$, we are in the subcritical case, that is we have
compactness in the Sobolev embedding $\huno \to L^p(\Omega; |x|^{-s}dx)$ 
and therefore $I_{p, \gamma}$ satisfies  the Palais-Smale condition. It is then standard to
find critical points $u_{p,k}$ for $I_{p, \gamma}$ at each level $c_{p,k}$  $($see for example the book \cite{Gh}$)$. Consider now the functional 
\begin{equation*}
I_{p, 0}(u)=\frac{1}{2}\int \limits_\Omega\vert \nabla u\vert^2\, dx- 
\frac{1}{p}\int \limits_{\Omega} b(x)\frac {|u|^p}{|x|^s}\, dx
\end{equation*}
and its critical values 
$$c^0_{p,k} = {\ds \inf_{g \in {\bf H}_k} \sup_{x\in E_k} I_{p, 0}(g(x))}.$$
In \cite{gr2} it has been shown that (1), (2) and (3) of Proposition \ref{critical.values} holds, with  
$c^0_{p,k}$ and $c^0_{k}$ replacing $c_{p,k}$ and  $c_{k}$ respectively and $b(x) \equiv 1$.  By similar arguments, $\lim\limits_{k\to \infty} c^0_k=\lim\limits_{k\to \infty} c^0_{\crits,k}=+\infty$.
\noindent On the other hand, with the coercivity \eqref{coer}, we obtain 
\begin{align*}
I_{p, \gamma}(u) \geq \Lambda_0^{\frac{p}{p-2}} I_{p,0}(v) \quad \hbox{for every $u\in \huno$,}
\end{align*}
where  $v=\Lambda_0^{-\frac{1}{p-2}}u$.  It then follows that $\lim\limits_{k\to \infty} c_k=\lim\limits_{k\to \infty} 
c_{\crits,k}=+\infty$.

\smallskip\noindent To complete the proof of Theorem  \ref{th:cpct:sc},   notice that since for each $k$, we
have $$ \lim\limits_{p_i\to \crits}I_{p_i, \gamma}(u_{p_i,k})=\lim\limits_{p_i\to\crits}c_{p_i,k}=c_k,$$
it follows that the sequence  $(u_{p_i,k})_i $ is
uniformly bounded in $\huno$. Moreover, since $I_{p_i}'(u_{p_i,k})=0$, letting $p_i\to \crits$ 
and using the  compactness Theorem \ref{th:cpct:sc} we get a solution $u_k$ of (\ref{main}) in such a way that
$I_{\crit(s), \gamma}(u_k)=\lim\limits_{p\to \crits}I_{p, \gamma}(u_{p,k})=\lim\limits_{p\to\crits}c_{p,k}=c_k$. Since the latter sequence goes to infinity, it follows 
that \eqref{main} has  an infinite number of critical levels with $2<p\leq \crits$. 
\medskip

\section{Proof of Theorems \ref{ground:sc:improv}, \ref{th:cpct:hyp} and  \ref{main.hyp:sc}}
\noindent
We reduce the problem to a smooth bounded domain of the Euclidean space $\rn$.

\medskip\noindent {\bf Case 1: $n\geq 5$.} It follows from Lemma \ref{Lemma:relation between Hyp and Euc solution on bounded domain} that finding solutions to \eqref{main.eq} amounts to finding solutions to
\begin{eqnarray*}\left\{ \begin{array}{lll}
-\Delta v  - \left( \frac{\gamma}{|x|^{2}} + h_{\gamma,\lambda}(x)\right)  v &=b(x) \dfrac{v^{\crits-1}}{|x|^{s}} &\hbox{ in }\Omega \\
\hfill v &=0   & \hbox{on }  \partial \Omega,
\end{array} \right.
\end{eqnarray*}
where $h_{\gamma,\lambda}\in C^{1}(\overline{\Omega})$ with $h_{\gamma,\lambda}(0) = 4\lambda+\dfrac{4(n-2)}{n-4}\left(\gamma-\dfrac{n(n-4)}{4}\right)$. Moreover in this case, $b\in C^2(\overline{\Omega})$, $\nabla b(0)=0$ and
\begin{equation*}
\frac{\Delta b(0)}{b(0)}=\frac{4n(2n-2-s)}{n-4}.
\end{equation*}
Since $b$ is radial, we have that $\partial_{ij}b(0)=\dfrac{\Delta b(0)}{n}\delta_{ij}$. Therefore with $\theta=0$, the limit \eqref{rate:1} in Proposition \ref{prop:rate:sc:2} reads
 \begin{equation}\label{rate:1:bis}
\lim \limits_{\eps \to 0} \frac{\pe}{\mu_{N,\eps}^{2}}  =- \frac{\mathcal{C}_{n,s}}{b(0)} \left(4\lambda+\frac{4(n-2)}{n-4}\left(\gamma-\frac{n(n-4)}{4}\right) +\frac{\Delta b(0)}{n~ \crits }\frac{\int \limits_{\rn}|X|^2\frac{|\tilde{u}_N|^{\crits}}{|X|^s}\, dX}{ \int \limits_{\R^{n}} \tu_{N}^{2} ~dx}\right). 
\end{equation}
when $\gamma<\dfrac{(n-2)^2}{4}-4$. Since $\Delta b(0)>0$, then arguing as in the proof of Theorem \ref{th:cpct:sc}, we get multiplicity of solutions as soon as we have: $4\lambda+\dfrac{4(n-2)}{n-4}\left(\gamma-\dfrac{n(n-4)}{4}\right) \geq 0$. 

\medskip\noindent We now consider the case of positive solutions. We have here $\theta=0$. Using the value of $\Delta b(0)$,  it then follows from Lemma \ref{lem:positive} and Proposition \ref{prop:rate:sc:2} that for $\gamma \geq 0$
 \begin{equation}\label{rate:1:ter}
\lim \limits_{\eps \to 0} \frac{\pe}{\mu_{N,\eps}^{2}}  \leq - 4\lambda \frac{\mathcal{C}_{n,s}}{b(0)} ~~\hbox{ when }\ap-\am>2.
\end{equation}
Therefore, for $\lambda>0$ we get a contradiction and hence compactness for positive subcritical solutions. This yields compactness of minimizers and the existence of minimisers when $\lambda>0$ for  $\ap-\am>2$ and $\gamma \geq 0$. When $\ap-\am=2$ and $\gamma \geq 0$, we get the same result by using \eqref{asymp:crit:log}. When $\ap-\am<2$ and $\gamma \geq 0$, we apply \eqref{est:mass} to get existence of minimizers when the mass $m_{\gamma,\lambda}(\Omega_{\Bn})$ is positive.
\medskip

\noindent {\bf Case 2: $n= 4$.} We prove point (2) of Theorem \ref{main.hyp:sc} for $n=4$. By Lemma \ref{Lemma:relation between Hyp and Euc solution on bounded domain} here $h_{\gamma,\lambda} \in C^{1}\left( \overline{\Omega} \setminus \{0\}\right)$ and   $h_{\gamma,\lambda}(r)= 8\gamma \ln \frac{1}{r}+c_4+O(r)$, as $ r \to 0$. In particular $h_{\gamma,\lambda}(x)=o(|x|^{-\theta})$ as $x\to 0$ for all $\theta>0$.

\smallskip\noindent$\bullet$  Assume that $\gamma<0$, so that $\ap-\am>2$. We fix $\theta\in (0,2)$ small, so that $\ap-\am>2-\theta$ and $K_{h_{\gamma,\lambda}}=0$. Then Propositions \ref{prop:rate:sc:4} and \ref{prop:rate:sc:2:4} do not permit to conclude.

\smallskip\noindent$\bullet$  Assume that $\gamma=0$, so that $\ap-\am=2$ and we have that $h_{\gamma,\lambda}(r)=c_4+O(r)$. Therefore we take $\theta=0$ and $K_{h_{\gamma,\lambda}}=c_4=2(\lambda-2)$. Then point (3) of Proposition \ref{prop:rate:sc:2:4} with $C_1=4(\crits+2)$ yields
\begin{align*}
\lim \limits_{\eps \to 0} \frac{\pe}{\mu_{N,\eps}^{2}\ln\frac{1}{\mu_{N,\eps}}}~=-\left(\frac{\mathcal{C}'_{n,s}}{b(0)} \cdot c_4+  \frac{\mathcal{D}_{n,s}}{b(0)} \cdot 4(\crits+2)\right)\hbox{ when }\ue>0.
\end{align*}
With the definitions of $\mathcal{C}'_{n,s}$ and $\mathcal{D}_{n,s}$, we have that
$$\frac{\mathcal{C}'_{n,s}}{b(0)} \cdot c_4+  \frac{\mathcal{D}_{n,s}}{b(0)} \cdot 4(\crits+2)=\frac{\mathcal{K}^{2} \omega_{3} c_4+\left(\frac{4(\crits+2)}{\crits} \right) \frac{1}{t_{N}^{\frac{4 }{\crits-2}}}  \int \limits_{\R^{4}} |x|^2\frac{|\tu_{N}|^{\crits}}{|x|^{s}} ~dx }{\frac{4-s}{\crits^2} b(0) \sum \limits_{i=1}^{N}  \frac{1}{t_{i}^{ \frac{2}{\crits-2}}} \int \limits_{\rr^4 }   \frac{|\tu_{i}|^{\crits }}{|x|^{s}} ~ dx}$$
with, according to \eqref{def:K},
$$\mathcal{K}:=\frac{L_{0,\Omega}}{t_N^{\frac{\ap}{\crits-2}}}\int \limits_{ \rr^4} \frac{\tilde{u}_{N}(z)^{\crits-1}}{|z|^{s+\am}}\,dz=\frac{L_{0,\Omega}}{t_N^{\frac{2}{\crits-2}}}\int \limits_{ \rr^4} \frac{\tilde{u}_{N}(z)^{\crits-1}}{|z|^{s}}\,dz$$
since $\ap=n-2=2$ and $\am=0$. We first compute $L_{0,\Omega}$. Since for $x,y\to 0$, we have that the Green's function for $-\Delta-h_{0,\lambda}$ behaves like $((n-2)\omega_{n-1})^{-1}|x-y|^{2-n}$ (see Robert \cite{r1}), we get that 
$$L_{0,\Omega}=\frac{1}{(n-2)\omega_{n-1}}=\frac{1}{2\omega_3}.$$
Arguing as in the proof of Lemma \ref{lem:positive} and using the same notations, we get that
$$\int \limits_{\R^{4}} |x|^2\frac{|\tu_{N}|^{\crits}}{|x|^{s}} ~dx =\frac{\mu^{\crits}\nu^2\omega_3}{2-s}I_{p+2}^p\hbox{ with }p=\frac{4}{2-s}$$
Integrating by parts as in Jaber \cite{jaber}, we get that for any $q\geq 0$
$$\left(1-\frac{q}{p+1}\right)\int_0^R\frac{t^{p+1}\, dt}{(1+t)^{q+1}}+\int_0^R\frac{t^{p}\, dt}{(1+t)^{q+1}}=\frac{R^{p+1}}{(p+1)(1+R)^{p+1}}\hbox{ for all }R>0.$$
Taking $q=p+1$ and letting $R\to+\infty$ yields $I_{p+2}^p=\frac{1}{p+1}$. With the same notations and the value of $-\Delta U$ given in the proof of Lemma \ref{lem:positive}, we have that
\begin{eqnarray*}
\int \limits_{ \rr^4} \frac{\tilde{u}_{N}(z)^{\crits-1}}{|z|^{s}}\,dz&=&\mu\nu^{\frac{n-2}{2}}\int_{\rr^4}(-\Delta U)\, dz=\lim_{R\to +\infty}\mu\nu^{\frac{n-2}{2}}\int_{B(0,R)}(-\Delta U)\, dz\\
&=&\lim_{R\to +\infty}\mu\nu^{\frac{n-2}{2}}\int_{B(0,R)}-\partial_\nu U\, d\sigma=\mu\nu^{\frac{n-2}{2}}(n-2)\omega_{n-1}\\
&=&2\mu\nu\omega_{3}
\end{eqnarray*}
Putting these expressions together, we get that
$$\frac{\mathcal{C}'_{n,s}}{b(0)} \cdot c_4+  \frac{\mathcal{D}_{n,s}}{b(0)} \cdot 4(\crits+2)=\frac{4\mu^2\nu^2\omega_3}{t_N^{\frac{4}{2-s}}\frac{4-s}{\crits^2} b(0) \sum \limits_{i=1}^{N}  \frac{1}{t_{i}^{ \frac{2}{\crits-2}}} \int \limits_{\rr^4 }   \frac{|\tu_{i}|^{\crits }}{|x|^{s}} ~ dx}\cdot \lambda$$
and then there exists a constant $\kappa_1>0$ such that
\begin{align*}
\lim \limits_{\eps \to 0} \frac{\pe}{\mu_{N,\eps}^{2}\ln\frac{1}{\mu_{N,\eps}}}~=-\kappa\cdot\lambda.
\end{align*}
Therefore we get compactness when $\lambda>0$.\par
\smallskip\noindent$\bullet$  Assume that $\gamma>0$, so that $\ap-\am<2$. We choose $\theta\in (0,2)$ such that $\ap-\am<2-\theta$. Then point (4) of Proposition \ref{prop:rate:sc:2:4} applies. Compactness for the sequence positive subcritical solutions then follows from the positivity of the mass as in the proof of Theorem \ref{th:cpct:sc:pos}. So in this case we get existence of minimizers when  $\gamma >0$ and the mass $m_{\gamma,\lambda}(\Omega_{\Bn})$ is positive. This proves point (2) of Theorem \ref{main.hyp:sc} for $n=4$.

%
%
%

\noindent {\bf Case 3: $n=3$.} We prove point (2) of Theorem \ref{main.hyp:sc} for $n=3$. In this case, by Lemma \ref{Lemma:relation between Hyp and Euc solution on bounded domain} we have $h_{\gamma,\lambda} \in C^{1}\left( \overline{\Omega} \setminus \{0\}\right)$ and  $h_{\gamma,\lambda}(r)= \dfrac{4 \gamma}{r} +{c_3+O(r)} $ as $ r \to 0$. Then $\theta=1$ and $K_{h_{\gamma,\lambda}}= 4 \gamma$ here. Moreover, we have $(x,\nabla b(x))=C_2|x|+O(|x|^2)$ with $C_2>0$.

\smallskip\noindent$\bullet$  Assume that $\gamma<0$, so that $\ap-\am>1$. Since $\theta=1$, it follows from Propositions \ref{prop:rate:sc:3} and \ref{prop:rate:sc:2:3} and $C_2:=(\crits+2)b(0)>0$ that
\begin{equation*} 
\lim \limits_{\eps \to 0} \frac{\pe}{\mu_{N,\eps}}  =-  \frac{\mathcal{E}_{n,s}}{b(0)}\cdot \left(\frac{16\gamma}{8-s}+\tilde{K}\right)
\end{equation*}
for some $\tilde{K}\geq 0$. Since $\gamma<0$, we cannot conclude.

\smallskip\noindent$\bullet$  Assume that $\gamma=0$, so that $\ap-\am=1$. It follows from Lemma \ref{Lemma:relation between Hyp and Euc solution on bounded domain} that we can take $\theta=0$. It then follows from Proposition \ref{prop:rate:sc:3:3} that
$$\lim_{\eps\to 0}\frac{p_\eps}{\mu_{N,\eps}^{\ap-\am}}=-\chi_1 \cdot m_{0,h_{0}}(\Omega)-\chi_2\cdot C_2$$
Compactness for the sequence positive subcritical solutions then follows from the positivity of the mass.

\smallskip\noindent$\bullet$  Assume that $\gamma>0$, so that $\ap-\am<1$.  Then point (4) of Proposition \ref{prop:rate:sc:2:3} applies. Compactness for the sequence positive subcritical solutions then follows from the positivity of the mass as in the proof of Theorem \ref{th:cpct:sc:pos}. So in this case we get existence of minimizers when  $\gamma >0$ and the mass $m_{\gamma,\lambda}(\Omega_{\Bn})$ is positive. This proves point (2) of Theorem \ref{main.hyp:sc} for $n=3$.

%
%

\smallskip\noindent A quick analysis as above yields no compactness result for sign-changing solutions for $n=3,4$. To be more precise, when estimating the different terms in the Pohozaev identity, we get compactness when a combination of several terms is positive: however, it is not possible to give an apriori criterion to get the positivity.

\section{Proof of the Non-existence result}\label{sec:nonex}  

\medskip\noindent{\bf Proof of Theorem \ref{thm:non:ter}:} We argue by contradiction. Let  $h_0\in \ctheta$ be such that $-\Delta-\gamma|x|^{-2}-h_0$ is coercive  and let $b(x)$ be such that \eqref{hyp:b}  holds.  We fix $\gamma<\dfrac{(n-2)^2}{4}$ and $\Lambda>0$. We assume that there is a family $(\ue)_{\eps>0}\in \huno$ of solutions to 
\begin{equation}\label{eq:non}\left\{ \begin{array}{cl}
-\Delta \ue-\gamma \dfrac{\ue}{|x|^2}-h_\eps \ue= b(x)\dfrac{\ue^{\crits-1}}{|x|^s}  &\text{in } \Omega\setminus\{0\},\\
\ue>0 &\hbox{ in }\Omega\\
\ue=0 &\hbox{ on }\partial\Omega
\end{array}\right.\end{equation}
with $\Vert \nabla \ue\Vert_2\leq\Lambda$, and \eqref{hyp:h0}, \eqref{hyp:he} holds.

\medskip\noindent We claim that $(\ue)_{\eps>0}$ is not pre-compact in $\huno$. Otherwise, up to extraction, there would be $u_0\in \huno$, $u_0\geq 0$, such that $\ue\to u_0$ in $\huno$ as $\eps\to 0$. Passing to the limit in the equation, we get that $u_0\geq 0$ and
\begin{equation}\label{eq:non:2}\left\{ \begin{array}{cl}
-\Delta u_0-\gamma \dfrac{u_0}{|x|^2}-h_0 u_0= b(x)\dfrac{u_0^{\crits-1}}{|x|^s}  &\text{in } \Omega\setminus\{0\},\\
u_0\geq 0 &\hbox{ in }\Omega\\
u_0=0 &\hbox{ on }\partial\Omega.
\end{array}\right.\end{equation}
The coercivity of $-\Delta u_0-\gamma |x|^{-2}-h_0$ and the convergence of $(h_\eps)_\eps$ yields
\begin{eqnarray*}
&&C\left(~\int \limits_{\Omega}\frac{\ue^{\crits}}{|x|^s}dx\right)^{{2}/{\crits}}\leq \int \limits_{\Omega} |\nabla \ue|^2\,dx-\int \limits_\Omega\left(\frac{\gamma}{|x|^2}+\he\right)\ue^2\, dx\leq\int \limits_{\Omega}b(x)\frac{\ue^{\crits}}{|x|^s}dx,
\end{eqnarray*}for some constant $C>0$ and $\epsilon>0$ small, and then, since $\ue>0$, there exists $c_{0}>0$ such that 
$$\int \limits_{\Omega}\frac{\ue^{\crits}}{|x|^s}dx\geq c_0$$
for all $\eps>0$ small. Passing to the limit yields $u_0\not\equiv 0$. Therefore, $u_0>0$ is a solution to \eqref{eq:non} with $\epsilon=0$. 
This is not possible by the hypothesis.

The family $(\ue)_\eps$ is therefore not pre-compact and hence it blows-up with bounded energy. Let $u_0\in \huno$ be its weak limit, which is necessarily a solution to \eqref{eq:non:2}, and hence must be the trivial solution $u_0\equiv 0$.  Proposition \ref{prop:rate:sc:2} then yields that 
either 
\begin{align}\label{res:blowup}
&\hbox{$\ap-\am\geq 2-\theta$ and therefore $K_{h_0}=0$,}\notag\\
\hbox{or } \qquad &\hbox{$\ap-\am<2-\theta$ and therefore  $m_{\gamma,h_0}(\Omega)=0$.}
\end{align}
Since $\gamma> \dfrac{(n-2)^2}{4}-\left(1-\dfrac{\theta}{2}\right)^2$ here, this contradicts our assumption on the mass. Therefore  no such a family of positive solutions $(\ue)_{\eps>0}$ exists, which proves the theorem. \qed  
\medskip

\noindent{\bf Proof of Corollary \ref{thm:non}:} First note that if $h_0$ satisfies 
\begin{equation}\label{hyp:h00}
h_0(x)+\frac{1}{2}(\nabla h_0(x), x)\leq 0\hbox{ for all }x\in \Omega, 
\end{equation}
then by differentiating for any $x\in \Omega$, the function $t\mapsto t^2 h_0(tx)$ (which is well defined for $t\in [0,1]$ since $\Omega$ is starshaped), we get that $h_0\leq 0$. Therefore $-\Delta-\gamma|x|^{-2}-h_0$ is coercive.\\
Assume now there is a positive variational solution $u_0$ of \eqref{eq:non}corresponding to $\epsilon=0$. The Pohozaev identity \eqref{PohoId} then gives
$$\int \limits_{\partial\Omega}(x,\nu)\frac{(\partial_\nu u_0)^2}{2}\, d\sigma-\int \limits_\Omega \left(h_0+\frac{1}{2}(\nabla h_0, x)\right)u_0^2\, dx=0.$$
Hopf's strong comparison principle yields $\partial_\nu u_0 <0$. Since $\Omega$ is starshaped with respect to $0$, we get that $(x,\nu)\geq 0$ on $\partial\Omega$. Therefore, with \eqref{hyp:h00}, we get that $(x,\nu)=0$ for all $x\in \Omega$, which is a contradiction since $\Omega$ is smooth and bounded. Since $\ap-\am<2-\theta$, we use Theorem 7.1 in Ghoussoub-Robert \cite{gr4} to find $\mathcal{K}\in C^2(\overline{\Omega}\setminus\{0\})$ and $A>0$ such that 
$$\left\{\begin{array}{ll}
-\Delta\mathcal{K}-\dfrac{\gamma}{|x|^2}\mathcal{K}-h_0 \mathcal{K}=0 &\hbox{ in }\Omega\\
\mathcal{K}>0 &\hbox{ in }\Omega\\
\mathcal{K}=0 &\hbox{ on }\partial\Omega\setminus\{0\}.
\end{array}\right.$$
and  such that 
$$\mathcal{K}(x)=A\left(\dfrac{\eta(x) }{|x|^{\ap}}+\beta(x)\right)\hbox{ for all }x\in\Omega,$$
where $\eta\in C^\infty_c(\rn)$ and $\beta\in \huno$ are as in Step \ref{step:p11}. We now apply the Pohozaev identity \eqref{PohoId} to $\mathcal{K}$ on the domain $U:=B_\delta(0)$ . Using that $\mathcal{K}^2\in L^1(\Omega)$ and $(\cdot,\nu)(\partial_\nu \mathcal{K})^2\in L^1(\partial\Omega)$ when $\ap-\am<2-\theta$, we get that
$$\int \limits_{\partial\Omega} (x,\nu)\frac{(\partial_\nu \mathcal{K})^2}{2}\, d\sigma-\int \limits_\Omega \left(h_0+\frac{1}{2}(\nabla h_0, x)\right)\mathcal{K}^2\, dx=M_{\delta}$$
where $M_\delta$ is defined in \eqref{est:inf:1}. With \eqref{est:Md}, we then get
$$\int \limits_{\partial\Omega} (x,\nu)\frac{(\partial_\nu \mathcal{K})^2}{2}\, d\sigma-\int \limits_\Omega \left(h_0+\frac{1}{2}(\nabla h_0, x)\right)\mathcal{K}^2\, dx=- 2\omega_{n-1}\left(\frac{(n-2)^2}{4}-\gamma\right)A^2\cdot m_{\gamma,h_0}(\Omega).$$
Since $\Omega$ is star-shaped and $h_0$ satisfies (\ref{hyp:h00}), it follows that $m_{\gamma,h_0}(\Omega)<0$ and Theorem \ref{thm:non:ter} then applies to complete our corollary. 

\section{Appendices}

\subsection{Appendix A: Pohozaev identity}\label{sec:app:poho}
\medskip

\noindent
We let $U\subset\rn$ be a smooth bounded domain in $\rn$. Let $h \in C^{1}(\overline{U})$ and let $b \in C^{1}(\overline{U})$. We suppose $u\in C^2(\overline{U})$ and $p \in [0,\crits-2)$. The classical Pohozaev identity and integration by parts yields for any $y_{0}\in \rn$
\begin{align}\label{PohoId}
&\int \limits_{U}\left((x-y_{0})^i\partial_i u+\frac{n-2}{2}u\right)\left(- \Delta u -\gamma \frac{u}{|x|^2}-hu -b(x)\frac{|u|^{\crits -2-p}}{|x|^{s}}u  \right)~ dx  \\
& - \int \limits_{U}  h(x) u^{2}~ dx  -  \frac{1}{2}\int \limits_{U}    \left( \nabla h, x-y_{0}\right) u^{2} ~ dx - \frac{p}{\crits} \left( \frac{n-s}{\crits -p}\right) \int \limits_{U} b(x)  \frac{|u|^{\crits-p }}{|x|^{s}} ~ dx \notag \\
&  -\frac{1}{\crits -p} \int \limits_{U} \left( x-y_{0}, \nabla b(x) \right) \frac{|u|^{\crits-p }}{|x|^{s}} ~dx - \gamma \int \limits_{U}\frac{(x,y_{0})}{|x|^{4}}u^{2} dx \notag\\
& -\frac{s}{\crits -p} \int \limits_{U}\frac{(x,y_{0})}{|x|^{s+2}}b(x)|u|^{\crits -p} ~dx \notag \\
&=  \int \limits_{\partial U}\left[( x-y_{0}, \nu) \left( \frac{|\nabla u|^2}{2}  -\frac{\gamma}{2}\frac{u^{2}}{|x|^{2}} - \frac{h(x)}{2} u^{2}  - \frac{b(x)}{\crits -p} \frac{|u|^{\crits-p }}{|x|^{s}} \right) \right]d\sigma \notag \\
& \qquad - \int \limits_{\partial U} \left[ \left((x-y_{0})^i\partial_i u+\frac{n-2}{2}u\right)\partial_\nu u \right]d\sigma, \notag
\end{align}
where $\nu$ is the outer normal to the boundary $\partial U$.

\subsection{Appendix B: Regularity}
As to the regularity of the solutions, this will follow from the following result established by Ghoussoub-Robert in \cites{gr3,gr4}.  Assuming that  $\gamma <\frac{(n-2)^2}{4}$, note that the  function $x\mapsto |x|^{-\beta}$ is a solution of  
\begin{equation}\label{homo.Rn}
(-\Delta -\frac{\gamma}{|x|^2})u=0 \qquad \hbox{ on }  \rn\setminus\{0\},
\end{equation}
 if and only if $\beta\in \{\beta_-(\gamma),\beta_+(\gamma)\}$, where 
 \begin{align}\label{def:beta}
 \alpha_{\pm}(\gamma):=\frac{n-2}{2}\pm\sqrt{\frac{(n-2)^2}{4}-\gamma}.
 \end{align}

Actually, one can show that any non-negative solution $u\in C^2(\rn\setminus \{0\})$ of (\ref{homo.Rn}) is of the form

\begin{equation}\label{prop:liouville}
 u(x)=C_- |x|^{-\am}+ C_+ |x|^{-\ap}\hbox{ for all }x\in \rnp,
\end{equation}
where $C_-,C_+\geq 0$. 
\medskip

We collect the following important results from the papers   \cites{gr3,gr4} which we shall use repeatedly in our work. 
\begin{theorem}[Optimal regularity and Hopf Lemma]
\label{th:hopf} 
Let  $\gamma<\frac{(n-2)^2}{4}$ and let $f: \Omega\times \rr\to \rr$ be a Caratheodory function such that
$$|f(x, v)|\leq C|v| \left(1+\frac{|v|^{\crits-2}}{|x|^s}\right)\hbox{\rm  for all }x\in \Omega\hbox{ \rm and }v\in \rr.$$
\begin{enumerate}
\item Let $u\in \huno$ be a weak solution of 
\begin{align}\label{regul:eq}
\Delta u-\frac{\gamma+O(|x|^\tau)}{|x|^2}u=f(x,u),  
\end{align}
for some $\tau>0$. Then, there exists $K\in\rr$ such that
\begin{align}
\label{eq:hopf}
\lim_{x\to 0}\frac{u(x)}{|x|^{-\am}}=K.
\end{align}
Moreover, if $u\geq 0$ and $u\not\equiv 0$, we have that $K>0$.

\item As a consequence, one gets that if $u \in D^{1,2}(\R^{n})$ is a weak solution for
$$-\Delta u-\frac{\gamma}{|x|^2} u= \frac{|u|^{\crits-2} u }{|x|^s} \qquad \hbox{ in } \rnp, $$ 
then there exists $K_{1}, K_{2} \in \R$ such that 
\begin{align*}
u(x)\sim_{|x|\to 0}  \frac{K_{1}}{|x|^{\am}}  \quad \hbox { and } \quad u(x)\sim_{|x|\to +\infty}  \frac{K_{1}}{|x|^{\ap}},  
\end{align*}
and therefore there exists a   constant $C>0$ such that for all $ x$ in $\rnp$, 
\label{prop:bdd:blowupsol}
\begin{align}
|u(x)| \leq  \frac{C}{|x|^{\am}+|x|^{\ap}}.
\end{align}
\end{enumerate}
\end{theorem}
Via the conformal transformation, the above regularity result will yield that the corresponding solutions  for equation \eqref{main.eq} in  the Hyperbolic Sobolev Space $ H_{0}^1(\Omega_{\Bn})$ will satisfy 
 \begin{equation}
\lim \limits_{ |x| \to 0}~ \frac{u(x)} {G(|x|)^{\alpha_{-}}} = K \in \R, 
\end{equation}
where 
 \begin{equation}
 \alpha_{-}(\gamma)= \frac{1}{2}-\sqrt{\frac{1}{4}-\frac{\gamma}{(n-2)^{2}}},
 \end{equation}
which amounts to the regularity claimed in Theorem \ref{main.hyp:sc}. 
 
\subsection{Appendix C: Green's function}
The next theorem  describes the properties of the Green's function of the Hardy-Schr\"odinger  operator in a bounded smooth domain.  To prove this theorem one argue as in the case $\theta=0$ in \cite{gr3}.  

\begin{theorem}[Green's function]\label{th:green:gamma:domain}
Let $\Omega$ be a smooth bounded domain of $\rn$ such that $0\in\Omega$ is an interior point. Let  $\gamma<\frac{(n-2)^2}{4}$, $0\leq \theta <2$ and let  $h_0\in\ctheta$ be such that the operator $-\Delta - \frac{\gamma}{|x|^2}-h(x)$ is coercive in $\Omega$. Then there exists $$G: (\Omega\setminus\{0\})^2\setminus \{(x,x)/\, x\in \Omega\setminus\{0\}\}\to \rr$$  such that 
 
\smallskip\noindent{\bf (i)} For any $p\in\Omega\setminus\{0\}$, $G_p:=G(p,\cdot)\in H_{1}^2(\Omega\setminus B_\delta(p))$ for all $\delta>0$, $G_p\in C^{2,\theta}(\overline{\Omega}\setminus\{0,p\})$

\smallskip\noindent{\bf (ii)} For all $f\in L^p(\Omega)$, $p>n/2$, and all $\varphi\in H_{1,0}^2(\Omega)$ such that
$$  -\Delta\varphi-\left(\frac{\gamma}{|x|^2}+h(x)\right)\varphi=f\hbox{ in }\Omega\; ; \; \varphi_{|\partial\Omega}=0,$$
we have 
\begin{equation}
\varphi(p)=\int \limits_{\Omega}G(p,x)f(x)\, dx
\end{equation}

\medskip\noindent In addition, $G>0$ is unique and 

\smallskip\noindent{\bf (iii)} For all $p\in\Omega\setminus\{0\}$, there exists $C_0(p)>0$ such that
\begin{equation}\label{asymp:G}
G_p(x)\sim_{x\to 0} \frac{C_0(p)}{|x|^{\am}} \hbox{ and }G_p(x)\sim_{x\to p}\frac{1}{(n-2)\omega_{n-1}|x-p|^{n-2}}
\end{equation}
\smallskip\noindent{\bf (iv)}  There exists $C>0$ such that
\begin{equation}\label{est:G:up}
0< G_p(x)\leq C \left(\frac{\max\{|p|,|x|\}}{\min\{|p|,|x|\}}\right)^{\am}|x-p|^{2-n}
\end{equation}
\smallskip\noindent{\bf (v)}  For all $\omega\Subset\Omega$, there exists $C(\omega)>0$ such that
\begin{equation}\label{est:G:low}
C(\omega) \left(\frac{\max\{|p|,|x|\}}{\min\{|p|,|x|\}}\right)^{\am}|x-p|^{2-n}\leq G_p(x)\hbox{ for all }p,x\in\omega\setminus\{0\}.
\end{equation}
\smallskip\noindent{\bf (v)}  There exists $L_{\gamma,\Omega}>0$ such that for any $(h_i) \in\ctheta$  satisfying
$$\lim \limits_{i \to + \infty}\sup_{x\in\Omega}|x|^\theta|h(x)-h_0(x)|+|x|^{\theta+1}|\nabla(h-h_0)(x)|<\epsilon,$$
 and for any sequences of points $(x_i)_i,(y_i)_i\in \Omega$ with
$$y_i=o(|x_i|)\hbox{ and }x_i=o(1)\hbox{ as }i\to +\infty,$$
we have  as $i\to +\infty$
\begin{equation}\label{est:G:4}
G_{h_i}(x_i,y_i)=\frac{L_{\gamma,\Omega}+o(1)}{|x_i|^{\ap} |y_i|^{\am}}.
\end{equation}
\end{theorem}

\end{document}